\documentclass[10pt]{amsart}
\usepackage{files/style}

\pdfoutput=1 

\begin{document}
\title{On tensor products with equivariant commutative operads}
\author{Natalie Stewart}
\date{\today}

\begin{abstract}
  We affirm and generalize a conjecture of Blumberg and Hill:
  unital weak $\mathcal{N}_\infty$-operads are closed under $\infty$-categorical Boardman-Vogt tensor products and the resulting tensor products correspond with \emph{joins of weak indexing systems};
  in particular, we acquire a natural $G$-symmetric monoidal equivalence
  \[
    \underline{\mathrm{CAlg}}^{\otimes}_{I} \underline{\mathrm{CAlg}}^{\otimes}_{J} \mathcal{C} \simeq \underline{\mathrm{CAlg}}^{\otimes}_{I \vee J} \mathcal{C}.
  \]
  We accomplish this by showing that $\mathcal{N}_{I\infty}^{\otimes}$ is $\obv$-idempotent and $\mathcal{O}^{\otimes}$ is local for the corresponding smashing localization if and only if $\mathcal{O}$-monoid $G$-spaces satisfy $I$-indexed Wirthm\"uller isomorphisms.

  Ultimately, we accomplish this by advancing the equivariant higher algebra of \emph{cartesian and cocartesian $I$-symmetric monoidal $\infty$-categories}. 
  Additionally, we acquire a number of structural results concerning $G$-operads, including a canonical lift of $\obv$ to a presentably symmetric monoidal structure and a general disintegration \& assembly procedure for computing tensor products of non-reduced unital $G$-operads.
  All such results are proved in the generality of atomic orbital $\infty$-categories.
  
  We also achieve the expected corollaries for (iterated) Real topological Hochschild and cyclic homology and construct a natural $I$-symmetric monoidal structure on right modules over an $\mathcal{N}_{I\infty}$-algebra.
 \end{abstract}

\maketitle
 
\toc
\newpage
\section*{Introduction}  
We're concerned with the relationship between homotopy-coherent interchange and equivariant commutative algebras, incarnated via \emph{$\cN_{I\infty}$-algebras} (henceforth \emph{$I$-commutative algebras}) in the sense of \cite{Blumberg-op,EBV}.
In particular, in \cite{EBV} we constructed a natural ``pointwise'' $G$-symmetric monoidal structure $\uAlg_{\cO}^{\otimes}(\cC)$ on $\Alg_{\cO}(\cC)$.
We hope to answer the following questions, where $\uCAlg_I^{\otimes}(\cC) \deq \uAlg_{\cN_{I\infty}}^{\otimes}(\cC)$.
\begin{questions*}
  Let $\cO^{\otimes}$ be a unital $G$-operad and $I,J \subset \FF_G$ a pair of unital weak indexing categories.
  \begin{enumerate}[label={(\Roman*)}]
    \item \label[ques]{q: absorption to I} When is the forgetful natural transformation $\Alg_{\cO} \uCAlg_I^{\otimes}(-) \implies \CAlg_I(-)$ an equivalence?
    \item \label[ques]{q: absorption to O} When is the forgetful natural transformation $\CAlg_I \uAlg_{\cO}^{\otimes}(-) \implies \Alg_{\cO}(-)$ an equivalence?
    \item \label[ques]{q: BH conjecture} What is the (unique) $G$-operad $\cO^{\otimes}$ with natural equivalence $\CAlg_I \uCAlg_J^{\otimes}(-) \simeq \Alg_{\cO}(-)$?\qedhere
  \end{enumerate}
\end{questions*}
Each of the left hand sides of these proposed equivalences are corepresented by \emph{Boardman-Vogt tensor products} of $G$-operads, so these \cref{q: absorption to I,q: absorption to O} are equivalent to the question of when distinguished maps $\cN_{I\infty}^{\otimes} \rightarrow \cO^{\otimes} \obv \cN_{I\infty}^{\otimes}$ and $\cO^{\otimes} \rightarrow \cN_{I\infty}^{\otimes} \obv \cO^{\otimes}$ are equivalences;
moreover \cref{q: BH conjecture} asks the value of the tensor product $\cN_{I\infty}^{\otimes} \obv \cN_{J\infty}^{\otimes}$.
In this form, \Blumberg{Conj.}{6.27} conjectured an answer.
\begin{conjecture*}[Blumberg-Hill]
  If $I$ and $J$ are \emph{indexing categories} then $\cN_{I\infty}^{\otimes} \obv \cN_{J\infty}^{\otimes} \simeq \cN_{I \vee J \infty}^{\otimes}$.
\end{conjecture*}

We begin by completely characterizing $G$-operad algebras in (co)cartesian $I$-symmetric monoidal $\infty$-categories:
cocartesian $I$-symmetric monoidal structures are characterized by the property that their $G$-objects canonically lift to $\cO$-algebras for any reduced $I$-operad $\cO^{\otimes}$, and cartesian $I$-symmetric monoidal structures are characterized by an \emph{$\cO$-monoid} formula generalizing \ha{Prop.}{2.4.2.5}.
Using this, we show that \cref{q: absorption to I} is true precisely when $\cO^{\otimes}$ is \emph{reduced} and $I$-commutative algebras admit underlying $\cO$-algebra structures.

We conclude that the unique map $\EE_{0}^{\otimes} \rightarrow \cN_{I \infty}^{\otimes}$ witnesses $\cN_{I \infty}^{\otimes}$ as an idempotent algebra in $\Op^{\uni}_G$;
\cref{q: absorption to O} asks to classify the associated smashing localization.
Indeed, we show that the equivalence holds whenever $\cO$-algebra $G$-spaces satisfy $I$-indexed Wirthm\"uller isomorphisms. 

Since $\uAlg_{\cO}(\cS_G)$-ambidextrous arities form a weak indexing category, we find that the intersection of the $\cN_{I\infty}^{\otimes}$- and $\cN_{J\infty}^{\otimes}$-smashing local categories is the $\cN_{I \vee J \infty}^{\otimes}$-smashing local category, constructing an equivalence $\cN_{I \infty}^{\otimes} \obv \cN_{J \infty}^{\otimes} \simeq \cN_{I \vee J \infty}^{\otimes}$ in full generality.
This answers \cref{q: BH conjecture} by affirming the evident unital extension of Blumberg-Hill's conjecture, constructing a (unique) natural equivalence
\[
  \CAlg_I \uCAlg_J^{\otimes}(\cC) \simeq \CAlg_{I \vee J}(\cC).
\]
This is the third part of an ongoing project  to develop the parameterized and equivariant higher algebra predicted in \cite{Barwick_Intro,Nardin} into simply usable foundations for equivariant homotopy theory and $K$-theory \cite{Windex,EBV};
as such, we spend the last third of the paper fleshing out higher algebraic corollaries.

These corollaries fall into two classes:
the first class gives $\Comm_G^{\otimes} \in \Op_G$ a unique idempotent algebra structure, which determines a unique compatible idempotent algebra structure on its $G$-symmetric monoidal envelope, enabling \emph{symmetric monoidality} of the equivariant equifibered perspective of \cite{Barkan,Haugseng_SMC,Barkan_equifibered}.
From this, we lift $\Op_G$ with the Boardman-Vogt tensor product to a canonical presentably symmetric monoidal $G$-$\infty$-category;
as an application, we develop equivariant operadic disintegration and assembly, and the associated distributivity of $\obv$ allows us to compute tensor products of unital $G$-operads whose underlying $G$-$\infty$-categories are $G$-spaces in terms of tensor products of reduced $G$-operads.

The second class simply applies \cref{q: absorption to I,q: absorption to O}:
by answering \cref{q: absorption to O} for $\cN_{I \infty}^{\otimes} \simeq \EE_\infty^{\otimes}$, we get an $\cO$-symmetric monoidal structure on left modules over an $\cO$-algebra;
for instance, specializing to $\cO^{\otimes} = \cN_{J \infty}^{\otimes}$ confirms a hypothesis of Hill 
\cite[Rmk.~\href{https://arxiv.org/pdf/1708.03017v1\#theorem.3.15}{3.15}]{Hill_chromatic}.

By answering \cref{q: absorption to I} for $\cO^{\otimes} \simeq \EE^{\otimes}_V$, we acquire an $I$-commutative algebra structure on (lax) $I$-symmetric monoidal $\EE_V$-algebra invariants of $I$-commutative algebras;
for instance, this constructs an $I$-commutative algebra structure on Real topological Hochschild homology and Real topological cyclic homology of an $I$-commutative algebra whenever $I$-commutative algebras have underlying $\EE_\sigma$-algebras.

We now move to a more careful account of the background, motivation, and main results of this paper.
\stoptocwriting 

\subsection*{Background and motivation}
Let $\cC$ be a 1-category with finite products.
Recall that a \emph{commutative monoid in $\cC$} is the data 
\[
  A \in \Ob(\cC), \hspace{30pt} \text{multiplication} \;\; \mu\cln A \times A \rightarrow A, \;\; \text{and} \hspace{30pt} \text{unit} \;\; \eta\cln * \rightarrow A,
\]
subject to the usual unitality, associativity, and commutativity assumptions;
more generally, if $(\cC,\otimes,1)$ is a symmetric monoidal 1-category, a \emph{commutative algebra in $\cC$} is the data of
\[
  R \in \Ob(\cC), \hspace{30pt} \text{multiplication} \;\; \mu\cln R \otimes R \rightarrow R, \;\; \text{and} \hspace{30pt} \text{unit} \;\; \eta\cln 1 \rightarrow R,
\]
satisfying analogous conditions.
When $\cC = \Set$, this recovers the traditional theory of commutative monoids, and when $\cC = \Mod_k$ with the tensor product of $k$-modules, this recovers the traditional theory of commutative $k$-algebras.
These have been the subject of a great deal of homotopy theory in three guises:
\begin{enumerate}[label={(\roman*)}]
  \item \label[perspective]{Lawvere perspective} We may define the $(2,1)$-category $\Span(\FF)$ to have objects the finite sets, morphisms from $X$ to $Y$ the spans of finite sets $X \leftarrow R \rightarrow Y$, 2-cells the isomorphisms of spans
    \[\begin{tikzcd}[row sep=0em]
	& R \\
  X && Y, \\
	& {R'}
	\arrow[from=1-2, to=2-1]
	\arrow[from=1-2, to=2-3]
	\arrow["\sim"{description}, from=1-2, to=3-2]
	\arrow[from=3-2, to=2-1]
	\arrow[from=3-2, to=2-3]
\end{tikzcd}\] 
    and composition the pullback of spans
    \[\begin{tikzcd}[row sep=0em]
	&& {R_{XZ}} \\
	& {R_{XY}} && {R_{YZ}} \\
  X && Y && Z.
	\arrow[from=1-3, to=2-2]
	\arrow[from=1-3, to=2-4]
	\arrow["\lrcorner"{anchor=center, pos=0.125, rotate=-45}, draw=none, from=1-3, to=3-3]
	\arrow[from=2-2, to=3-1]
	\arrow[from=2-2, to=3-3]
	\arrow[from=2-4, to=3-3]
	\arrow[from=2-4, to=3-5]
\end{tikzcd}\]
  If $\cC$ is an $\infty$-category, then we define the \emph{$\infty$-category of commutative monoids in $\cC$} as the $\cC$-valued models of the associated Lawvere theory;
  that is, we define the product-preserving functor category
  \[
    \CMon(\cC) \deq \Fun^{\times}(\Span(\FF), \cC),
  \]
  noting that products in $\Span(\FF)$ correspond with disjoint unions of finite sets.
  Indeed, if $\cC$ is a 1-category and $A$ a commutative monoid in $\cC$, we flesh this out with the dictionary
  \def\mapstowithspace{\hspace{30pt} \longmapsto \hspace{30pt}}
  \begin{align*}
    ([2] = [2] \rightarrow [1])             &\mapstowithspace \, \mu\cln A^{\times 2} \rightarrow A;\\
    (\emptyset = \emptyset \rightarrow [1]) &\mapstowithspace \, \eta\cln * \simeq A^{\times 0} \rightarrow A;\\
    ([1] \leftarrow [2] = [2])              &\mapstowithspace \Delta\cln A \rightarrow A^{\times 2}\\
    ([1] \leftarrow \emptyset = \emptyset)  &\mapstowithspace \;\;\; !\cln A \rightarrow A^{\times 0} \simeq *.
  \end{align*}
  Unitality, associativity, and commutativity are conveniently packaged by functoriality.
  This turns out to be equivalent to Graeme Segal's \emph{special $\Gamma$ spaces} \cite{Segal} when $\cC = \cS$, and for general $\cC$, it recovers the anologously defined theory in $\cC$ (see \cite[Thm.~\href{https://arxiv.org/pdf/1109.1598v1\#thm.5.7}{5.7}]{Cranch}).
  \item \label[perspective]{Semiadditive perspective} We say that a pointed $\infty$-category is \emph{semiadditive} if it has finite products and coproducts and for all finite sets $S$, the ``identity matrix'' natural transformation $\coprod_{s \in S}(-) \implies \prod_{s \in S}(-)$ is an equivalence.
    The full subcategory $\Pr^{L,\oplus} \subset \Pr^L$ of \emph{semiadditive presentable $\infty$-categories} possesses a localization functor $L_{\oplus}\cln \Pr^L \rightarrow \Pr^{L,\oplus}$, which we study.
  \item \label[perspective]{Operad perspective} Let $\Op$ denote the $\infty$-category of operads.\footnote{This is unambiguous \cite{Hinich}, but we will tend to model these as $\infty$-operads in the sense of \cite{HA}.}
    Then, there is a terminal operad $\Comm^{\otimes} \simeq \EE_\infty^{\otimes}$;
    given $\cC$ a symmetric monoidal $\infty$-category, we may form the $\infty$-category of \emph{commutative algebra objects}
    \[
      \CAlg(\cC) \deq \Alg_{\Comm}(\cC) \simeq \Alg_{\EE_\infty}(\cC).
    \]
    We study this and its specialization to the cartesian symmetric monoidal structure.
\end{enumerate}
These three perspectives each present the same $\infty$-category, i.e. \cite{Cranch,Gepner} show that
\begin{equation}\label{Forgetful functors}
  \CMon(\cC) \simeq \CAlg(\cC^\times) \simeq L_{\oplus} \cC.
\end{equation}
As a result, translating between these perspectives has proved invaluable;
for instance, \cite{Gepner} uses \cref{Semiadditive perspective,Operad perspective} to construct an essentially unique symmetric monoidal structure on $\CMon(\cC)$ and \cite{Cnossen_bispans} uses \cref{Lawvere perspective,Operad perspective} to model commutative algebras in $\CMon(\cC)^{\otimes}$ as models for the Lawvere theory of \emph{commutative semirings}.

Crucially, \cref{Lawvere perspective,Operad perspective} may be used to construct homotopical lifts of the \emph{Eckmann-Hilton argument};
for instance, in \cite{Schlank}, it is shown that for \emph{any} reduced operad $\cO^{\otimes}$, the forgetful functors
\[
  \CAlg \Alg_{\cO}^{\otimes}(\cC) \rightarrow \CAlg(\cC) \leftarrow \Alg_{\cO} \CAlg^{\otimes}(\cC),
\]
are equivalences for the ``pointwise'' symmetric monoidal structure on algebras.
Such an equivalence may be exhibited by recognizing the far left and far right side each as algebras over the \emph{Boardman-Vogt tensor product} $\cO^{\otimes} \obv \Comm^{\otimes}$ and each arrow as pullback along the canonical map
\[
  \Comm^{\otimes} \simeq \triv^{\otimes} \obv \Comm^{\otimes} \xrightarrow{\mathrm{can} \otimes \id} \cO^{\otimes} \obv \Comm^{\otimes};
\]
that \cref{Forgetful functors} consists of equivalences reduces to the well-known fact that $\cO^{\otimes} \obv \Comm^{\otimes} \in \Op$ is terminal, which one can quickly prove via \cref{Semiadditive perspective,Operad perspective}.

This result is used ubiquitously to replace (lax) symmetric monoidal functors $\Alg_{\cO}^{\otimes}(\cC) \rightarrow \cC^{\otimes}$ with (lax) symmetric monodial endofunctors 
\[
  \CAlg^{\otimes}(\cC) \simeq \CAlg^{\otimes}\Alg_{\cO}^{\otimes}(\cC) \rightarrow \CAlg^{\otimes}(\cC);
\]
for instance, this underlies the symmetric monoidal structure on left-modules \cite{HA} and the multiplicative structure on factorization homology \ha{Thm.}{5.5.3.2}, $\mathrm{TC}$ \cite[\S~\href{https://arxiv.org/pdf/1707.01799v2\#section.4.2}{4.2}]{Nikolaus}, and algebraic $K$-theory \cite{Barwick_K}.

This paper concerns the analogs of \cref{Lawvere perspective,Semiadditive perspective,Operad perspective} in the equivariant theory of algebra stemming from Hill-Hopkins-Ravanel's use of \emph{norms of $G$-spectra} on the Kervarire invariant one problem, as well as the resulting theory of \emph{indexed tensor products and (co)products} (see \cite{Hill_SMC,Barwick_Intro,Nardin}).

For the rest of this introduction, fix $G$ a finite group.
In $G$-equivariant homotopy theory, the point is replaced with elements of the \emph{orbit category} $\cO_G \subset \Set_G$, whose objects are homogeneous $G$-sets $[G/H]$;
indeed, Elmendorf's theorem \cite{Elmendorf} realizes $G$-spaces as coefficient systems $\cS_G \simeq \Fun\prn{\cO_G^{\op}, \cS}$.\footnote{
  Maps $[G/K] \rightarrow [G/H]$ may equivalently be presented as elements of $g$ such that $gKg^{-1} \subset H$, modulo $K$;
see e.g. \cite{Dieck} for details.}
In $G$-equivariant higher category theory, $\infty$-categories are thus replaced with $G$-$\infty$-categories
\[
  \Cat_G \deq \Fun\prn{\cO_G^{\op}, \Cat}.
\]
In $G$-equivariant higher algebra, following \cref{Lawvere perspective}, we may form the effective Burnside 2-category $\Span(\FF_G)$ whose objects are finite $G$-sets, whose morphisms are spans, whose 2-cells are isomorphisms of spans, and whose composition is pullback;
the following central definition is the heart of this subject.
\begin{definition*}
  The \emph{$\infty$-category of $G$-commutative monoids in $\cC$} is the product-preserving functor $\infty$-category
  \[
      \CMon_G(\cC) \deq \Fun^\times(\Span(\FF_G),\cC);
  \]
  the \emph{$\infty$-category of small $G$-symmetric monoidal $\infty$-categories} is
  \[
    \Cat_G^{\otimes} \deq \CMon_G(\Cat).\qedhere
  \]  
\end{definition*}
These are a homotopical lift of Dress' \emph{semi-Mackey functors} \cite{Dress} (c.f. \cite{Lindner}).
Indeed, given $\cC^{\otimes} \in \Cat^{\otimes}_G$ a $G$-symmetric monoidal $\infty$-category,
pullback along the product-preserving functor 
\[
    \iota_H\colon \Span(\FF) \xrightarrow{* \mapsto G/H} \Span(\FF_G)
\]
constructs a symmetric monoidal $\infty$-category $\cC_H^{\otimes} \deq \iota_H^* \cC^{\otimes}$ whose underlying $\infty$-category $\cC_H$ is the value of $\cC^{\otimes}$ on the orbit $[G/H]$.
For all subgroups $K \subset H \subset G$, the covariant and contravariant functoriality of $\cC^{\otimes}$ then yield symmetric monoidal \emph{restriction} and \emph{norm} functors
\begin{align*}
    \Res_K^H\colon \cC^{\otimes}_H &\rightarrow \cC^{\otimes}_K,\\
    N_K^H\colon \cC^{\otimes}_K &\rightarrow \cC^{\otimes}_H,
\end{align*}
which satisfy a form of Mackey's \emph{double coset formula}.
\begin{example*}[{\cite{Cnossen_tambara,Bachmann}}]
  There is a presentably $G$-symmetric monoidal $\infty$-category $\uSp_G^{\otimes}$ with:
  \begin{itemize}
    \item $H$-value given by the symmetric monoidal $\infty$-category $\prn{\uSp^{\otimes}_G}_H \simeq \Sp_H^{\otimes}$ of \emph{genuine $H$-spectra},
    \item restriction functors $\Res_K^H\colon \Sp^{\otimes}_H \rightarrow \Sp^{\otimes}_K$ given by the usual restriction functors,  and
    \item norm functors $N_K^H\colon \Sp_K^{\otimes} \rightarrow \Sp_H^{\otimes}$ given by the \emph{HHR norm} of \cite{HHR}.
  \end{itemize}
  In fact, this structure is completely determined by its unit object $\SS_G \in \Sp_G^{\otimes}$.
\end{example*}
Fix $\cC^{\otimes} \in \Cat_G^{\otimes}$.
If $H \subset G$ is a subgroup and $S \in \FF_H$ a finite $H$-set, we may form the induced $G$-set $\Ind_H^G S \rightarrow [G/H]$, and the covariant and contravariant functoriality then yield an \emph{$S$-indexed tensor product} and \emph{$S$-indexed diagonal}
\[
  \bigotimes^S_K\cln \cC_S \rightarrow \cC_H, \hspace{60pt} \Delta^S\cln \cC_H \rightarrow \cC_S.
\]
where $\cC_S \deq \hspace{-10pt} \prod\limits_{[H/K] \in \Orb(S)} \hspace{-10pt} \cC_K$.
Note that $N_H^K$ is the $[H/K]$-indexed tensor product and $\Res_K^H$ the $[H/K]$-indexed diagonal.
As explained in \cite{EBV}, the ``orbit collapse'' factorization $\Ind_H^G S \rightarrow \coprod_{[H/K] \in \Orb(S)} [G/H] \rightarrow [G/H]$ yields natural equivalences
\[
  \bigotimes^S_K X_K \simeq \bigotimes_{[H/K] \in \Orb(S)} N_K^H X_K, \hspace{60pt} \Delta^S (X) = \prn{\Res_K^H X}_{[H/K] \in \Orb(S)},
\]
so we may often reduce arguments about $S$-indexed tensor products to to binary tensor products and norms.
Similarly, we define the \emph{$S$-indexed tensor power}
\[
  X^{\otimes S}_H \deq \bigotimes^S_K \prn{\Delta^S X_H}  
  \simeq \bigotimes^S_K \Res_H^K X_H 
  \simeq \bigotimes_{[H/K] \in \Orb(S)} N_K^H \Res_K^H X_H.
\]
If it exists, the pointwise left-adjoint to $\Delta^S$ is the \emph{indexed coproduct}
\[
  \coprod_K^S X_K \simeq \coprod_{[H/K] \in \Orb(S)} \Ind_K^H S,
\]
where $\Ind_K^H$ is the left adjoint to the restriction map $\cC_H \rightarrow \cC_K$.
The \emph{indexed products} are defined analogously.

Given $H \subset G$ a subgroup, we say that $\cC$ is \emph{$H$-pointed} if $\cC_K$ is pointed for all $(K) \subset (H)$.
Given $S \in \FF_H$, we say that $S$ is \emph{$\cC$-ambidextrous} if $\cC$ is $H$-pointed, $\cC$ admits $S$-indexed products and coproducts, and the \emph{Wirthm\"uller} natural transformation
\[
  W_{S}\colon \coprod_K^S(-) \implies \prod_K^S(-)
\]
(often called the \emph{norm}, first defined in \nar{Def.}{5.3})
is an equivalence.
We say that $\cC$ is \emph{$G$-semiadditive} if $S$ is $\cC$-ambidextrous for all $S \in \FF_H$ and $H \subset G$.
More generally, if $\uFF_I \subset \uFF_G$ is a weak indexing system corresponding with the weak indexing category $I \subset \FF_G$ (see \cite{Windex} or our review in \cref{Recollections on T-operads subsection}), we say that $\cC$ is \emph{$I$-semiadditive} if $S$ is $\cC$-ambidextrous whenever $S \in \FF_{I,H}$.

In this level of generality, \cref{Lawvere perspective,Semiadditive perspective} are known to present equivalent $\infty$-categories of $I$-commutative monoids;
indeed, the \emph{semiadditive closure} theorem of \cite[Thm.~\href{https://arxiv.org/pdf/2403.07676v3\#introthm.2}{B}]{Cnossen_semiadditive} demonstrates that $\Pr_G^{L,I-\oplus} \subset \Pr_G^{L}$ is a smashing localization implemented by
\[
  L_{I-\oplus}(\cC) \simeq \uCMon_I(\cC) \deq \uFun_G^{\times}\prn{\Span_I(\uFF_G), \cC}, 
\]  
and in particular, when $\cC$ is a $G$-$\infty$-category of coefficient systems
\[
  \uCoFr^G(\cD)_H \deq \Fun\prn{\cO_H^{\op},\cD},
\]
\cite[Thm.~\href{https://arxiv.org/pdf/2403.07676v3\#introthm.3}{C}]{Cnossen_semiadditive} yields the formula
\[
  \uCMon_I\prn{\uCoFr^G(\cD)}_H \simeq \Fun^\times\prn{\Span_I(\FF_H),\cD},
\]
where $\Span_I(\FF_H) \subset \Span(\FF_H)$ is the wide subcategory of spans whose forward maps lie in the restriction of $I$ to $\FF_H$.
Thus, we set the notation $\CMon_I(\cD) \deq \uCMon_I\prn{\uCoFr^G(\cD)}_G \simeq \Fun^\times\prn{\Span_I(\FF_G),\cD}$ and make the following definition.
\begin{definition*}
  For $I$ is a weak indexing category, the \emph{$\infty$-category of small $I$-symmetric monoidal $\infty$-categories} is
  \[
    \Cat_I^{\otimes} \deq \Fun^{\times}\prn{\Span_I(\FF_G), \Cat}.\qedhere
  \]
\end{definition*}

Following through on \cref{Operad perspective}, algebraic objects $X_\bullet$ in a $G$-symmetric monoidal $\infty$-category should possess collections of \emph{$S$-ary operations} $X_H^{\otimes S} \rightarrow X_H$ subject to various coherences, controlled by a theory of \emph{genuine equivariant operads};
we use Nardin-Shah's $\infty$-category $\Op_G$, whose objects we call \emph{$G$-operads}.
Given $\cO^{\otimes} \in \Op_G$ a $G$-operad, $K \subset H \subset G$ a pair of subgroups, $S \in \FF_H$ a finite $H$-set, and $T_i$ a finite $K_i$-set for all orbits $[H/K_i] \subset S$, in \cite{EBV} we constructed a \emph{space of $S$-ary operations} $\cO(S)$, \emph{operadic composition maps}
\begin{equation}\label{Operadic composition equation}
  \gamma\cln \cO(S) \otimes \hspace{-5pt} \bigotimes_{[H/K_i] \in \Orb(S)} \hspace{-5pt} \cO(T_i) \rightarrow \cO\prn{\coprod_{[H/K_i] \in \Orb(S)} \Ind_{K_i}^H T_i},
\end{equation}
\emph{operadic restriction maps}
\begin{equation}\label{Operadic restriction equation}
  \Res\cln \cO(S) \rightarrow \cO\prn{\Res_K^H S},
\end{equation}
and \emph{equivariant symmetric group action}
\begin{equation}\label{Operadic symmetric action equation}
    \rho\cln \Aut_H(S) \times \cO(S) \rightarrow \cO(S).
\end{equation}
More generally, given a map $T \rightarrow S$ and objects $(\bB,\bC) \in \cO(T;S)$, we constructed a space of multi-operations $\cO(\bB;\bC)$.
We made the following simplifying definition.
\begin{definition*}
  A $\cO^{\otimes} \in \Op_G$ \emph{has one color} if $\cO(*_H) = *$ for all $H \subset G$;
  these span a full subcategory $\Op^{\oc}_G \subset \Op_G$.
\end{definition*}
We showed in \cite[Thm.~\href{https://arxiv.org/pdf/2501.02129v1\#cooltheorem.1}{A}]{EBV} that \cref{Operadic restriction equation,Operadic symmetric action equation} lift to a conservative functor $\Op_G^{\oc} \rightarrow \Fun(\tot \Sigma_G, \cS)$;
in other words, $\prn{S \mapsto \cO(S) \;\; \middle| \;\; H \subset G, \,S \in \FF_H}$ are jointly conservative.

When $\cO^{\otimes}$ has one color, an $\cO$-algebra in a $G$-symmetric monoidal $\infty$-category $\cC^{\otimes}$ can intuitively be viewed as a tuple $\prn{X_H \in \cC^{BW_G(H)}_H}_{G/H \in \cO_G}$ with $X_K \simeq \Res_K^H X_H$ for all $K \subset H \subset G$, together with $\cO(S)$-actions
\begin{equation}\label{Mu action equation}
  \mu_S\colon \cO(S) \rightarrow \Map_{\cC_H}\prn{X_H^{\otimes S}, X_H}
\end{equation}
for all $H \subset G$ and $S \in \FF_H$, homotopy-coherently compatible with \cref{Operadic composition equation,Operadic restriction equation,Operadic symmetric action equation}.\footnote{Here, $W_G(H) = N_G(H)/H$ is the \emph{Weyl group} of $H \subset G$, i.e. the automorphism group of the homogeneous $G$-set $[G/H]$. The restriction-compatible data specified above may be more familiarly referenced as a \emph{$G$-object};
it's canonically extended from a choice $X_G \in \cC_G$.}
We are concerned with the following examples.

\begin{example*}
  There exists a terminal $G$-operad $\Comm_G^{\otimes}$, which is characterized up to (unique) equivalence by the property that $\Comm_G(S)$ is contractible for all $S \in \FF_H$;
  its algebras are endowed with contractible spaces of maps $X_H^{\otimes S} \rightarrow X_H$ for all $S \in \FF_H$, as well as coherent homotopies witnessing their compatibility.
  We call these \emph{$G$-commutative algberas}.

  On one hand, we saw in \cite[\S~\href{https://arxiv.org/pdf/2501.02129v1\#subsection.2.7}{2.7}]{EBV} that $\Comm_G$-algebras present a homotopical lift of Hill-Hopkins' \emph{$G$-commutative monoids} \cite[\S~\href{https://arxiv.org/pdf/1610.03114v1\#subsection.3.3}{3.3}]{Hill_SMC}, though we prefer to reserve this name for the Cartesian case, following the convention of \cite{HA}.
  On the other hand, our model agrees with that of \cite{Cnossen_tambara}, so the recent \emph{homotopical Tambara functor theorem} of Cnossen, Lenz, and Linskens \cite[Thm.~\href{https://arxiv.org/pdf/2407.08399v1\#introthm.2}{B}]{Cnossen_tambara} presents $G$-commutative algebra objects in $\uSp_G^{\otimes}$ (i.e. \emph{$G$-commutative ring spectra}) as a form of \emph{homotopical $G$-Tambara functors}. 

  Additionally, the recent rectification theorem of Lenz, Linskens, and P\"utzst\"uck \cite{Lenz_norms} establishes $G$-commutative ring spectra as a Dwyer-Kan localization of strict commutative algberas in symmetric (or orthogonal) $G$-spectra at the weak equivalences transferred from a ``positive stable'' model structure.
\end{example*}

\begin{example*}
  Let $V$ be a real orthogonal $G$-representation.
  There is a \emph{little $V$-disks $G$-operad} $\EE^{\otimes}_V$ whose structure spaces are \emph{spaces of equivariant configurations}:
  \[
    \EE_V(S) \simeq \Conf_S^H(V)
  \]
  (see \cite{Hill_disk,Horev}).
  This is modelled by the \emph{Steiner graph $G$-operad}, so e.g. pointed $G$-spaces of the form $X = \Omega^V Y \deq \Map_*\prn{S^V,Y}$ lift to $\EE_V$-spaces by composition of loops \cite{Guillou,Horev_poincare}; 
  moreover, many $\EE_V$-ring spectra may be constructed as Thom $G$-spectra of $V$-fold loop maps \cite{Horev_poincare}. 
\end{example*}

\begin{example*}
  Given $I \subset \FF_G$ a weak indexing category, in \cite{EBV} we constructed a \emph{weak $\cN_\infty$ $G$-operad} $\cN_{I \infty}^{\otimes}$ which is characterized up to (unique) equivalence by its structure spaces
  \begin{equation}\label{Ninfty structure spaces equation}
    \cN_{I \infty}(S) \simeq \begin{cases}
      * & S \in \uFF_I \\ 
      \emptyset & S \not \in \uFF_I 
    \end{cases}	
  \end{equation}
  These recover the $\cN_\infty$-operads of \cite{Blumberg-op} when $I$ is an indexing category, i.e. $\cN_{I\infty}(n \cdot *_G) \simeq *$ for $n \in \NN$;
  in general, they are identified as the sub-terminal objects of $\Op_G$ \cite[Thm.~\href{https://arxiv.org/pdf/2501.02129v1\#cooltheorem.3}{C}]{EBV}.
\end{example*}
For instance, we verify in \cref{EV is ninfty} that the condition $V \oplus V \simeq V$ for an orthogonal $G$-representation $V$ implies that $\EE_V$ is a weak $\cN_\infty$-operad, which is an $\cN_\infty$-operad precisely when $V^G$ is positive-dimensional;
In particular, $\Comm_G^{\otimes} \simeq \EE_{\infty \rho}^{\otimes} \simeq \cN_{\FF_{\cT} \infty}^{\otimes}$.
Moreover, $\EE^{\otimes}_\infty$ presents the initial $\cN_\infty$-operad, and its algebras are \emph{naive} commutative algebra objects \EBVex{3.24}:
\[
  \Alg_{\EE_\infty}(\cC) \simeq \CAlg(\cC_G).
\]
If $I$ is an indexing category, the structure of an $\cN_{I \infty}$-ring spectrum is intuitively viewed as commutative ring structures on each spectrum $X_H$, connected by multiplicative $I$-indexed norms, suitably compatible with the restriction and (additive) transfer structures inherent to $G$-spectra.
We refer to $\cN_{I\infty}$-algebras in general as \emph{$I$-commutative algebras} and $\cN_{I \infty}$-ring spectra as \emph{$I$-commutative ring spectra}, writing
\[ 
    \CAlg_I(\cC) \deq \Alg_{\cN_{I\infty}}(\cC).
\]

In this paper, we are primarily concerned with homotopy coherently interchanging $\cO-$ and $\cP$-algebra structures, which are implemented as algebras over \emph{Boardman-Vogt tensor product} $\cO^{\otimes} \obv \cP^{\otimes}$ of \cite{EBV};
in particular, we are concerned with computing $\cN_{I \infty}^{\otimes} \obv \cN_{J\infty}^{\otimes}$, which corepresents pairs of interchanging $I$- and $J$-commutative algebra structures.

To start, in \EBVcor{2.91} we characterized $I \mapsto \cN^{\otimes}_{I\infty}$ as right adjoint to the \emph{arity support} construction
\[
  A\cO \deq \cbr{T \rightarrow S \;\; \middle| \;\; \prod_{[H/K] \in \Orb(S)} \cO\prn{T \times_S [H/K]} \neq \emptyset} \subset \FF_{G};
\]
when $\cO^{\otimes},\cP^{\otimes}$ have one object, we will show that $A\prn{\cO \otimes \cP} = A\cO \vee A\cP$, the latter denoting the join in the poset of weak indexing category.
This constructs a unique pairing
$\cN_{I \infty}^{\otimes} \obv \cN_{J \infty}^{\otimes} \rightarrow \cN^{\otimes}_{I \vee J \infty}$.

Intuitively, given an algebra with $I \vee J$-indexed norms, we may separate these into $I$- and $J$-indexed norms together with coherent homotopies witnessing interchange between the two.
Now, the transfer system for $I \vee J$ consists of those inclusions $K \subset H$ which can be factored as
\[
  K \subset K_{I1} \subset K_{J1} \subset K_{I2} \subset \cdots \subset K_{Jn} \subset H
\]
where $K_{I\ell} \subset K_{J\ell}$ is in $I$ and $K_{J \ell} \subset K_{I (\ell+1)}$ is in $J$ (see \cite[Prop~3.1]{Rubin_N});
intuition suggests that we may combine interchanging $I$- and $J$-commutative algebra structures to construct an $I\vee J$-commutative algebra structure.
Indeed, Blumberg and Hill conjectured that there is an equivalence $\cN_{I \infty}^{\otimes} \obv \cN_{J \infty}^{\otimes} \simeq \cN_{I \vee J \infty}^{\otimes}$ 
\Blumberg{Conj.}{6.27};
the main theorem of this paper confirms their conjecture in $\Op_G$, as well as characterizing exactly how far we may weaken $I$ and $J$.

\subsection*{Summary of main results}
Recall that a weak indexing category $I \subset \FF_G$ is \emph{almost essentially unital} if whenever a non-isomorphism $T \sqcup T' \rightarrow S$ lies in $I$, the factor map $T \rightarrow S$ lies in $I$, and \emph{almost-unital} if additionally $*_G \in I$.
We begin with a rigidity result for (co)cartesian $I$-symmetric $\infty$-categories under almost-unitality.
\begin{cooltheorem}\label{Cocartesian rigidity main theorem}
  When $I$ is almost-unital, there are fully faithful embeddings $(-)^{I-\sqcup}$ and $(-)^{I-\times}$ making the following commute:
  \[\begin{tikzcd}[column sep = huge]
	{\Cat_{I}^{\sqcup}} & {\Cat_{I}^{\otimes}} & {\Cat_I^{\times}} \\
	& {\Cat_G}
	\arrow["{(-)^{I-\sqcup}}", hook', from=1-1, to=1-2]
	\arrow["U"', from=1-1, to=2-2]
	\arrow["U"{description}, from=1-2, to=2-2]
	\arrow["{(-)^{I-\times}}"', hook, from=1-3, to=1-2]
	\arrow["U", from=1-3, to=2-2]
\end{tikzcd}\]
  The essential image of $(-)^{I-\sqcup}$ is spanned by the $I$-symmetric monoidal $\infty$-categories whose $I$-indexed tensor products are indexed coproducts, and $(-)^{I-\times}$ by those whose $I$-indexed tensor products are indexed products.
\end{cooltheorem}
\begin{remark*}
    After this introduction, we replace $\cO_G$ with an atomic orbital $\infty$-category $\cT$; we prove \cref{Cocartesian rigidity main theorem} as well as the other theorems in this introduction in this setting, greatly generalizing the stated results at the cost of ease of exposition.
\end{remark*}
We refer to $I$-symmetric monoidal $\infty$-categories of the form $\cC^{I-\times}$ as \emph{cartesian}, and $\cC^{I-\sqcup}$ \emph{cocartesian}.
\begin{remark*}
  Given $I$-symmetric monoidal $\infty$-categories $\cC^{\otimes}$ and $\cD^{\otimes}$ and an $I$-product-preserving functor $F\colon \cC \rightarrow \cD$ between their underlying $G$-$\infty$-categories, we may define the \emph{$\infty$-category of $I$-symmetric monoidal lifts}
  \[\begin{tikzcd}[ampersand replacement=\&]
	{\Fun_{G}^{\otimes, F}\prn{\cC^{\otimes}, \cD^{\otimes}}} \& {\Fun_{G}^{\otimes}\prn{\cC^{\otimes}, \cD^{\otimes}}} \\
	{\cbr{F}} \& {\Fun_{G}\prn{\cC, \cD}}
	\arrow[from=1-1, to=1-2]
	\arrow[from=1-1, to=2-1]
	\arrow["\lrcorner"{anchor=center, pos=0.125}, draw=none, from=1-1, to=2-2]
	\arrow[from=1-2, to=2-2]
	\arrow[from=2-1, to=2-2]
\end{tikzcd}\]
To interpret \cref{Cocartesian rigidity main theorem} as a rigidity theorem, note that it directly implies that whenever $\cC^{\otimes}$ and $\cD^{\otimes}$ are cartesian (resp. cocartesian), the core space $\Fun_G^{\otimes, F}\prn{\cC^{\otimes}, \cD^{\otimes}}^{\simeq}$ is contractible if $F$ is $I$-product preserving ($I$-coproduct preserving) and empty otherwise.
Moreover, we confirm this fact without taking cores in \cref{Cocartesian rigidity pre-prop,Cocartesian rigidity dual pre-prop}.
\end{remark*}
To state our remaining theorems, we need the following definition.
\begin{definition*}
  An $I$-operad $\cO^{\otimes}$ is \emph{unital} if the unique map $f\colon \cO^{\otimes} \rightarrow \cN_{I\infty}$ induces an equivalence 
  \[
      \cO(\emptyset_H;C) \simeq \cN_{I\infty}(\emptyset_H) 
  \]
  for all $H \subset G$ and $C \in \cO_H$ (c.f. \cref{Ninfty structure spaces equation});
  an $I$-operad is \emph{reduced} if additionally $f$ induces an equivalence
  \[
    \cO(*_H) \simeq \cN_{I\infty}(*_H).
  \]
  A $G$-operad $\cO^{\otimes}$ is \emph{almost essentially unital} (resp \emph{almost essentially reduced}) if it's unital (reduced) as an $A\cO$-operad and $A\cO$ is almost essentially unital.
\end{definition*}

Algebraically, we identify cartesian $I$-commutative algebras with $I$-commutative monoids and cocartesian (unital) $I$-commutative algebras with $G$-objects, 
identifying \cref{Lawvere perspective,Semiadditive perspective,Operad perspective}.
\begin{cooltheorem}\label{Cartesian algebras main theorem}
  If $I$ is almost-unital, $\cC^{\otimes}$ is a cartesian $I$-symmetric monoidal $\infty$-category, and $\cO^{\otimes}$ is an $I$-operad, then the forgetful functor
  \[
    U\colon \Alg_{\cO}(\cC) \longrightarrow \Fun_G(\Tot_{G} \cO^{\otimes},\cC)
  \]
  is fully faithful with image spanned by the $G$-functors $\Tot_G \cO^{\otimes} \rightarrow \cC$ sending $S$-indexed tuples to $S$-indexed products;
  in particular, this specializes to an equivalence
  \[
    \CAlg_I\prn{\cC^{I-\times}} \xrightarrow{\;\; \sim \;\;} \uFun_G^{I-\oplus}\prn{\uFF_{I,*}, \cC}.
  \]
  In particular, in the case of coefficient systems, we acquire an equivalence
  \[
    \Alg_{\cO}\prn{\uCoFr^G \cD^{I-\times}} \simeq \Seg_{\Tot \Tot_{G} \cO}(\cD) \simeq \Seg_{\Tot \cO}(\cD),
  \]
  where $\Seg_{(-)}(-)$ refers to \emph{Segal objects} in the sense of \cite{Chu}.
  Hence there is an additional equivalence
  \[
    \CAlg_I\prn{\uCoFr^G \cD^{I-\times}} \simeq \CMon_I(\cD). 
  \]
  Moreover, for all unital $I$-operads $\cO^{\otimes}$, the forgetful functor yields an equivalence
  \[
    \Alg_{\cO}\prn{\cC^{I-\sqcup}} \xrightarrow{\;\;\;\; \sim \;\;\;\;} \Fun_G(U\cO,\cC).
  \]
\end{cooltheorem}
\begin{proof}[References]
  This is \cref{Cartesian algebras,Tot tot corollary,CMon cor 1,CMon is CAlg corollary,Cocartesian algebras computation}.
\end{proof}

In this theorem, $\Tot \Tot_{G} \cO$ is the total $\infty$-category of the fibration over $\uFF_{G,*}$ and $\Tot \cO$ is the total $\infty$-category of the fibration over $\Span(\FF_G)$.
\begin{remark*}
  The composed equivalence $\uFun_G^{I-\oplus}\prn{\uFF_{I,*}, \uCoFr^G \cD} \simeq \Fun^{\times}\prn{\Span_I(\FF_G), \cD}$ is not new;
  indeed, it was claimed for the complete weak indexing system as far back as \cite{Nardin-Stable}, it was proved in greater generality than this article in \cite{Cnossen_semiadditive}, and we verified in \EBVthm{A.1} that it also follows from \cite{Barkan}, as well as the more general comparison between the two Segal object models for cartesian algebras.
  The new content is the identification of these notions with $G$-operad algebras.

  Moreover, in the case that $\cD = \cS$ and that $I$ is an indexing category, this is a direct analog to \cite[Thm.~\href{https://arxiv.org/pdf/2402.12447v1\#theoo.1}{A}]{Marc} in the $\infty$-categorical setting;
  the reader should interpret this relationship as a lift of Pavlov-Scholbach's comparison result \cite[Thm~1.3]{Pavlov} for a particularly nice choice of $G$-operad and value category.
\end{remark*}

In \cref{Cartesian subsection} we verify that $\uAlg_{\cO}^{\otimes}(\cC)$ is cartesian when $\cC$ is.
Following this, in \cref{Cocartesian subsection} we show that $I$-indexed tensor products in $\uCAlg_I^{\otimes} \cC$ are indexed coproducts (i.e. its underlying $I$-symmetric monoidal $\infty$-category is \emph{cocartesian}) and that this completely characterizes $\cN_{I\infty}^{\otimes}$.
The heart of our strategy uses the explicit monadic description of \EBVcor{2.73} to reduce to the case of $G$-spaces $\cC^{\otimes} \simeq \ucS_{G}^{G-\times}$;
in this case, we see that the cartesian $I$-symmetric monoidal $\infty$-category $\uCAlg_I^{\otimes}(\ucS_G^{G-\times}) \simeq \uCMon_I(\ucS_G)^{I-\times}$ is cocartesian, as its underlying $G$-$\infty$-category is $I$-semiadditive by 
\cite[Thm.~\href{https://arxiv.org/pdf/2403.07676v3\#introthm.2}{B}-\href{https://arxiv.org/pdf/2403.07676v3\#introthm.3}{C}]{Cnossen_semiadditive}
We conclude the following.
\begin{cooltheorem}\label{Cocartesian main theorem}
    Let $\cO^{\otimes}$ be an almost essentially reduced $G$-operad.
    Then, the following conditions are equivalent.
    \begin{enumerate}[label={(\alph*)}]
        \item The $G$-$\infty$-category $\uAlg_{\cO} \ucS_G$ is $A\cO$-semiadditive.
        \item The unique map $\cO^{\otimes} \rightarrow \cN^{\otimes}_{A\cO \infty}$ is an equivalence.
    \end{enumerate}
    Moreover, for all almost essentially unital weak indexing categories $I$ and $I$-symmetric monoidal $\infty$-categories $\cC^{\otimes}$, the $I$-symmetric monoidal $\infty$-category $\uCAlg_{I}^{\otimes} \cC$ is cocartesian.
\end{cooltheorem}
\cref{Cartesian algebras main theorem,Cocartesian main theorem} together with conservativity of $\uAlg_{(-)}(\ucS_G)$ \EBVcor{2.74} yields the following.
\begin{coolcorollary}\label{Absorption main corollary}
    $\cN_{I\infty}^{\otimes} \obv \cN_{I \infty}^{\otimes}$ is a weak $\cN_\infty$-operad if and only if $I$ is almost essentially unital.
    In this case, if $\cO^{\otimes}$ is a reduced $I$-operad, then the unique map $\cO^{\otimes} \otimes \cN_{I\infty}^{\otimes} \rightarrow \cN_{I \infty}^{\otimes}$ is an equivalence.
\end{coolcorollary}

In particular, whenever $I$ is almost unital, there exists a map $\triv_G^{\otimes} \rightarrow \cN_{I \infty}^{\otimes}$ witnessing $\cN_{I \infty}^{\otimes}$ as an idempotent object in $\Op_G$.
We verified in \cite[Thm.~\href{https://arxiv.org/pdf/2501.02129\#cooltheorem.4}{D}]{EBV} that $\Env\cln \Op_G \rightarrow \Cat_G^{\otimes}$ is compatible with the unit and tensor products under the mode symmetric monoidal structure on $\Cat_G^{\otimes}$;
this yields a $\circledast$-idempotent algebra structure on $\uFF_{G}^{G-\sqcup} = \Env(\Comm_G) \in \Cat_G^{\otimes}$, and hence a symmetric monoidal structure on $\uCat^{\otimes}_{G, /\uFF_G^{G-\sqcup}}$.
We acquire an equivariantization of a modification of \cite[Thm.~\href{https://arxiv.org/pdf/2301.08650v3\#thmA.5}{E}]{Barkan_Segal}.

\begin{coolcorollary}\label{Operad main theorem}
   There exists a unique symmetric monoidal structure $\uOp_G^{\otimes}$ on $\uOp_G$ attaining a (necessarily unique) symmetric monoidal structure on the fully faithful $G$-functor
   \[
     \Env^{/\uFF_{G}^{G-\sqcup}}\cln \uOp_{G}^{\otimes} \longrightarrow \uCat^{\otimes}_{G, /\uFF_{G}^{G-\sqcup}}
   \]
   of \cite{Nardin,Barkan} with respect to $\circledast$;
   the tensor product of this structure is $\obv$.
\end{coolcorollary}

Idempotent objects correspond with smashing localizations, i.e. they classify particular properties \cite[\S~\href{https://people.math.harvard.edu/~lurie/papers/HA.pdf\#subsection.4.8.2}{4.8.2}]{HA};
in \cref{The smashing localization theorem}, we conclude that the smashing localization corresponding with $\cN^{\otimes}_{I \infty} \in \Op_J^{\red}$ classifies the property of \emph{having $I$-indexed Wirthm\"uller isomorphisms}
\begin{align*}
  \cO^{\otimes} \obv \cN_{I \infty}^{\otimes} \simeq \cO^{\otimes} \;\;\; 
  &\iff \;\;\; \forall \cC^{\otimes} \in \Cat_J^{\otimes},\;\;\; \forall S \in \FF_{I,V}, \;\;\; \forall (X_U)_S \in \uAlg_{\cO}(\cC)_S \;\;\; W_{S}\colon \coprod\limits_U^S X_U \xrightarrow{\;\; \sim \;\;} \bigotimes_U^S X_U \\
  &\iff \;\;\; \uAlg_{\cO}(\ucS_{G}) \text{ is } I \text{-semiadditive}.
\end{align*}
Recall that tensor products of idempotents algebras are idempotent algebras, classifying the intersection of the associated smashing localizations \csy{Prop.}{5.1.8};
conveniently, indexed semiadditivity is classified by a weak indexing category \EBVprop{1.52}, so $\uAlg_{\cO}(\ucS_G)$ is $I \vee J$-semiadditive if and only if it is $I$-semiadditive and $J$-semiadditive.
This allows us to affirm
Blumberg and Hill's conjecture with respect to $\obv$.
\begin{cooltheorem}\label{Main theorem}
    $\cN_{(-)\infty}^{\otimes}\colon \wIndSys_G \rightarrow \Op_G$ restricts to a fully faithful symmetric monoidal $G$-right adjoint
    \[
        \begin{tikzcd}
          {\uwIndSys^{aE\uni}_G} && {\uOp^{aE\uni}_G}
            \arrow[""{name=0, anchor=center, inner sep=0}, "{\cN_{(-)\infty}^{\otimes}}"', curve={height=20pt}, hook, from=1-1, to=1-3]
            \arrow[""{name=1, anchor=center, inner sep=0}, "A"', curve={height=20pt}, from=1-3, to=1-1]
            \arrow["\dashv"{anchor=center, rotate=-90}, draw=none, from=0, to=1]
        \end{tikzcd}
    \]
    Furthermore, the resulting tensor product of weak $\cN_\infty$-operads is computed by the Borelified join
    \begin{align*}
      \cN_{I\infty}^{\otimes} \obv \cN_{J\infty}^{\otimes} \simeq \cN^{\otimes}_{\Bor_{cI \cap cJ}^{G} \prn{I \vee J}\infty}.
    \end{align*}
    Hence when $I,J$ are almost-unital weak indexing categories and $\cC^{\otimes}$ is an $I \vee J$-symmetric monoidal $\infty$-category, there is a canonical equivalence of $I \vee J$-symmetric monoidal $\infty$-categories
    \[
      \uCAlg^{\otimes}_I \uCAlg^{\otimes}_J(\cC) \simeq \uCAlg^{\otimes}_{I \vee J}(\cC).
    \]
\end{cooltheorem}
For instance, using \chll{Thm.}{4.3.6} to identify $I$-Tambara functors in an $\infty$-category $\cC$ with $I$-commutative algebras in Mackey functors, this confirms that $I \vee J$-Tambara functors are equivalent to $I$-commutative algebras in $J$-Tambara functors with respect to the box product.
\begin{remark*}
  The reader interested in computing tensor products of $G$-operads may benefit from reading the combinatorial characterization of joins of weak indexing systems in terms of \emph{closures} in \Windprop{2.22};
  there, we prove that the join of weak indexing systems $\FF_{I} \vee \FF_{J}$ is computed by closing the union $\FF_I \cup \FF_J$ under iterated $I$ and $J$-indexed coproducts.
\end{remark*}

\subsection*{Relationship to the literature}
There are three main bodies of literature which present results in homotopy-coherently equivariant algebra:
the model categorical, the atomic orbital, and the global.
We now attempt to give a bit of a Rosetta stone to connect our definitions to the model categorical and global settings.

We established in \cite{Windex} that our weak indexing categories specialize to Blumberg-Hill's indexing categories \cite{Blumberg_incomplete} in the case $\cT = \cO_G$ and $n \cdot *_G \rightarrow *_G$ lies in $I$ for all $n \in \NN$, and our weak indexing systems to the indexing systems of \cite{Blumberg-op} when $n \cdot *_G \in \FF_{I,G}$ for all $n \in \NN$;
moreover, this was shown to be compatible with Bonventre's nerve in \cite{Bonventre-nerve,EBV}, which is intertwines with the underlying $G$-symmetric sequence and restricting to an equivalence on at-most-one-color $G$-operads with $0$-truncated structure spaces, showing that our weak $\cN_\infty$-operads specialize to those of \cite{Blumberg_incomplete,Rubin_N,Gutierrez,Bonventre}.
Additionally, we saw in \cite{Windex} that our weak indexing categories specialize the \emph{weakly extensive span pairs} of \cite{Cnossen_tambara} to the case that the larger category is $\FF_{\cT}$.

We saw in \cite{EBV} that our algebras agree with Blumberg-Hill's in the discrete setting, and combining \cref{CMon is CAlg corollary} with the main result of \cite{Marc} identifies a Dwyer-Kan localization of the latter with the former in the case $\ucS_{G}^{G-\times}$.
Of course, the recent results of Lenz, Linskens, and P\"utszst\"uck \cite[Thm~A]{Lenz_norms} have established \emph{rectification} to $G$-commutative ring spectra, establishing them as presented by Hammock localization of a right-transferred structure on commutative algebras,
$\CAlg(\Sp_G^{\Sigma})$, with respect to the positive stable model structure on symmetric $G$-spectra (or equivalently, on orthogonal $G$-spectra).
The author is not aware of any study into rectification for the incomplete case.

Our $I$-symmetric monoidal $\infty$-categories and $I$-operads generalize \cite{Nardin} and specialize \cite{Lenz_norms} by \cite[\S~A]{EBV} and by definition, respectively.
Moreover, there is a homotopical \emph{operadic nerve} construction mapping the (equivalent) settings of \cite{Blumberg-op,Pereira,Bonventre} to $\Op_I$ and $\Cat_I^{\otimes}$ by \cite{Bonventre-nerve,EBV}.

The author is not aware of a comparison result between $\obv$ and the point-set Boardman-Vogt tensor product appearing in \cite{Blumberg-op};
moreover, the ``derived'' tensor product appearing in \cite{Rubin_N} is only defined on an $\infty$-category which is equivalent to $\IndSys_G$, so it's not clear that it makes sense to ask for a comparison to $\obv$ other than confirming that $\obv$ confirms Blumberg-Hill's conjecture, as demanded by the results of \cite[Thm~A]{Rubin_N} (after which Rubin explicitly claimed that the conjecture remained open).

For (co)cartesian $I$-symmetric monoidal $\infty$-categories, we show in \cref{Cartesian subsection} that our definitions generalize \cite{Nardin} and agree with \cite{Cnossen_tambara} when both are defined;
in particular, we recover a non-$(\infty,2)$-categorical reproof of the identification theorem of the two, as conjectured in \cite{Cnossen_tambara} and verified in \cite{Cnossen_unfurling}.
The author is not aware of a definition to be compared with in the model-categorical setting,
but comparisons with such constructions will be as easy as verifying that those structures have indexed tensor products which present (derived) indexed products.

Some versions of our results on cocartesianness and algebras are proved independently in the literature;
though it is not clear that Nardin-Shah's $\cT$-operad of algebras \cite{Nardin} agrees with ours, they confirmed that their version of $\uCAlg_{\cT}^{\otimes}(\cC)$ is cocartesian and claimed that $\uCAlg_I^{\otimes}(\cC)$ is $I$-cocartesian when $I$ is an indexing category.
Moreover, it is shown in \cite[Prop~2.26]{Lenz_norms} that $(-)^{I-\sqcup}$ \emph{admits} a left adjoint, but this left adjoint is not computed therein (whereas we confirm it to be $U$ in the almost-unital case as \cref{Left adjoint to U corollary}).  

Additionally, after the first preprint version of this article appeared on the arxiv, the author learned of independent accounts of some aspects the cartesian and cocartesian structure in \cite[\S~A]{Hilman_McDuff} and \cite[\S~4.1]{Yang_HKR}, respectively.
In particular, the former recovers \cref{CMon is CAlg corollary} and the latter \cref{Cocartesian algebras computation} in the case $I = \cT = \cO_G$.

\subsection*{Notation and conventions}
We assume that the reader is familiar with the technology of higher category theory and higher algebra as developed in \cite{HTT} and \cite[\S~2-3]{HA}, though we encourage the reader to engage with such technologies via a ``big picture'' perspective akin to that of \cite[\S~1-2]{Gepner_HA} and \cite[\S~1-3]{Haugseng}.
In particular, our treatment is almost entirely model agnostic--we only pierce the veil in \cref{Quasicat subsection} and use quasicategorical language in order to verify that a few functors are exponentiable.

We additionally assume that the reader is familiar with \emph{parameterized} higher category theory over an $\infty$-category as developed in \cite{Shah,Shah2};
the material reviewed in the prequel \cite[\S~1]{EBV} will be enough.
In particular,
\begin{itemize}
  \item $\cT$ will always be an atomic orbital $\infty$-category in the sense of \cite{Nardin}, $\FF_{\cT}$ its corresponding $\infty$-category of finite $\cT$-sets, and $\uFF_{\cT}$ its corresponding $\cT$-1-category of finite $\cT$-sets.
  \item $\cF \subset \cT$ will always be a $\cT$-family in the sense of \cite{Windex}.
  \item $I \subset \FF_{\cT}$ and $\uFF_I$ will always be a weak indexing category and corresponding weak indexing system in the sense of \cite{Windex}. $c(I)$ will be its color family and $\upsilon(I)$ its unit family.
  \item $\Cat$ will always be the $\infty$-category of small $\infty$-categories, $\Cat_{\cT}$ of small $\cT$-$\infty$-categories, and $\Cat_I^{\otimes}$ of small $I$-symmetric monoidal $\infty$-categories.
    All $\infty$-categories will be assumed to be small unless otherwise mentioned.
  \item $\cT$-operad will always mean $\cT$-$\infty$-operad in the sense of \cite{Nardin} and $\Op_{\cT}$ the $\infty$-category of $\cT$-operads.
\end{itemize}
\subsection*{Acknowledgements}
I'm indebted to Maxime Ramzi for disillusioning me to a fatally flawed strategy on work related to this paper, leading me to the drawing board;
it wasn't so clear at the time, but it was in subsequent conversation with him that the main idea of this paper emerged.
I am additionally grateful to Piotr Pstr\k{a}gowski, who pointed out a mistake in my early strategy in this paper, leading to the condition of \emph{almost essential unitality} on the main theorem.
Also, I owe Mike Hill for pointing out to me that pullback-stable subcategories are replete, obviating one of the assumptions on weak indexing categories.

Additionally, I would like to thank Andy Senger, Clark Barwick, and Dhilan Lahoti, with whom I had enlightening (to me) conversations about the topic of this paper. 
Of course, none of this work would be possible without the help of my advisor, Mike Hopkins, who I'd like to thank for many helpful conversations.

While developing this material, the author was supported by the NSF Grant No. DGE 2140743.

\resumetocwriting 
\section{\tI-symmetric monoidal categories and \tI-operads}\label{First section}
We begin in \cref{Recollections on CommI subsecttion} by recalling results of \cite{Nardin-Stable,Nardin,Cnossen_semiadditive,Windex,EBV} concerning the theory of $I$-commutative monoids and $I$-symmetric monoidal $\infty$-categories.
Moving on, in \cref{Recollections on T-operads subsection} we recall results of \cite{Nardin,EBV} concerning $\cT$-operads;
in either case, all reviewed information was used in the preceding article \cite{EBV}.
We then go on in \cref{Arity borelification section} to begin to carefully study the interactions of restriction, arity-borelification, arity-support, and Boardman-Vogt tensor products.
We finish the section in \cref{Cartesian subsection}, where we develop a number of foundational results on (co)cartesian $I$-symmetric monoidal $\infty$-categories, ultimately elaborating on the technical minutiae of \cref{Cocartesian proof subsubsection}.
\subsection{Recollections on \tI-commutative monoids and \tI-symmetric monoidal \tinfty-categories}\label{Recollections on CommI subsecttion}
For the rest of this paper, we fix $\cT$ an atomic orbital $\infty$-category.
\subsubsection{Weak indexing systems and semiadditivity}
We will use the following machinery of \cite{Windex}.
\begin{definition}\label{Windex definition}
  A \emph{$\cT$-weak indexing category} is a subcategory $I \subset \FF_{\cT}$ satisfying the following conditions:
  \begin{enumerate}[label={(IC-\alph*)}]
    \item \label[condition]{Restriction stable condition} (restrictions) $I$ is stable under arbitrary pullbacks in $\FF_{\cT}$, and
    \item \label[condition]{Windex segal condition} (segal condition) $T \rightarrow S$ and $T' \rightarrow S$ are both in $I$ if and only if $T \sqcup T' \rightarrow S \sqcup S'$ is in $I$.
  \end{enumerate}
  A \emph{$\cT$-weak indexing system} is a full $\cT$-subcategory $\uFF_{I} \subset \uFF_{\cT}$ satisfying the following conditions:
  \begin{enumerate}[label={(IS-\alph*)}]
    \item \label[condition]{Point in condition} whenever the $V$-value $\FF_{I,V} \deq \prn{\uFF_I}_V$ is nonempty, we have $*_V \in \uFF_{I,V}$, and
    \item \label[condition]{Coproducts condition}$\uFF_{I} \subset \uFF_{\cT}$ is closed under $\uFF_I$-indexed coproducts.\qedhere
  \end{enumerate}
  We say that a $\cT$-weak indexing system $\uFF_I$:
  \begin{enumerate}[label={(\roman*)}]
    \item has one color if for all $V \in \cT$, we have $\FF_{I,V} \neq \emptyset$,
    \item is almost essentially unital (or a$E$-unital) if whenever $\uFF_I$ has a non-contractible $V$-set, $\emptyset_V \in \FF_{I,V}$,
    \item is almost-unital (or a-unital) it's almost essentially unital and has one color,
    \item is unital if $\emptyset_V \in \FF_{I,V}$ for all $V \in \cT$, and
    \item is an \emph{indexing system} if the subcategory $\uFF_{I,V} \subset \FF_V$ is closed under finite coproducts for all $V \in \cT$.
  \end{enumerate}
  These occupy embedded sub-posets
  \[
    \IndSys_{\cT} \subset \wIndSys_{\cT}^{\uni} \wIndSys_{\cT}^{a\uni}  \subset \wIndSys_{\cT}^{aE\uni} \subset \wIndSys_{\cT}.\qedhere
  \]
\end{definition}

Given a weak indexing category $I \subset \FF_{\cT}$, we denote the \emph{$I$-admissible $V$-sets} by
\[
  \FF_{I,V} \deq \cbr{S \in \FF_{I,V} \mid \Ind_V^{\cT} S \rightarrow V \in I} \subset \FF_V;
\]
\cref{Restriction stable condition} guarantees that these assemble into a full $\cT$-subcategory $\uFF_I \subset \uFF_{\cT}$, which contains all of the information of $I$ by \cref{Windex segal condition}.
In \cite[Thm~A]{Windex} we proved the following and expressed the conditions of \cref{Windex definition} in the language of weak indexing categories.
\begin{proposition}[{Generalized \cite[Thm~1.4]{Blumberg_incomplete}}]
  The assignment $I \mapsto \uFF_I$ implements an equivalence between the posets of $\cT$-weak indexing categories and $\cT$-weak indexing systems.
\end{proposition}

Intuitively, \cref{Point in condition} corresponds with identity arrows in $I$ and \cref{Coproducts condition} with composition.
We will need the following invariants of weak indexing systems.
\begin{construction}
  Given $\uFF_I$ a weak indexing system, we define the \emph{color and unit families}
  \begin{align*}
    c(I) &\deq \cbr{V \in \cT \mid *_V \in \FF_{I,V}} \subset \cT;\\
    \upsilon(I) &\deq \cbr{V \in \cT \mid \emptyset_V \in \FF_{I,V}} \subset \cT.
  \end{align*}
  Indeed, we saw that these are families in \cite{EBV}.
\end{construction}

One reason to study this is \emph{indexed semiadditivity} in the sense of the following definitions of Nardin.
\begin{definition}\label{Wirth defn}
  A $\cT$-$\infty$-category $\cC$ is said to be \emph{$V$-pointed} if $\cC_U$ is a pointed $\infty$-category for all $U \rightarrow V$. 
  Given $S \in \FF_V$ is a finite $V$-set, $\cC$ $V$-pointed admitting $S$-indexed products and coproducts, and $(X_U)_S \in \cC_S$ an $S$-tuple, we define the \emph{$S$-indexed Wirthmuller map} $W_{S, (X_U)}\colon \coprod_U^S X_U \rightarrow \prod_U^S X_U$ to extend the following maps via the universal property for $S$-indexed coproducts:
  \[
    W_{S,(X_U),W}\colon X_W \simeq X_W \times \prod_{U'}^{\Res_W^V S - W} *_{U'} \xrightarrow{\;\;\;\; (\id;!)\;\;\;\;} X_W \times \prod_{U'}^{\Res_U^V S - W} \Res_{U'}^V X_{o(U')} \simeq \Res_W^V \prod_U^S X_U
  \]
  where $o(U') \in \Orb(S)$ is the orbit whose restriction contains $U'$.
  We say that $S$ is \emph{$\cC$-ambidextrous} if $\cC$ is $V$-pointed and $W_{S, (X_U)}$ is an equivalence for all $(X_U) \in \cC_S$;
  given $\uFF_I$ a weak indexing system, we say that $\cC$ is \emph{$I$-semiadditive} if $S$ is $\cC$-ambidextrous for all $S \in \uFF_I$.
\end{definition}

\begin{remark}
  The map $W_{S, (X_U),W}$ is determined via the universal property for $S$-indexed products by its projections $W_{S,(X_U),W,W'}\colon X_W \rightarrow \Res_W^V \CoInd_{W'}^V X_{W'}$, which are zero when $W \neq W'$, and otherwise they are the map induced under functoriality of products by the dashed arrow \[\begin{tikzcd}[ampersand replacement=\&, row sep=small]
	{\Ind_W^{V} *_W} \\
	\& {\Ind_W^V \Res_W^V \Ind_W^V *_W} \& {\Ind_W^V*_W} \\
	\& {\Ind_W^{V} *_W} \& V
	\arrow[from=1-1, to=2-2]
	\arrow[curve={height=-12pt}, equals, from=1-1, to=2-3]
	\arrow[curve={height=12pt}, equals, from=1-1, to=3-2]
	\arrow[from=2-2, to=2-3]
	\arrow[from=2-2, to=3-2]
	\arrow["\lrcorner"{anchor=center, pos=0.125}, draw=none, from=2-2, to=3-3]
	\arrow[from=2-3, to=3-3]
	\arrow[from=3-2, to=3-3]
\end{tikzcd}\]
In particular, they match the norms constructed in \cite{Nardin-Stable}.
\end{remark}

In \cite{EBV} we proved that the collection of $\cC$-ambidextrous finite $V$-sets form a weak indexing system and concluded the following important observation.
\begin{proposition}[{\cite[\S~1.2]{EBV}}]\label{Semiadditivity over joins}
  Let $\vee$ denote the join in $\wIndCat_{\cT}$.
  Then, $\cC$ is $I$-semiadditive and $J$-semiadditive if and only if $\cC$ is $I \vee J$-semiadditive.
\end{proposition}

\subsubsection{$I$-commutative monoids}
In \cite{Barwick1}, the notion of \emph{adequate} triple was defined, consisting of triples $(\cC,\cC_b,\cC_f)$ with $\cC_f,\cC_B \subset \cC$ a pair of wide subcategories satisfying pullback-stability and distributivity conditions;
if $I$ is a weak indexing category, then $(\FF_{c(I)},\FF_{c(I)},I)$ is an adequate triple.

Adequate triples form a full subcategory $\Trip^{\Adeq} \subset \Fun(\bullet \rightarrow \bullet \leftarrow \bullet, \Cat)$;
\cite{Barwick1} constructed a functor
\[
  \Span_{-,-}(-)\cln \Trip^{\Adeq} \rightarrow \Cat,
\]
called the \emph{effective Burnside category}.
In the case that $c(I)$ is a 1-category (e.g. $\cT$ has a terminal object, see \cite[Prop~2.5.1]{Nardin}), $\FF_{c(I)}$ is a 1-category, so the effective Burnside category 
\[
  \Span_I(\FF_{\cT}) \deq \Span_{\FF_{c(I)},I}(\FF_{c(I)})
\]
is a $(2,1)$-category with objects agreeing with $\FF_{c(I)}$, morphisms the spans $X \leftarrow R \xrightarrow{f} Y$ with $f$ in $I$, 2-cells the isomorphisms of spans, and composition of morphisms computed by pullbacks in $\FF_{c(I)}$ (which are guaranteed to be morphisms in $\Span_I(\FF_{\cT})$ by pullback-stability of $I$).

Much of the technical work of \cite{Barwick1,Barwick2} has been extended by \cite{Haugseng_two}, so we generally refer the reader there.
At any rate, we recall this in order to define \emph{homotopical incomplete semi-Mackey functors for $I$}, which we call \emph{$I$-commutative monoids}.
\begin{definition}
  If $\cC$ is an $\infty$-category with finite products, then an \emph{$I$-commutative monoid in $\cC$} is a product-preserving functor $\Span_I(\FF_{\cT}) \rightarrow \cC$.
  More generally, if $\cD$ is a $\cT$-$\infty$-category with $I$-indexed products, then an \emph{$I$-commutative monoid in $\cD$} is an $I$-product-preserving $\cT$-functor $\Span_I(\uFF_{\cT}) \rightarrow \cD$.
  We write
  \begin{align*}
    \uCMon(\cD) &\deq \uFun^{I-\times}_{\cT}\prn{\Span_I(\uFF_{\cT}), \cD}\\
    \CMon(\cD)  &\deq \Gamma^{\cT} \uCMon(\cD)\\
    \uCMon(\cC) &\deq \uCMon\prn{\uCoFr^{\cT} \cC}\\
    \CMon(\cC)  &\deq \CMon\prn{\uCoFr^{\cT} \cC}.\qedhere
  \end{align*}
\end{definition} 
An important result of Cnossen-Lenz-Linskens resolves the notational clash.
\begin{proposition}[{\cite[Thm~C]{Cnossen_semiadditive}}]\label{CMon in coeff}
  When $\cC$ is an $\infty$-category, restriction furnishes an equivalence
  \[
    \CMon(\cC) \simeq \Fun^{\times}\prn{\Span_I(\FF_{\cT}), \cC},
  \]
  and more generally, we have $\uCMon(\cC)_V \simeq \Fun_V^{\times}\prn{\Span_I(\FF_V),\cC}$ with restriction given by pullback along $\Span_I(\FF_V) \rightarrow \Span_I(\FF_W)$.
\end{proposition}

Let $I$ be a one-object weak indexing category and let $\Cat_{\cT}^{I-\times} \subset \Cat_{\cT}$ be the (non-full) subcategory whose objects are $\cT$-$\infty$-categories admitting $I$-indexed products and functors preserving $I$-indexed products.
Let $\Cat_I^{I-\oplus} \subset \Cat_{\cT}^{I-\times}$ be the full subcategory spanned by $I$-semiadditive $\cT$-$\infty$-categories.
The following result is fundamental in the theory of equivariant semiadditivity and equivariant higher algebra.
\begin{theorem}[{\cite[Thm~B]{Cnossen_semiadditive}}]\label{Semiadditivization theorem}
  The inclusion $\Cat_{\cT}^{I-\oplus} \subset \Cat_{\cT}^{I-\times}$ has left adjoint $\uCMon(-)$.
\end{theorem}

\subsubsection{$I$-symmetric monoidal $\infty$-categories}
The following definition is central to equivariant higher algebra.
\begin{definition}
  The \emph{$\infty$-category of small $I$-symmetric monoidal $\infty$-categories} is $\Cat_I^{\otimes} \deq \CMon_I(\Cat)$.
\end{definition}  
We refer to maps in $\Cat_I^{\otimes}$ as \emph{$I$-symmetric monoidal functors}.
An important lemma is the following.
\begin{lemma}[{\cite[Cor~8.2]{Chu}}]\label{Underlying conservative proposition}
  If $\cC$ is an $\infty$-category and $I$ a one-object weak indexing category, then the \textbx{underlying coefficient system} functor $\CMon_I(\cC) \rightarrow \CoFr^{\cT} \cC$ is conservative;
  in particular, if a $I$-symmetric monoidal functor's underlying $\cT$-functor is an equivalence, then it is a $\cT$-symmetric monoidal equivalence.
\end{lemma}
Now, these are defined for the following notation's sense.
\begin{notation}
  Given an $I$-symmetric monoidal $\infty$-category and an $I$-map $\Ind_V^{\cT} S \rightarrow V$ with $V \in \cT$, we denote the covariant functorialty of $\cC^{\otimes}$ by  $\bigotimes^S_U\colon \cC_S \rightarrow \cC_V$ and the contravariant functoriality by $\Delta^S\colon \cC_V \rightarrow \cC_S$.
\end{notation}
We will also need \emph{presentability}.
In his thesis, Nardin defined a $\cT$-symmetric monoidal $\infty$-category $\uPr^{L,\otimes}_{\cT}$ of \emph{presentable $\cT$-$\infty$-categories}, whose $S$-ary tensor products are characterized by mapping $\infty$-categories
\[
  \Fun^{L}_{\cT}\prn{\bigotimes_U^S \cC_U, \cD}
  \deq \Fun^{S-\partial}_{\cT}\prn{\prod_U^S \cC_U, \cD} 
\]
where $\Fun^{S-\partial}_{\cT}$ consists of the ``$S$-distributive $\cT$-functors.''
We will not care too much about the details of this in general, and instead shunt the interested the reader to \cite[Def~5.14]{Quigley_tate}.
Nevertheless, we care about the following case.
\begin{example}
  $\Fun^{2*_V-\partial}_{V}\prn{\cC \times \cD, \cE} \subset \Fun\prn{\cC\times \cD, \cE}$ is the full subcategory of $\cT$-functors whose fibers $\cC \times \cbr{D}, \cbr{C} \times \cD \rightarrow \cE$ all strongly preserve $\cT$-colimits.
\end{example}

Now, we make the following definition.
\begin{definition}
  A \emph{distributive $I$-symmetric monoidal $\infty$-category} is an $I$-symmetric monoidal $\infty$-category $\cC^{\otimes}$ whose $S$-tensor functors $\otimes_U^S\colon \cC_S \rightarrow \cC_V$ are $S$-distributive.
  A \emph{presentably $I$-symmetric monoidal $\infty$-category} is a distributive $I$-symmetric monoidal $\infty$-category whose underlying $c(I)$-$\infty$-category is presentable.
\end{definition}

We will also need the following construction.
\begin{construction}
  The ``opposite category'' construction $\op\colon \Cat \rightarrow \Cat$ is an equivalence, so in particular it is product-preserving.
  Hence postcomposition with $(-)^{\op}$ yields \emph{fiberwise opposite} functor
  \[
    (-)^{\vop}\colon \Cat_I^{\otimes} \rightarrow \Cat_I^{\otimes}.
  \]
  Note that the underlying category $\prn{\cC^{\otimes}}^{\vop}$ is the traditional \emph{vertical opposite $\cT$-$\infty$-category } $\cC^{\vop}$.
\end{construction}

\subsection{Recollections on \tcT-operads}\label{Recollections on T-operads subsection} 
\subsubsection{$I$-operads} 
In \cite{EBV}, we made the following definition.
\begin{definition}\label{T-operad definition}
  An \emph{$I$-operad} is a functor $\pi\colon \Tot \cO^{\otimes} \rightarrow \Span_I(\FF_{\cT})$ satisfying the following conditions. 
  \begin{enumerate}[label={(\alph*)}]
    \item \label[condition]{Inert cocartesian lifts} $\Tot `\cO^{\otimes}$ has $\pi$-cocartesian lifts for backwards maps in $\Span_I(\FF_{\cT})$;
    \item \label[condition]{Color Segal condition} (Segal condition for colors) writing $\cO_S \simeq \pi^{-1}(S)$, for every $S \in \FF_\cT$, cocartesian transport along the $\pi$-cocartesian lifts lying over the inclusions $\prn{S \leftarrow U = U \mid U \in \mathrm{Orb}(S)}$ together induce an equivalence
      \[
        \cO_S \simeq \prod_{U \in \mathrm{Orb}(S)} \cO_U; 
      \]
    \item \label[condition]{Multimorphism Segal condition} (Segal condition for multimorphisms)  for every map of orbits $T \rightarrow S$ in $I$ and pair of objects $(\bC,\bD) \in \cO_T \times \cO_U$, postcomposition with the $\pi$-cocartesian lifts $\bD \rightarrow D_U$ lying over the inclusions 
    $\left(S \leftarrow U = U \mid U \in \mathrm{Orb}(S)\right)$ induces an equivalence
      \[
        \Map^{T \rightarrow S}_{\Tot \cO^{\otimes}}(\bC,\bD) \simeq \prod_{U \in \mathrm{Orb}(S)} \Map^{T \leftarrow T_U \rightarrow  U}_{\Tot \cO^{\otimes}}(\bC,D_U).
      \]
      where $T_U \deq T \times_S U$. 
  \end{enumerate}
  The $\infty$-category of $I$-operads is defined to be a localizing subcategory 
\[\begin{tikzcd}[ampersand replacement=\&]
	{\Op_I} \& {\Cat^{\Int-\cocart}_{/\Span_I(\FF_{\cT})};}
	\arrow[curve={height=16pt}, hook, from=1-1, to=1-2]
	\arrow["{L_{\Op_I}}"', curve={height=16pt}, from=1-2, to=1-1]
\end{tikzcd}\]
  that is, a morphism of $I$-operads is a functor $\Tot \cO^{\otimes} \rightarrow \Tot \cP^{\otimes}$ over $\Span_I(\FF_{\cT}) $sending $\pi_{\cO}$-cocartesian morphisms to $\pi_{\cP}$-cocartesian morphisms.
  We also call these \emph{$\cO$-algebras in $\cP$} and we let 
  \[
    \Alg_{\cO}(\cP) \deq \Fun_{/\Span_I(\FF_{\cT})}^{\Int-\cocart}(\cO^{\otimes},\cP^{\otimes}) \subset \Fun_{/\Span_I(\FF_{\cT})}(\cO^{\otimes},\cP^{\otimes})
  \]
  be the full subcategory spanned by $\cO$-algebras in $\cP$.
\end{definition}
This doesn't obviously recover the notion of \cite{Nardin}.
To discuss the comparison, we temporarily assume the reader is familiar with \emph{fibrous patterns} in the sense of \cite{Barkan} (which are essentially \emph{weak Segal fibrations} in the sense of \cite{Chu}).
\begin{construction}\label{FI construction}
  Let $\Tot \uFF_{\cT}^{\vee}$ be the cartesian unstraightening of the functor $V \mapsto \FF_V$, so that its objects are $\cT$-arrows $ S \rightarrow U$ with $V \in \cT$ and its morphisms $f\colon \prn{T \rightarrow V} \longrightarrow \prn{S \rightarrow U}$ are commutative diagrams between arrows, i.e.
  \[\begin{tikzcd}[ampersand replacement=\&]
	T \& {S \times_U V} \& S \\
	\& V \& U
	\arrow["{f^{\circ}}"{description}, from=1-1, to=1-2]
	\arrow["{f_{s}}", curve={height=-15pt}, from=1-1, to=1-3]
	\arrow[from=1-1, to=2-2]
	\arrow[from=1-2, to=1-3]
	\arrow[from=1-2, to=2-2]
	\arrow["\lrcorner"{anchor=center, pos=0.125}, draw=none, from=1-2, to=2-3]
	\arrow[from=1-3, to=2-3]
	\arrow["{f_t}", from=2-2, to=2-3]
\end{tikzcd}\]
  We say that $f$ is \emph{s.i.} if $f^{\circ}$ is a summand inclusion and \emph{$I$-tdeg} if $f_t$ is an identity arrow and $f_s$ lies in $I$.
  Then, we define the algebraic pattern
  \[
    \Tot \uFF_{I,*} \deq \Span_{s.i., I-tdeg}\prn{\Tot \FF_{\cT}^\vee}.
  \]
  The map of triples $(\Tot \uFF_{\cT}^{\vee}, s.i., I-tdeg) \rightarrow (\FF_{\cT}, all, I)$ induces a Segal morphism $s\colon \Tot \uFF_{I,*} \rightarrow \Span_I(\FF_{\cT})$.  
\end{construction}
We recover Nardin-Shah's notion of $\cT$-$\infty$-operads by the following result.
\begin{proposition}[{\cite{Barkan,EBV}}]
  The $\infty$-category $\Op_{\cT,\infty}$ of \cite{Nardin} is equivalent to fibrous $\Tot \uFF_{\cT,*}$-patterns, and $s$ induces an equivalence
  \[
    \Op_{I} \simeq \Fbrs\prn{\Span_I(\FF_{\cT})} \xrightarrow{\;\;\;\; \sim \;\;\;\;} \Fbrs\prn{\Tot \uFF_{I,*}}.
  \]
\end{proposition}
In particular, we get a composite functor
\[
  \Tot_{\cT}\colon \Op_I \hookrightarrow \Cat^{\Int-\cocart}_{/\Span_I(\FF_{\cT})} \xrightarrow{\;\; s^* \;\;} \Cat_{\cT, /\uFF_{\cT,I}} \xrightarrow{\;\;U\;\;} \Cat_{\cT}
\]
The following observation about this composite functor is key.
We greatly strengthen  it in \cref{sec:morita}.
\begin{observation}
  $\Tot_{\cT}\colon \Op_I \rightarrow \Cat_{\cT}$ is conservative, since each of the component arrows are conservative.
\end{observation} 
Now, it follows by unwinding definitions that a cocartesian fibration $\pi\colon \Tot_{\cT} \cO^{\otimes} \rightarrow \Span(\FF_{\cT})$ is an $I$-operad if and only if its unstraightening $\Span_I(\FF_{\cT}) \rightarrow \Cat$ is an $I$-symmetric monoidal category.
\cite{Barkan} and \cite{Nardin} thus independently construct an adjunction
\[
  \begin{tikzcd}
    {\Op_{I}} & {\Cat_{I}^{\otimes}.}
    \arrow[""{name=0, anchor=center, inner sep=0}, "\Env_I", curve={height=-12pt}, from=1-1, to=1-2]
    \arrow[""{name=1, anchor=center, inner sep=0}, "U", curve={height=-12pt}, from=1-2, to=1-1]
    \arrow["\dashv"{anchor=center, rotate=-90}, draw=none, from=0, to=1]
  \end{tikzcd}
\]
Now, $\Op_I$ has a terminal object $\cN_{I\infty}^{\otimes}$, and in \cite{EBV} we computed $\Env_I\cN_{I\infty}^{\otimes} \simeq \uFF_{I}^{I-\sqcup}$, i.e. it is the weak indexing system for $I$ with indexed tensor products given by indexed coproducts;
\cite[Prop~4.21]{Barkan} then verifies that the \emph{sliced} left adjoint $\Env_I^{/\uFF_{I}^{I-\sqcup}} \cln \Op_{I} \rightarrow \Cat_{I, /\uFF_{I}^{I-\sqcup}}^{\otimes}$ is fully faithful and identifies its image, i.e. $\Op_{I}$ is a colocalizing subcategory of $I$-symmetric monoidal $\infty$-categories over $\uFF_{I}^{I-\sqcup}$ consisting of the \emph{equifibrations}.

Now, the following construction is very occasionally important.
\begin{construction}
  Given $I$ a weak indexing category and $V \in \cT$, define the $\cT_{/V}$-weak indexing category $I_V$ to consist of those maps over $V$ which lie in $I$.
  Define
  \[
    \Span_I(\FF_{V}) \deq \Span_{I_V} (\FF_V).
  \]
  This is evidently functorial under \emph{unslicing} functors;
  in particular, pullback along $\Span_I(\FF_V) \rightarrow \Span_I(\FF_W)$ yields a functor
  \[
    \Res_V^W\colon \Op_{I_W} \rightarrow \Op_{I_V}.
  \]
  We refer to the associated $\cT$-$\infty$-category as $\uOp_I$.
\end{construction}

\subsubsection{The underlying $\cT$-symmetric sequence}
\begin{definition}
  The \emph{underlying $\cT$-$\infty$-category} $U\cO$ of an $I$-operad $\cO^{\otimes}$ is the straightening of the pullback
  \[\begin{tikzcd}[ampersand replacement=\&]
	{\Tot U\cO} \& {\cO^{\otimes}} \\
	{\cT^{\op}} \& {\Span_I(\FF_{\cT})}
	\arrow[from=1-1, to=1-2]
	\arrow[from=1-1, to=2-1]
	\arrow["\lrcorner"{anchor=center, pos=0.125}, draw=none, from=1-1, to=2-2]
	\arrow[from=1-2, to=2-2]
	\arrow[from=2-1, to=2-2]
\end{tikzcd}\]\
A $\cT$-operad \emph{has at most one color} if each value $U\cO_V$ is either empty or contractible, \emph{has at least one color} if $U\cO_V$ is nonempty for each $V$, has \emph{has one color} if $U\cO \simeq *_{\cT}$.
These occupy full subcategories
\[
  \Op_{I}^{\oc} \subset \Op_{I}^{\leq \oc}, \Op_I^{\geq \oc} \subset \Op_I.\qedhere
\]
\end{definition}

In \cite[\S~2.3]{EBV} we defined an \emph{underlying $\cT$-symmetric sequence} functor and proved the following.
\begin{theorem}[{\cite[Thm~A]{EBV}}]
  The underlying $\cT$-symmetric sequence functor $\sseq\cln \Op^{\leq oc}_{\cT} \rightarrow \Fun(\Tot \uSigma_{\cT}, \cS)$ is monadic;
  in particular, it is conservative.
\end{theorem}
The $V$-objects in $\uSigma_{\cT} \simeq \uFF_{\cT}^{\simeq}$ are finite $V$-sets;
given $S \in \Sigma_V \simeq \FF_V^{\simeq}$, writing $\cO(S)$ for $\sseq \cO^{\otimes}(S)$, we remember this as saying that at-most-one-color $\cT$-operads are identified conservatively by their $S$-ary structure spaces.
Using this, we defined the full subcategory of \emph{$\cT$-$d$-operads} as those with $(d-1)$-truncated structure spaces:
\[
  \Op_{\cT,d} \deq \cbr{\cO^{\otimes} \mid \forall S, \,\, \cO(S) \in \cS_{\leq (d-1)}} \subset \Op_{\cT} 
\]
In \cite[\S~2.5]{EBV}, we verified the following.
\begin{proposition}[{\cite{EBV}}]
  The inclusion $\Op_{\cT,d} \subset \Op_{\cT}$ has a left adjoint $h_d\colon \Op_{\cT} \rightarrow \Op_{\cT,d}$, and given $\cP^{\otimes} \in \Op_{\cT,d}$, the $\infty$-category 
  $\Alg_{\cO}(\cP)$ is a $d$-category;
  moreover, if $\cP^{\otimes} \in \Op_{\cT,0}$, then $\Alg_{\cO}(\cP)$ is either empty or contractible.
  In particular, $\Op_{\cT,d}$ is a $(d+1)$-category and $\Op_{\cT,0}$ is a poset.
\end{proposition}
We call $h_d \cO^{\otimes}$ the \emph{homotopy $\cT$-$d$-operad of $\cO^{\otimes}$}.
We went on to compute the \emph{free $\cO$-algebra monad};
for algebras in a cartesian structure on coefficient systems in a cocomplete cartesian closed $\infty$-category $\cC$, this sends $X \in \CoFr^{\cT} \cC$ to the coefficient system $T_{\cO} X$ with
\[
  \prn{T_{\cO} X}^V \simeq \coprod_{S \in \FF_V} \prn{\Fr_{\cC} \cO(S) \times \prod_{U \in \Orb(S)} X^U}_{h\Aut_V(S)},
\]
where $\Fr_{\cC}\colon \cS \rightarrow \cC$ is the unique symmetric monoidal left adjoint.
In particular, given $S \in \FF_V$, in \cite{EBV} we found a natural splitting $\Fr_{\cC} \cO(S) \oplus J \simeq \prn{T_{\cO} S}^V$.
A multiple-color version of this argument yielded the following.
\begin{proposition}[{\cite[\S~2.4]{EBV}}]\label{Conservativity corollary}\label{Unital conservativity corollary}\label{n-connected prop}
  A map of $\cT$-operads $\varphi\cln \cO^{\otimes} \rightarrow \cP^{\otimes}$ is an $h_d$-equivalence if and only if:
  \begin{enumerate}[label={(\alph*)}]
    \item the underlying $\cT$-functor $U\varphi \cln U\cO \rightarrow U\cP$ is essentially surjective,  and 
    \item the pullback functor $\varphi^* \cln \Alg_{\cP}(\ucS_{\cT, \leq (d-1)}) \rightarrow \Alg_{\cO}\prn{\ucS_{\cT, \leq (d-1)}}$ is an equivalence.
  \end{enumerate}
  In particular, $\varphi$ is an equivalence if and only if it is $U$-essentially surjective and induces an equivalence on algebras in $\ucS_{\cT}$.
\end{proposition}

\subsubsection{Rudiments of weak $\cN_\infty$-operads}
In \cite[\S~2.2]{EBV}, we constructed a family of $\cT$-operads:
\begin{proposition}[{\cite{EBV}}]
  Let $I \subset J \subset \FF_{\cT}$ be pullback-stable subcategories.
  Then, $\Span_I(\FF_{\cT}) \rightarrow \Span_J(\FF_{\cT})$ presents a $J$-operad if and only if $I$ is a weak indexing category.
\end{proposition}
These are called weak $\cN_\infty$-operads;
in the case that $I$ is an indexing category, these are called $\cN_\infty$-operads.
To state their universal property, we defined the \emph{arity support} subcategory
\[
  A\cO \deq \cbr{T \rightarrow S \;\; \middle| \;\; \prod_{U \in \Orb(S)} \cO(T \times_S U) \neq \emptyset } \subset \FF_{\cT}, 
\]
\begin{theorem}[{\cite[\S~2.6]{EBV}}]\label{Arity support reftheorem}
  The arity support of a $\cT$-operad is a weak indexing category, and the associated essential surjection admits a fully faithful right adjoint
  \[
    \begin{tikzcd}[ampersand replacement=\&]
	{\Op_{\cT}} \& {\wIndCat_{\cT}}
	\arrow[""{name=0, anchor=center, inner sep=0}, "A", curve={height=-12pt}, from=1-1, to=1-2]
	\arrow[""{name=1, anchor=center, inner sep=0}, "{\cN_{(-)\infty}}", curve={height=-12pt}, hook', from=1-2, to=1-1]
	\arrow["\dashv"{anchor=center, rotate=-90}, draw=none, from=0, to=1]
\end{tikzcd}
  \]
  The essential image of $\cN_{(-)\infty}$ is spanned by $\cT$-operads $\cO^{\otimes}$ satisfying the following equivalent conditions.
  \begin{enumerate}[label={(\alph*)}]
    \item $\cO^{\otimes}$ is a weak $\cN_\infty$-operad.
    \item $\cO^{\otimes}$ is a $\cT$-0-operad.
    \item The map of $\cT$-operads $\cO^{\otimes} \rightarrow \Comm_{\cT}^{\otimes}$ is a monomorphism.
  \end{enumerate}
\end{theorem}
In particular, this isolates the weak $\cN_\infty$-operads as those possessing a fully faithful unslicing functor
\[
  \Op_{\cT, /\cN_{I\infty}} \hookrightarrow \Op_{\cT}.
\]
This yields functors for change of weak indexing systems, which we use to generically specialize to $I = \cT$.
\begin{proposition}[{\cite{EBV}}]
  Postcomposition along $\Span_I(\FF_{\cT}) \rightarrow \Span_J(\FF_{\cT})$ yields a fully faithful embedding
  \[
    \Op_I \simeq \Op_{J, /\cN_{I\infty}^{\otimes}} \hookrightarrow \Op_J.
  \]
\end{proposition}
We denote the associated push-pull adjunction as
\begin{equation}\label{E bor adjunction}
  \begin{tikzcd}[ampersand replacement=\&]
    {\Op_I} \& {\Op_J;}
    \arrow[""{name=0, anchor=center, inner sep=0}, "{E_I^J}", curve={height=-12pt}, hook', from=1-1, to=1-2]
    \arrow[""{name=1, anchor=center, inner sep=0}, "{\Bor_I^J}", curve={height=-12pt}, from=1-2, to=1-1]
    \arrow["\dashv"{anchor=center, rotate=-90}, draw=none, from=0, to=1]
  \end{tikzcd}
\end{equation}
it follows from \cref{Arity support reftheorem} that $A\Bor_I^J \cO \simeq I \cap A\cO$ and $AE_I^J \cO \simeq A\cO$.
We have several examples:

\begin{example}\label{Triv example}
Let $I^{\triv}$ be the initial one-color weak indexing category. 
The corresponding weak $\cN_\infty$-operad $\triv_{\cT}^{\otimes} \deq \cN_{I^{\triv}\infty}^{\otimes}$ is called the \emph{trivial $\cT$-operad}, and it is characterized by its algebras \cite{EBV,Nardin}
\[
  \uCAlg_{I^{\triv}}(\cO) \simeq U\cO;
\]
in particular, the restriction of the underlying $\cT$-$\infty$-category construction to $I^{\triv}$-operads yields an equivalence $\Op_{I^{\triv}} \xrightarrow\sim \Cat_{\cT}$ and $E_{I^{\triv}}^J$ is compatible with $U$ \cite{Nardin};
that is, \cref{E bor adjunction} takes the form of an adjunction 
\[\begin{tikzcd}[ampersand replacement=\&]
	{\Cat_{\cT}} \& {\Op_{\cT},}
	\arrow[""{name=0, anchor=center, inner sep=0}, "{\triv(-)^{\otimes}}", curve={height=-12pt}, hook', from=1-1, to=1-2]
	\arrow[""{name=1, anchor=center, inner sep=0}, "U", curve={height=-12pt}, from=1-2, to=1-1]
	\arrow["\dashv"{anchor=center, rotate=-90}, draw=none, from=0, to=1]
\end{tikzcd}\]
i.e. given a $\cT$-$\infty$-category $\cC$, we acquire a $\cT$-operad characterized by its algebras \cite{Nardin}
\[
  \Alg_{\triv(\cC)}(\cO) \simeq \Fun_{\cT}(\cC, U\cO);
\]
this formula allows for an alternative construction for $\triv(\cC)$ is as the operadic localization \cite{EBV}
\[
  \triv(\cC)^{\otimes} \simeq L_{\Op_{\cT}}\prn{\cC \rightarrow \cT^{\op} \rightarrow \Span(\FF_{\cT})}.\qedhere
\]
\end{example}
\begin{example}\label{Inflated example}
  Given $\cF \subset \cT$ a family, define 
  \[
    \FF^0_{\cT,\cF, V} \deq \begin{cases}
      \cbr{\emptyset_V,*_V} & V \in \cF \\ 
      \cbr{*_V} & V \not \in \cF.
    \end{cases}
  \]
  Let $I^0_{\cT,\cF}$ be the associated weak indexing category;
  this is the initial one-color weak indexing category with $\upsilon(I) \supset \cF$.
  We define $\EE_{0,\cF}^{\otimes} \deq \cN_{I^0_{\cT,\cF} \infty}^{\otimes}$;
  in particular, we write $\EE_0^{\otimes} \deq \EE_{0,\cT}^{\otimes}$.

  It was shown in \cite[Thm~5.2.10]{Nardin} (and re-shown in \cite[\S~3.3]{EBV}) that there is an equivalence
  \[
    \Gamma^{\cT} \uAlg_{\EE_0}^{\otimes}(\cC) \simeq \prn{\Gamma^{\cT} \cC}^{\otimes}_{1/}.\qedhere
  \]
\end{example}

\begin{example}
  Let $\uFF^\infty$ be the initial indexing system, i.e.
  \[
    \FF^\infty_V \deq \cbr{n \cdot *_V \mid n \in \NN}.
  \]
  Let $I^\infty$ be its indexing category.
  We write $\EE_\infty^{\otimes} \deq \cN_{I^\infty \infty}^{\otimes}$.
  In \cite[\S~3.3]{EBV} we constructed an equivalence
  \[
    \Gamma^{\cT} \uAlg_{\EE_\infty}^{\otimes}(\cC) \simeq \CAlg^{\otimes} \prn{\Gamma^{\cT} \cC}.\qedhere
  \]
\end{example}
\begin{example}
  The terminal $\cT$-operad $\Comm_{\cT}^{\otimes} = \prn{\Span(\FF_{\cT}) = \Span(\FF_{\cT})}$ is the $\cN_\infty$-operad for the terminal indexing category $\FF_{\cT} = \FF_{\cT}$.
\end{example}

\subsubsection{The Boardman-Vogt tensor product}
In \cite[Thm~D]{EBV}, in the case that $\cT$ has a terminal object, we equipped $\Op_{\cT}$ with a closed \emph{Boardman-Vogt tensor product}
\[
  \cO^{\otimes} \obv \cP^{\otimes} \deq L_{\Op}\prn{\cO^{\otimes} \times \cP^{\otimes} \longrightarrow \Span(\FF_{\cT}) \times \Span(\FF_{\cT}) \xrightarrow{\;\; \wedge\;\; } \Span(\FF_{\cT})},
\]
Its internal hom is denoted $\uAlg_{\cO}^{\otimes}(\cP)$;
its underlying $\cT$-$\infty$-category is denoted $\uAlg_{\cO}(\cP)$, and it has values
\[
  \uAlg_{\cO}(\cP)_V \simeq \Alg_{\Res_V^{\cT} \cO}(\Res_V^{\cT} \cP).
\]
We verified several properties in \cite{EBV}; for instance, $\uAlg_{\cP}(\cC)$ is an $I$-symmetric monoidal $\infty$-category when $\cC$ is, functorially for $I$-symmetric monoidal maps in $\cC^{\otimes}$ and $\cT$-operad maps in $\cP^{\otimes}$. 
We interpret $\cO \otimes \cP$-algebras as \emph{homotopy-coherently interchanging pairs of $\cO$-algebra and $\cP$-algebra structures} via the following.
\begin{recollection}\label{1-categorical recollection}
  Suppose $\cC^{\otimes}$ is an $I$-symmetric monoidal 1-category and $\cO^{\otimes},\cP^{\otimes}$ are one-color $\cT$-operads.
  We saw in \cite{EBV} that an $\cO^{\otimes} \obv \cP^{\otimes}$-algebra structure on a $\cT$-object $X \in \Gamma^{\cT} \cC$ is equivalently viewed as a pair of $\cO$-algebra and $\cP$-algebra structures 
  subject to the \emph{interchange relation} that, for all $\mu_S \in \cO(S)$ and $\mu_T \in \cP(T)$, the following diagram commutes.
\[\begin{tikzcd}[ampersand replacement=\&, column sep=tiny]
  {\bigotimes\limits_U^S X^{\otimes \Res_U^V T}_V} \& {X_{V}^{\otimes S \times T}} \& {\bigotimes\limits_W^T X^{\otimes \Res_W^V S}_V} \&\&\&\&\&\&\&\&\&\& {X^{\otimes T}_V} \\
	{X^{\otimes S}_V} \&\&\&\&\&\&\&\&\&\&\&\& {X_V}
	\arrow["\simeq"{description}, draw=none, from=1-1, to=1-2]
	\arrow["{\prn{\Res_U^V \mu_T}}", from=1-1, to=2-1]
	\arrow["\simeq"{description}, draw=none, from=1-2, to=1-3]
	\arrow["{\prn{\Res_W^V \mu_S}}", from=1-3, to=1-13]
	\arrow["{\mu_{T}}", from=1-13, to=2-13]
	\arrow["{\mu_{S}}", from=2-1, to=2-13]
\end{tikzcd}\]
  A morphism of $\cO^{\otimes} \obv \cP^{\otimes}$-algebras is simply a morphism  of $\cT$-objects which is simultaneously an $\cO$-algebra map and a $\cP$-algebra map.
\end{recollection}
The following proposition exhibited a key role played by $\triv_\cT^{\otimes}$.
\begin{proposition}[{\cite[Thm~D.(3)]{EBV}}]\label{Triv algebras prop}
  $\triv_{\cT}^{\otimes}$ is the $\obv$-unit;
  hence there is an equivalence of $\cT$-operads
  \[
    \uAlg_{\triv_{\cT}}(\cO) \simeq \cO^{\otimes}
  \]
\end{proposition}
We also saw that $\obv$ is compatible with the \emph{Mode} (i.e. Day coonvolution or box product) structure.
\begin{proposition}[{\cite[Thm~D.(7)]{EBV}}]\label{Env prop}
  The $\cT$-symmetric monoidal envelope intertwines the mode symmetric monoidal structure on $\Cat_{\cT}^{\otimes}$ with Boardman-Vogt tensor products, i.e.
  \[
    \Env\prn{\cO^{\otimes} \obv \cP^{\otimes}} \simeq \Env \prn{\cO^{\otimes}} \circledast \Env \prn{\cP^{\otimes}}.
  \]
  Furthermore, $\Env\prn{\triv_{\cT}^{\otimes}}$ is the $\circledast$-unit.
\end{proposition}
To use all of these results, for the remainder of \cref{First section,Algebras section} we will make the following assumption.
In \cref{Assumption (b) is attained}, we will establish \cref{Bifunctor assumption} in full generality,
so the results of \cref{First section,Algebras section} will apply for arbitrary $\cT$ after that point.
\begin{assumption}
  We assume one of the following things is true.
  \begin{enumerate}[label={(\alph*)}]
    \item $\cT$ has a terminal object.
    \item \label[assumption]{Bifunctor assumption} $\obv\colon \uOp_{\cT} \times \uOp_{\cT} \rightarrow \uOp_{\cT}$ is a $\cT$-bifunctor whose restriction $\Op_V \times \Op_V \rightarrow \Op_V$ is $\obv$ over $\cT_{/V}$.\qedhere
  \end{enumerate}
\end{assumption}

\subsection{Restriction and arity-borelification}\label{Arity borelification section}
We now expand on $\Res_U^V$ and $\Bor_I^J$.
\subsubsection{Operadic restriction and (co)induction}
Recall from \cite[\S~2.3]{EBV} that the underlying $\cT$-symmetric sequence forms a $\cT$-functor $\usseq\colon \uOp_{\cT}^{\red} \rightarrow \uFun_{\cT}(\uSigma_{\cT},\ucS_{\cT})$;
in particular, restriction of $\uV$-operads lies over restriction of $\uV$-symmetric sequences.
This upgrades \cref{Arity support reftheorem} to an adjunction of $\cT$-$\infty$-categories.
\begin{proposition}\label{Windex T-adjunction corollary}
    $\Res_V^W \cN_{I \infty}^{\otimes} \simeq \cN_{\Res_V^W I \infty}^{\otimes}$;
    more generally, $A \dashv \cN_{(-)\infty}^{\otimes}$ lifts to a $\cT$-adjunction
    \[
        \begin{tikzcd}
    	{\uOp_{\cT}} && {\uwIndSys_{\cT}}
    	\arrow[""{name=0, anchor=center, inner sep=0}, "A" above, curve={height=-15pt}, from=1-1, to=1-3]
    	\arrow[""{name=1, anchor=center, inner sep=0}, "{\cN^{\otimes}_{(-)\infty}}" below, curve={height=-15pt}, hook', from=1-3, to=1-1]
    	\arrow["\dashv"{anchor=center, rotate=-90}, draw=none, from=0, to=1]
        \end{tikzcd}
    \]
\end{proposition}
\begin{proof}
  Restriction compatiblility of the underlying symmetric sequence implies that $\Res_V^W A\cO = A \Res_V^W \cO$, lifting $A$ to a $\cT$-functor $\uOp_{\cT} \rightarrow \uwIndSys_{\cT}$ whose $V$-value is $A:\Op_V \rightarrow \wIndSys_V$.
  The right adjoints $\cN^{\otimes}_{(-)\infty}$ uniquely lift to a right $\cT$-adjoint to $\cN^{\otimes}_{(-)\infty}$ by \cite[Prop~7.3.2.6]{HA}, completing the proposition.
\end{proof}

Since $A$ is a $\cT$-left adjoint, it is compatible with $\cT$-colimits.
Applying this for indexed coproducts, we immediately acquire the following convenient properties of $A$.
\begin{corollary}\label{Induced support corollary}
    If $\cO,\cP$ are $\cT$-operads, then we have
    \[
        A(\cO \sqcup \cP) = A\cO \vee A\cP.
    \]
    If $\cQ$ is a $V$-operad, then we have
    \[
        A \Ind_V^{W} \cQ = \Ind_V^{W} A\cQ.
    \]
\end{corollary}
We may use an analogous argument to that of \cite[Lem~4.1.13]{Barkan} to show that $\uOp_{\cT}$ strongly admits $\cT$-limits;
since the fully faithful $\cT$-functor $\uOp_{\cT} \rightarrow \uCat^{\mathrm{int-cocart}}_{/\Span(\uFF_{\cT})}$ possesses pointwise left adjoints (given by $L_{\Fbrs}$), it possesses a $\cT$-left adjoint;
in particular, we may compute $\cT$-limits of $\cT$-operads in $\uCat^{\mathrm{int-cocart}}_{/\Span(\uFF_{\cT})}$.
Then, an analogous argument using \cite[Prop~2.3.7]{Barkan} constructs $\cT$-limits in $\uCat^{\mathrm{int-cocart}}_{/\Span(\uFF_{\cT})}$ in $\uFun_{\cT}(\Span(\FF_{\cT}), \uCat_{\cT})_{/\uFF_{\cT}^{\cT-\sqcup}}$, which strongly admits $\cT$-limits, as its a slice $\cT$-$\infty$-category of a functor $\cT$-$\infty$-category into a $\cT$-$\infty$-category which strongly admits $\cT$-limits.
In particular, this constructs a right adjoint to $\Res_V^W\colon \Op_W \rightarrow \Op_V$, which we call $\CoInd_V^W$.
\begin{proposition}\label{Truncated prop}
  If $\cO^{\otimes}$ is a $V$-$d$-operad, then $\CoInd_V^W \cO^{\otimes}$ is a $W$-$d$-operad.
\end{proposition}
\begin{proof}
  This follows simply by taking right adjoints within the following diagram
  \begin{equation}\label{Res op diagram}
    \begin{tikzcd}
      {\Op_{W}} && {\Op_{V}} \\
      {\Op_{W,d}} && {\Op_{V,d}}
      \arrow["{\Res_V^W}", from=1-1, to=1-3]
      \arrow[from=1-1, to=2-1]
      \arrow[from=1-3, to=2-3]
      \arrow["{\Res_V^W}", from=2-1, to=2-3]
    \end{tikzcd}
  \end{equation}
  \shortqed
\end{proof}

\begin{corollary}\label{CoInd support corollary}
  There exist equivalences  
  \begin{align*}
    \sseq \CoInd_V^W \cO^{\otimes} &\simeq \CoInd_V^W \sseq \cO^{\otimes};\\
        A\CoInd_V^{W} \cO &= \CoInd_V^{W} A \cO.
  \end{align*}
\end{corollary}
\begin{proof}
    The first statement equivalence by noting that $\Fr \Res_V^{W} = \iota_V^{W*} \Fr$ and taking right adjoints.
    The second follows by taking right adjoints of \cref{Res op diagram} in the case $d = 0$.
\end{proof}

We care about $\CoInd_{V}^{W} \cO^{\otimes}$ because it is a structure borne by \emph{norms of algebras} as follows.
\begin{construction}\label{Norms of algebras construction}
  Let $I$ be a $W$-weak indexing category containing the map $V \rightarrow W$, let $\cP^{\otimes} \rightarrow \CoInd_V^W \cO^{\otimes}$ be a functor of one-object $I$-operads, and let $\cC^{\otimes}$ be a $I$-symmetric monoidal $\infty$-category.
  Then, the adjunct map $\varphi\colon \Res_V^W \cP \rightarrow \cO^{\otimes}$ participates in a commutative diagram of symmetric monoidal functors
  \[
      \begin{tikzcd}
        \Alg_{\cO}(\Res_V^W \cC) \arrow[r, "\varphi^*"] \arrow[d,"U_V"]
          &\Alg_{\Res_V^W \cP}(\Res_V^W \cC) \arrow[r,"N_V^W"] \arrow[d,"U_V"]
          &\Alg_{\cP}(\cC) \arrow[d,"U_W"]\\
          \cC_V \arrow[r,equals]
          & \cC_V \arrow[r,"N_V^W"]
          & \cC_W
      \end{tikzcd}
  \]
  Intuitively, $\CoInd_V^{W} \cO^{\otimes}$ bears the \emph{universal} natural structure on $N_V^W X$ for all $X \in \Alg_{\cO}(\cC)$. 
\end{construction}

\subsubsection{Color and arity  Borelification}%
Let $\cF \subset \cT$ be a $\cT$-family.
There is a terminal $\cF$-colored weak indexing category $\FF_{\cF}$;
we refer to $\FF_{\cF}$-Borelification as \emph{$\cF$-Borelification} and write $\Bor_\cF^{\cT} \deq \Bor_{\FF_{\cF}}^{\cT}$.
Note that
\[
  \Op_{\cF} \simeq \Op_{I_{\cF}}.
\]
Let $\triv_{\cF}^{\otimes} \deq \triv(*_{\cF})^{\otimes}$;
this is the weak $\cN_\infty$-operad for the initial $\cF$-colored weak indexing category $I_{\cF}^{\triv} = \FF_{\cF}^{\simeq} = E_{\cF}^{\cT} I^{\triv}$, and in particular, there is an equivalence
\begin{equation}\label{Triv is E equation}
  \triv^{\otimes}_{\cF} \simeq E_{\cF}^{\cT} \triv_{\cF}^{\otimes}.
\end{equation}
This is our first nontrivial example of a $\obv$-idempotent $\cT$-operad.
\begin{proposition}[Color-borelification]\label{Smashing colocalization corollary}
    Given $\cF \in \Fam_{\cT}$ a $\cT$-family, there is a natural equivalence
    \[
        \Alg_{\triv_{\cF}}(\cO) \simeq \Gamma^{\cF}\cO;
    \]
    hence there is a natural equivalence
    \[
        \triv_{\cF}^{\otimes} \obv \cO^{\otimes} \simeq E_{\cF}^{\cT} \Bor_{\cF}^{\cT} \cO^{\otimes}.
    \]
\end{proposition}
\begin{proof}
    The first statement follows by using \cref{Triv is E equation,Triv algebras prop} to construct equivalences
    \[
        \Alg_{\triv_{\cF}}(\cO) \simeq \Alg_{\triv_{\cF}}(\Bor_{\cF}^{\cT}(\cO)) \simeq \Gamma^{\cF} \cO.
    \]
    The second statement then follows by Yoneda's lemma, noting that
    \begin{align*}
         \Alg_{\triv_{\cF} \otimes \cO}(\cP) 
         &\simeq \Alg_{\triv_{\cF}} \uAlg_{\cO}^{\otimes}(\cP)\\
         &\simeq \Gamma^{\cF} \Alg_{\cO}(\cP)\\
         &\simeq \Alg_{\Bor_{\cF}^{\cT} \cO}(\Bor_{\cF}^{\cT} \cP)\\
         &\simeq \Alg_{E_{\cF}^{\cT} \Bor_{\cF}^{\cT} \cO}(\cP).\qedhere
    \end{align*}
\end{proof}
Given $\cO^{\otimes} \in \Op_{\cT}$, we define the color family $c(\cO) \deq c(A\cO) = \cbr{V \mid \cO_V \neq \emptyset}$.
\begin{remark}
  \cref{Smashing colocalization corollary} exhibits $\mathrm{Im} E_{\cF}^{\cT} = \cbr{\cO^{\otimes} \in \Op^{\cT} \mid c(\cO) \subset \cF}$ as a \emph{$\otimes$-ideal}, i.e. if $c(\cO) \subset \cF$, and $\cP^{\otimes}$ is arbitrary, then $c\prn{\cO \obv \cP} \subset \cF$.
\end{remark}

This is important in part because it reduced $\obv$ computations to the at-least one color case.
\def\spotbv{{\hspace{0pt} \obv \hspace{0pt}}}
\begin{observation}\label{One-color tensor product observation}
    There is a string of natural equivalences
    \begin{align*}
        \cO^{\otimes} 
            \obv           
            \cP^{\otimes} 
        &\simeq \cO^{\otimes} 
            \obv 
            \triv_{c\cO}^{\otimes} 
            \spotbv 
            \triv_{c\cP}^{\otimes} 
            \obv 
            \cP^{\otimes},\\
        &\simeq \cO^{\otimes} 
            \spotbv
            \triv_{c\cO \cap c\cP}^{\otimes} 
            \spotbv 
            \cP^{\otimes},\\
        &\simeq \cO^{\otimes} 
            \obv
            \triv_{c\cO \cap c\cP}^{\otimes} 
            \spotbv
            \triv_{c\cO \cap c\cP}^{\otimes} 
            \obv
            \cP^{\otimes},\\
            &\simeq E_{c\cO \cap c\cP}^{\cT} \Bor_{c\cO \cap c\cP}^{\cT} \prn{\cO^{\otimes}} 
            \spotbv 
            E_{c\cO \cap c \cP}^{\cT} \Bor_{c\cO \cap c\cP}^{\cT} \prn{\cP^{\otimes}},\\
            &\simeq E_{c\cO \cap c\cP}^{\cT}\prn{ \Bor_{c\cO \cap c\cP}^{\cT} \prn{\cO^{\otimes}} 
            \obv
          \Bor_{c\cO \cap c\cP}^{\cT} \prn{\cP^{\otimes}}}.
    \end{align*}
    Moreover, the $c\cO \cap c\cP$-operads $\Bor_{c\cO \cap c\cP}^{\cT} \prn{\cO^{\otimes}}$ and $\Bor_{c\cO \cap c\cP}^{\cT} \prn{\cP^{\otimes}}$ both have at least one color.
\end{observation}
Having done this, we may compute arity-supports of arbitrary tensor products of $\cT$-operads.
\begin{proposition}\label{Support prop}
    Suppose $\cO^{\otimes},\cP^{\otimes}$ are $\cT$-operads.
    Then, 
    \[
      A\prn{\cO \obv \cP} = E_{c\cO \cap c\cP}^{\cT} \Bor_{c\cO \cap c\cP}^{\cT} \prn{A\cO \vee A\cP}.
    \]
\end{proposition}
\begin{proof}
    By \cref{One-color tensor product observation}, we have equivalences 
    \begin{align*}
        A\prn{\cO^{\otimes} \otimes \cP^{\otimes}} 
        &\simeq E_{c\cO \cap c\cP}^{\cT} A\prn{\Bor_{c\cO \cap c\cP}^{\cT} \prn{\cO^{\otimes}} \obv \Bor_{c\cO \cap c\cP}^{\cT} \prn{\cP^{\otimes}}},
    \end{align*}
    so it suffices to prove the proposition in the case that $\cO^{\otimes}$ and $\cP^{\otimes}$ have at least one color.
    
    In this case, first note that there exist maps
    \[
        \cO^{\otimes} \otimes \triv^{\otimes}_{\cT}, \triv_{\cT}^{\otimes} \otimes \cP^{\otimes}  \rightarrow \cO^{\otimes} \otimes \cP^{\otimes},
    \]
    and applying $A$ together with the universal property for joins yields an inequality
    \[
        A\cO \vee A\cP \leq A(\cO \vee \cP).
    \]
    To provide the inequality in the other direction, by \cref{Windex T-adjunction corollary}, in the case of \cref{Bifunctor assumption} we may pass to a restriction and assume that $\cT$ has a terminal object;
    in this case, 
    there exists a composite map
    \[
        \cO^{\otimes} \otimes \cP^{\otimes}  
            \rightarrow 
        \cN_{A\cO \infty}^{\otimes} \otimes \cN_{A\cP\infty}^\otimes 
            \rightarrow
        \cN_{A\cO \vee A\cP \infty}^{\otimes} \otimes \cN_{A\cO \vee A\cP \infty}^\otimes 
            \rightarrow 
        \cN_{A\cO \vee A\cP\infty}^{\otimes},
    \]
    whose last map is presented by the bifunctor 
    \[
      \begin{tikzcd}[ampersand replacement=\&]
        {\Span_I(\FF_{\cT}) \times \Span_I(\FF_{\cT})} \& {\Span_I(\FF_{\cT})} \\
        {\Span(\FF_{\cT}) \times \Span(\FF_{\cT})} \& {\Span(\FF_{\cT});}
        \arrow["\wedge", from=1-1, to=1-2]
        \arrow[from=1-1, to=2-1]
        \arrow[from=1-2, to=2-2]
        \arrow["\wedge", from=2-1, to=2-2]
      \end{tikzcd}
    \]
    here, the top map is defined canonically by the fact that weak indexing categories $I \subset \FF_{\cT}$ are closed under cartesian products \cite{Windex}.
    Applying $A$ to this map yields $A(\cO \vee \cP) \leq A\cO \vee A\cP$, as desired.
\end{proof}

We immediately acquire the following corollary.
\begin{corollary}
  $\Op_I \subset \Op_{\cT}$ is closed under binary tensor products;
  if $I$ has one color, then $\triv_{\cT}^{\otimes} \in \Op_I$.
\end{corollary}

\subsection{(Co)cartesian \tI-symmetric monoidal \tinfty-categories}\label{Cartesian subsection}
Fix $I$ an almost-unital weak indexing system.
In this section, we characterize \emph{cartesian and cocartesian} $I$-symmetric monoidal $\infty$-categories, in part as examples of interest and in part as universal construction.

We defer the minutiae of these to \cref{Cocartesian proof subsubsection}, where we construct $\infty$-categories  
$\cC^{I-\times},\cC^{I-\sqcup}$ over $\Tot \uFF_{\cT,*}$, verifying that $\cC^{I-\times}$ is an $I$-symmetric monoidal $\infty$-category precisely when $\cC$ has $I$-indexed products, that $\cC^{I-\sqcup}$ is always an $I$-operad, and that that $\cC^{I-\sqcup}$ is an $I$-symmetric monoidal $\infty$-category precisely when $\cC$ has $I$-indexed coproducts.
Most of this abuts to unenlightening technicalities about parameterized higher category theory, which we defer to \cref{Cocartesian proof subsubsection}, summarizing the outcomes as they become relevant.
 
\subsubsection{(Co)cartesian rigidity}
Denote by $\Cat_{I}^{I-\sqcup}, \Cat^{I-\times}_{I} \subset \Cat_{\cT}$ the replete subcategories with objects given by $\cT$-$\infty$-categories attaining $I$-indexed coproducts (resp. products) and with morphisms given by $\cT$-functors which preserve $I$-indexed coproducts (products).
In \cref{Cocartesian proof subsubsection}, we prove the following.

\begin{repmaintheorem}{cooltheorem}{cooltheorem}{Cocartesian rigidity main theorem}\label{Cocartesian rigidity main theorem generalized}
  There are fully faithful embeddings $(-)^{I-\sqcup},(-)^{I-\times}$ making the following commute:
  \[\begin{tikzcd}[column sep = large]
	{\Cat_{I}^{I-\sqcup}} & {\Cat_{I}^{\otimes}} & {\Cat_I^{I-\times}} \\
	& {\Cat_\cT}
	\arrow["{(-)^{I-\sqcup}}", hook', from=1-1, to=1-2]
	\arrow["U"', from=1-1, to=2-2]
	\arrow["U"{description}, from=1-2, to=2-2]
	\arrow["{(-)^{I-\times}}"', hook, from=1-3, to=1-2]
	\arrow["U", from=1-3, to=2-2]
\end{tikzcd}\]
The image of $(-)^{I-\sqcup}$ is spanned by the $I$-symmetric monoidal $\infty$-categories whose $I$-admissible indexed tensor functors $\otimes^S\colon \cC_S \rightarrow \cC_V$ are left adjoint to the indexed diagonal $\Delta^S\colon \cC_V \rightarrow \cC_S$ (i.e. whose indexed tensor products are are indexed coproducts), and the image of $(-)^{I-\times}$ is spanned by those whose $I$-admissible indexed tensor functors $\otimes^S$ are right adjoint to $\Delta^S$.
\end{repmaintheorem}
We call $I$-symmetric monoidal $\infty$-categories of the form $\cC^{I-\sqcup}$ \emph{cocartesian}, and $\cC^{I-\times}$ \emph{cartesian}.
\begin{philremark}
  In higher category theory, a fundamental rigidity result is that of \emph{adjoints};
  by \httprop{5.2.6.2}, the full subcategory of $\Fun(\cD,\cC)^{\op}$ spanned by functors right adoint to a fixed functor $L\colon \cC \rightarrow \cD$ is contractible.
  Moreover, adjointness itself as a property requires only finite data to test per-object, separately (c.f. \httpropprop{5.2.2.9}{5.2.2.12}).

  It is this rigidity which we are leveraging in \cref{Cocartesian rigidity main theorem generalized}:
  in essence, the coherent data witnessing $I$-symmetric monoidality of a functor $\cC^{I-\times} \rightarrow \cD^{I-\times}$ is constructed (up to contractible ambiguity) by the \emph{property} that the underlying $\cT$-functor $\cC \rightarrow \cD$ is compatible with the adjunctions $\Delta^S \dashv \bigotimes^S$.
\end{philremark}

Many similar definitions have been made in the literature.
Luckily, they all agree.
\begin{remark}
  \cite{Nardin} constructed a pair of structures $\cC^{\coprod}, \cC^{\prod}$ which, after unwinding definitions, satisfy the conditions of \cref{Cocartesian rigidity main theorem generalized} in the case $I = \cT$.
  In particular, there are unique $I$-symmetric monoidal equivalences $\Bor_I^{\cT} \cC^{\coprod} \simeq \cC^{I-\sqcup}$ and $\Bor_I^{\cT} \cC^{\prod} \simeq \cC^{I-\times}$ lying over the identity whenever $\cC$ admits finite indexed (co)products.

  Moreover, \cite{Cnossen_tambara} introduced another structure, specifically a \emph{cartesian $I$-symmetric monoidal structure} on $\uFF_{\cT}$, and conjectured it to be equivalent to Nardin-Shah's construction;
  this conjecture was recently verified in \cite{Cnossen_unfurling} by verifying universality of the unfurling construction used in Nardin-Shah's construction.
  We acquire an independent proof of this result without serious $(\infty,2)$-category theory:
  Cnossen-Haugseng-Lenz-Linsken's construction satisfies the condition of \cref{Cocartesian rigidity main theorem generalized}, so there is a unique pair of $I$-symmetric monoidal equivalences $\Bor_I^{\cT} \uFF_{\cT}^{\prod} \simeq \uFF_{\cT}^{I-\times} \simeq \uFF_{\cT,\times}$
  lying over the identity.

  Moreover, after drafts of this article were made public, \llp{Def.}{2.15} constructed another structure, which specializes to a \emph{cocartesian $I$-symmetric monoidal structure} $\cC^{I-\coprod}$;
  by \llp{Lem.}{2.16}, in the case that $\cC$ has $I$-indexed coproducts, $\cC^{I-\coprod}$ satisfies the conditions of \cref{Cocartesian rigidity main theorem generalized}, so there is a unique $I$-symmetric monoidal equivalence $\cC^{I-\coprod} \simeq \cC^{I-\sqcup}$ lying over the identity.
\end{remark}
We prove the following in \cref{Cocartesian proof subsubsection} as a precursor to \cref{Cocartesian rigidity main theorem generalized}, though it also follows from it.
\begin{proposition}\label{vop cart observation}
    There is a unique equivalence $\prn{\cC^{I-\times}}^{\vop} \simeq \prn{\cC^{\vop}}^{I-\sqcup}$ lying over the identity.
\end{proposition} 

Before characterizing the algebras in $\cC^{I-\sqcup}$ and $\cC^{I-\times}$, we point out that they are often presentable.
\begin{proposition}\label{Cartesian presentability prop}
  Suppose $\cC$ is a presentable $\infty$-category with $I$-indexed products and coproducts. 
  \begin{enumerate}
    \item $\uCoFr^{\cT} \cC^{I-\sqcup}$ is presentably $I$-symmetric monoidal.
    \item If finite products in $\cC$ commute with colimits separately in each variable (i.e. it is Cartesian closed), then $\uCoFr^{\cT}\cC^{I-\times}$ is presentably $I$-symmetric monoidal,
  \end{enumerate}
\end{proposition}
\begin{proof}
  It follows from Hilman's characterization of parameterized presentability \cite[Thm.~\href{https://arxiv.org/pdf/2202.02594v2\#subsubsection.6.1.2}{6.1.2}]{Hilman} that $\uCoFr^{\cT}$ is presentable, so we're tasked with proving that the $\cT$-symmetric monoidal structures are distributive.
  The first case is just commutativity of colimits with colimits, and the second is \ns{Prop.}{3.2.5}.
\end{proof}

\subsubsection{$\cO$-monoids}
We will identify algebras in $\cC^{I-\times}$ with the following.
\begin{definition}
  Fix $\cO^{\otimes}$ an $I$-operad and $\cC$ a $\cT$-$\infty$-category.
  Then, an \emph{$\cO$-monoid in $\cC$} is a $\cT$-functor $M\colon \Tot_{\cT}\cO^{\otimes} \rightarrow \cC$ satisfying the condition that, for each orbit $V \in \cT$, each finite $V$-set $S \in \FF_V$, and each $S$-tuple $X = (X_U) \in \cO_S$, the canonical maps $M(X) \rightarrow \CoInd_U^V M(X_U)$ realize $M(X)$ as the indexed product
  \[
    M(X) \simeq \prod_U^S M(X_U).\qedhere
  \]
\end{definition}
Indeed, we prove the following equivariant lift of \ha{Prop.}{2.4.2.5} as \cref{Cartesian algebras generalized}. 
\begin{proposition}\label{Cartesian algebras}
  Given $\cO^{\otimes}$ an $I$-operad and $\cC$ a $\cT$-$\infty$-category with $I$-indexed products, the forgetful functor
  \[
    \Alg_{\cO}\prn{\cC^{I-\times}} \longrightarrow \Fun_{\cT}\prn{\Tot_{\cT} \cO^{\otimes},\cC}
  \]  
  is fully faithful with image spanned by the $\cO$-monoids.
\end{proposition}
\begin{corollary}\label{Tot tot corollary}
  Given $\cO^{\otimes}$ an $I$-operad and $\cD$ an $\infty$-category with finite products, the forgetful functor
  \[
    \Alg_{\cO}\prn{\uCoFr^G(\cD)^{I-\times}} \longrightarrow \Fun\prn{\Tot\Tot_{\cT} \cO^{\otimes}, \cD}
  \]
  is fully faithful with image spanned by $\Seg_{\Tot \Tot_{\cT} \cO^{\otimes}}(\cD)$.
\end{corollary}
\begin{proof}
  After \cref{Cartesian algebras}, it suffices to characterize the image of $\cO$-monoids under the equivalence 
  \[
    \Fun(\Tot \Tot_{\cT} \cO^{\otimes}, \cD) \simeq \Fun_{\cT}(\Tot_{\cT} \cO^{\otimes}, \uCoFr_G(\cD)).
  \]
  By \cite[Ex~1.17]{Nardin_thesis}, given a finite $V$-set $S \in \FF_V$ and writing $\Tot S \simeq \coprod_{U \in \Orb(S)} \cT_{/U}$ for the total $\infty$-category of the associated $V$-category, the above identification turns $S$-indexed products into right Kan extensions:
  \[
    \begin{tikzcd}[ampersand replacement=\&, row sep=small]
      {\Fun_{\cT}(S,\uCoFr^{\cT}(\cD))} \& {\CoFr^{\cT}(\cD)} \\
	{\Fun\prn{\Tot S, \cD}} \& {\Fun(\cT^{\op},\cD)}
	\arrow[from=1-1, to=1-2, "\prod^S"]
	\arrow["\simeq"{marking, allow upside down}, draw=none, from=1-1, to=2-1]
	\arrow["\simeq"{marking, allow upside down}, draw=none, from=1-2, to=2-2]
  \arrow[from=2-1, to=2-2, "\mathrm{RKE}"]
    \end{tikzcd}
  \]
  Thus the image of $\Mon_{\cO}(\uCoFr^{\cT} \cD)$ is those functors $\Tot \Tot_{\cT} \cO^{\otimes} \rightarrow \cD$ whose image of an object $((X_U),S) \in \Tot \Tot_{\cT} \cO^{\otimes}$ is right Kan extended along elementary maps, which is exactly the relevant Segal condition.
\end{proof}

\begin{corollary}\label{CMon cor 1}
  Given $\cO^{\otimes}$ an $I$-operad and $\cD$ an $\infty$-category with finite products, the forgetful functor
  \[
    \Alg_{\cO}(\uCoFr^G(\cD)^{I-\times}) \longrightarrow \Fun(\Tot \cO^{\otimes}, \cD)
  \]
  is fully faithful with image spanned by $\Seg_{\Tot \cO^{\otimes}}(\cD)$.
\end{corollary}
\begin{proof}
  Apply \cref{Tot tot corollary} and the equivalence $\Seg_{\Tot \Tot_{\cT} \cO^{\otimes}}(\cC) \simeq \Seg_{\Tot \cO^{\otimes}}(\cC)$ constructed in \cite[\S~\href{https://arxiv.org/pdf/2501.02129v1\#appendix.A}{A}]{EBV}.
\end{proof}

We finally identify \cref{Lawvere perspective,Semiadditive perspective,Operad perspective} from the introduction.
\begin{corollary}[``$\CMon =\CAlg$'']\label{CMon is CAlg corollary}
  There is a canonical equivalence $\uCMon_I(\cC) \simeq \uCAlg_I(\cC^{I-\times})$ over $\cC$.
\end{corollary}
\begin{proof}
  Our proof is similar to that of \nar{Thm.}{6.5};
  there is a pullback square over $\cC$
  \[\begin{tikzcd}[ampersand replacement=\&, column sep=tiny]
	{\CMon_I(\cC)} \&\& {\CAlg_I(\cC^{I-\times})} \& {\Fun^{I-\times}(\Span(\uFF_{\cT}), \cC)} \\
	{\Fun_{\cT}(\cC^{\op}, \uCMon_I(\ucS_{\cT}))} \&\& {\Fun_{\cT}\prn{\cC^{\op}, \uCAlg_I(\cC^{I-\times})}} \& {\Fun_{\cT}(\cC^{\op}, \uFun^{I -\times}(\Span(\uFF_{\cT}), \ucS_{\cT}))}
	\arrow[from=1-1, to=1-3]
	\arrow[from=1-1, to=2-3, "\lrcorner" very near start, phantom]
	\arrow[from=1-1, to=2-1]
	\arrow["\simeq"{description}, draw=none, from=1-3, to=1-4]
	\arrow[from=1-3, to=2-3]
	\arrow[from=1-4, to=2-4]
	\arrow[from=2-1, to=2-3]
	\arrow["\simeq"{description}, draw=none, from=2-3, to=2-4]
\end{tikzcd}\]
  so it suffices to prove this in the case $\cC = \ucS_{\cT}$.
  There, we simply compose equivalences as follows
  \[\begin{tikzcd}[ampersand replacement=\&]
	{\CMon_I(\ucS_{\cT})} \& {\CMon_I(\cS)} \& {\CAlg_I(\ucS_{\cT}^{I-\times})}
	\arrow["{\ref{CMon in coeff}}", from=1-1, to=1-2]
	\arrow["{\ref{CMon cor 1}}", from=1-2, to=1-3]
\end{tikzcd}\qedhere\]
\end{proof}
\begin{remark}
  As with much of the rest of this subsection, \cref{CMon is CAlg corollary} possesses an alternative strategy where both are shown to furnish the $I$-semiadditive closure, the latter using \cite[Thm~\href{https://arxiv.org/pdf/2403.07676v3\#introthm.2}{B}]{Cnossen_semiadditive}.
  The above argument was chosen for brevity, as its requisite parts are also needed elsewhere.
\end{remark}
\begin{remark}
  In the case $\cC \simeq \ucS_G$, and $I$ is an indexing category, the analogous result was recently proved in \cite{Marc} for a Dwyer-Kan localization of algebras over the corresponding \emph{graph $G$-operads}.
  To the knowledge of the author, this is one of the first concrete higher-categorical indications that the genuine operadic nerve of \cite{Bonventre-nerve} may induce equivalences between $\infty$-categories of algebras.
\end{remark}

\subsubsection{\texorpdfstring{$\cF$}{F}-unitality} \label{F-unitality subsection}
We now study $I$-operadic unitality, beginning with the following definition.
\begin{definition}
  We say that an $I$-operad $\cO^{\otimes}$ is \emph{unital} if $\cO(\emptyset_V;C) = *$ for all $C \in \cO_V$ with $V \in \upsilon(I)$, and \emph{reduced} if also $\cO(*_V) = *$ for all $V \in c(I)$.
  More generally if $\cF \subset \upsilon(I)$ is a family, we say that $\cO^{\otimes}$ is \emph{$\cF$-unital} if $\cO(\emptyset_V;C) \simeq *$ for all $C \in \cO_V$ and $V \in \cF$ and \emph{$\cF$-reduced} if also $\cO(*_V) = *$ for all $V \in \cF$;
  equivalently, $\cO^{\otimes}$ is $\cF$-unital (resp. $\cF$-reduced) if and only if $\Bor_{I \cap \FF_{\cF}}^{I} \cO^{\otimes}$ is unital (reduced).
\end{definition}
An under-appreciated case of unitality is the (equivariantly) \emph{symmetric monoidal case}.
\begin{observation}\label{Unital SMC observation}
  If $\cC^{\otimes}$ is an $I$-symmetric monoidal $\infty$-category with unit $\upsilon(I)$-object $1_\bullet$ and $X \in \cC_V$, then the Segal conditon for multimorphisms constructs an equivalence 
  \[
    \Map_{\cO^{\otimes}}(\emptyset_V, X) \simeq \Map_{\cC_V}(1_V,X);
  \]
  hence $\cC^{\otimes}$ is $\cF$-unital if and only if $1_\bullet \in \Gamma^{\cF} \cC$ is initial.
  In particular, if $\Bor_{I \cap \FF_{\cF}}^{I} \cC^{\otimes}$ is cartesian, then $\cC^{\otimes}$ is $\cF$-unital if and only if $U \cC$ is $\cF$-pointed.
\end{observation}
In fact, we may reverse this, characterizing unitality in terms of the $I$-symmetric monoidal envelope.
\begin{observation}
  We may identify an object $(T,\bC) \in \Env_I(\cO)_V$ with an $I$-admissible finite $V$-set $T$ and an $S$-color $\bC$;
  the space of maps $(S,\bB) \rightarrow (T, \bC)$ lying over a fixed map $\psi\colon T \rightarrow S$ is precisely the multimorphism space $\Mul_{\cO}^{\psi}(\bC;\bB)$.
  In particular, letting $\emptyset_V$ denote the unique $\emptyset_V$-color of $\cO$ for $V \in \cF$, we find that
  \[
    \Map_{\Env_I(\cO)}(\emptyset_V, (S,\bC)) \simeq \prod_{U \in \Orb(S)} \cO(\emptyset_V; C_U).
  \]
  Products of spaces are contractible if and only if their factors are contractible;
  hence $\cO^{\otimes}$ is $\cF$-unital if and only if $\Env_I(\cO)^{\otimes}$ is $\cF$-unital.
\end{observation}

We can identify this via an algebraic mapping-in property as follows.
\begin{lemma}\label{F-unitality lemma}
  If $\cC^{\otimes}$ is an $I$-symmetric monoidal $\infty$-category, then $\cC^{\otimes}$ is $\cF$-unital if and only if the forgetful $\cT$-functor $U\colon \uAlg_{\EE_{0,\cF}}(\cC) \rightarrow \cC$ is an equivalence.
\end{lemma}
\begin{proof}
  The forward implication follows from the computation \cref{Kan extension lemma} in the case $I^0_{\cF}$, so assume $U$ is an equivalence.
  Then, for all $V \in \cF$, $\cC_V^{1_V/} \simeq \Alg_{\EE_{0,\cF}}(\cC)_V \rightarrow \cC$ is an equivalence, so $1_V \in \cC_V$ is initial.
  Thus \cref{Unital SMC observation} implies the lemma.
\end{proof}
We can replace $\EE_{0,\cF}^{\otimes}$ with an \emph{arbitrary} $\cF$-unital $I$-operad, retaining the above property.
\begin{lemma}[{Incomplete \ns{Thm.}{5.2.11}}]\label{Algebras are unital lemma}
  If $\cO^{\otimes}$ is an $\cF$-unital $I$-operad and $\cC^{\otimes}$ is an $I$-symmetric monoidal $\infty$-category, then $\uAlg_{\cO}(\cC)$ is $\cF$-unital.
\end{lemma}
\begin{proof}
  Using the same trick as \cref{F-unitality lemma}, we prove this when $\cF = \upsilon(I) = \cT$.
  Then, in light of \cref{Unital SMC observation}, this is simply \ns{Thm.}{5.2.11}.
\end{proof}

\subsubsection{$\cO$-comonoids, indexed semiadditivity}
We prove the following fundamental fact as \cref{Kan extension lemma}.
\begin{proposition}\label{Cocartesian algebras computation}
  Let $\cC$ be a $\cT$-$\infty$-category and $\cO^{\otimes}$ a unital $I$-operad.
  The forgetful functor is an equivalence
  \[
    \uAlg_{\cO}\prn{\cC^{I-\sqcup}} \xrightarrow{\;\;\;\; \sim \;\;\;\;} \uFun_{\cT}(U\cO, \cC).
  \]
\end{proposition}

We use this to construct another recognition result for cocartesian $I$-symmetric monoidal $\infty$-categories.
\begin{construction}\label{Wirthmuller construction}
  Let $\cC^{\otimes}$ be an $I$-symmetric monoidal $\infty$-category satisfying the property that $1_\bullet \in \Gamma^{\upsilon(I)} \cC$ is initial and let $(X_U) \in \cC_S$ be an $S$-tuple for some $S \in \uFF_I$.
  The \emph{$\otimes$-Wirthm\"uller map for $(X_U)$} is the map
  \[
    W_{S,(X_U)}\colon \coprod_U^S X_U \longrightarrow \bigotimes_U^S X_U
  \]
  classified by the summand maps
  \[\begin{tikzcd}[ampersand replacement=\&]
{ X_W } \& {\Res_U^V \bigotimes_W^S X_W.} \\
{X_W \otimes \bigotimes_{U'}^{\Res_U^V S - W} 1_{U'}} \& {X_W \otimes \bigotimes_{U'}^{\Res_U^V S - W} \Res_{U'}^{o(U')} X_{o(U')}}
	\arrow["{W_{S, (X_U), U'}}", from=1-1, to=1-2]
	\arrow["\simeq"{marking, allow upside down}, draw=none, from=1-1, to=2-1]
	\arrow["\simeq"{marking, allow upside down}, draw=none, from=1-2, to=2-2]
	\arrow["(\id;!)", from=2-1, to=2-2]
\end{tikzcd}\]
where $o(U') \in \Orb(S)$ is the orbit whose restriction contains $U'$.
$1_U$ exists and is initial by almost-unitality.
Dually, if $\cC^{\otimes}$ is an $I$-symmetric monoidal $\infty$-category such that $1_\bullet \in \Gamma^{\upsilon(I)} \cC$ the \emph{$\otimes$-co-Wirthm\"uller map} $W^{co}_{S,(X_U)}\colon \bigotimes_U^S X_U \rightarrow \prod_U^S X_U$ is the Wirtm\"uller map for $(X_U)_S$ in the fiberwise opposite $\cC^{\otimes, \vop}$.
\end{construction}

\begin{lemma}\label{Wirthmuller cocartesianness}
  Let $\cC^{\otimes}$ be an $I$-symmetric monodial $\infty$-category.
  \begin{enumerate}
    \item $\cC^{\otimes}$ is cocartesian if and only if $1_\bullet \in \Gamma^{\upsilon(I)} \cC$ is initial and $W_{S,(X_U)}$ is an equivalence for all $S \in \uFF_{I}$ and $(X_U) \in \cC_S$;
    \item $\cC^{\otimes}$ is cartesian if and only if $1_\bullet \in \Gamma^{\upsilon(I)} \cC$ is terminal and $W_{S, (X_U)}^{co}$ is an equivalence for all $S \in \uFF_I$ and $(X_U) \in \cC_S$.
\end{enumerate}
\end{lemma}
\begin{proof}
  (1) is \cref{Appendix wirthmuller}, and part (2) follows directly from part (1) under \cref{vop cart observation}.
\end{proof}

\begin{observation}\label{Wirthmuller is norm obs}
  When $\cC^{\otimes}$ is cartesian, the assumption that $1_\bullet$ is initial is precisely the assumption that $\cC$ is $\upsilon(I)$-pointed;
  moreover, unwinding definitions, $W_{S, (X_U)}$ matches the Wirthm\"uller map of \cref{Wirth defn}.
\end{observation}

Finally, we find the indexed semiadditivity of \cite{Nardin-Stable,Cnossen_semiadditive} within equivariant higher algebra.
\begin{corollary}[{Equivariant \cite[Prop.~\href{https://arxiv.org/pdf/1305.4550v1\#prop.2.3}{2.3}]{Gepner}}]\label{Semiadditivity is cocartesianness}
  Suppose $\cC$ is a $\cT$-$\infty$-category with $I$-indexed products.
  Then, the following conditions are equivalent.
  \begin{enumerate}[label={(\alph*)}]
    \item $\cC$ is $I$-semiadditive.
    \item There exists an $I$-symmetric monoidal equivalence $\cC^{I-\times} \simeq \cC^{I-\sqcup}$ lying over the identity.
    \item The forgetful $\cT$-functor $\uCMon_I(\cC) \rightarrow \cC$ is an equivalence.
  \end{enumerate}
\end{corollary}
\begin{proof}
  We have proved this redundantly.
  (a) $\iff$ (b) is \cref{Wirthmuller cocartesianness,Wirthmuller is norm obs}.
  (b) $\implies$ (c) is \cref{CMon is CAlg corollary,Cocartesian algebras computation}.
  (c) $\iff$ (a) is \cllsemicor{7.8}.
\end{proof}
 
\subsubsection{Pointwise indexed tensor products of $\cO$-monoids}
Last, we characterize $\uAlg^{\otimes}_{\cO}(\cC^{I-\times})$ categorically.
\begin{lemma}\label{Conservative product preserving symmetric monoidal funcor}
  Fix $F\colon \cD^{\otimes} \rightarrow \cC^{\otimes}$ an $I$-symmetric monoidal functor.
  \begin{enumerate}
    \item If $F$ is $I$-coproduct preserving then there is a homotopy $F W_{S, (X_U)} \sim W_{S, F(X_U)}$.
    \item If $F$ is $I$-product preserving then there is a homotopy $F W^{co}_{S, (X_U)} \sim W^{co}_{S, F(X_U)}$.
    \item If $F$ is $I$-coproduct-preserving and conservative and $\cC^{\otimes}$ is cocartesian then $\cD^{\otimes}$ is cocartesian.
    \item If $F$ is $I$-product preserving and conservative and $\cC^{\otimes}$ is cartesian then $\cD^{\otimes}$ is cartesian.
    \item If $F$ is a fiberwise-monadic right $\cT$-adjoint and $\cC^{\otimes}$ is cartesian then $\cD^{\otimes}$ is cartesian.
  \end{enumerate}
\end{lemma}
\begin{proof}
  For (1), since $F$ is $I$-coproduct preserving, we're tasked with constructing a homotopy making the following diagram commute
\[\begin{tikzcd}[ampersand replacement=\&]
	{F\prn{X_U \otimes \bigotimes_W^{\Res_U^V S - U} 1_W}} \& {F\prn{X_U \otimes \bigotimes_W^{\Res_U^V S - U} X_W}} \\
	{FX_U \otimes \bigotimes_W^{\Res_U^V S - U} 1_U} \& {FX_U \otimes \bigotimes_W^{\Res_U^V S - U} FX_W}
	\arrow["{(\id;!)}", from=1-1, to=1-2]
	\arrow["\simeq"{marking, allow upside down}, draw=none, from=1-1, to=2-1]
	\arrow["\simeq"{marking, allow upside down}, draw=none, from=1-2, to=2-2]
	\arrow["{(\id;!)}", from=2-1, to=2-2]
\end{tikzcd}\]
  In fact, there is a contractible space of such choices.
  (2) follows by applying fiberwise-opposites.

  For (3), applying (1) and \cref{Wirthmuller cocartesianness} shows that $FW_{S, (X_U)}$ is an equivalence for all $S \in \uFF_I$ and $(X_U) \in \cD_S$, so conservativity implies that $W_{S,(X_U)}$ is an equivalence and \cref{Wirthmuller cocartesianness} concludes that $\cD^{\otimes}$ is cocartesian.
  (4) follows similarly from (2).
  (5) is directly implied by (4), since fiberwise monadicity implies conservativity and right $\cT$-adjoints are $I$-product preserving.
\end{proof}

Applying (5) of \cref{Conservative product preserving symmetric monoidal funcor} to $U\colon \uAlg_{\cO}^{\otimes}\prn{\cC^{I-\times}} \rightarrow \cC^{I-\times}$ immediately yields the following.
\begin{corollary}\label{Algebras are cartesian corollary}
  If $\cC$ has $I$-indexed products, then $\uAlg_{\cO}^{\otimes}\prn{\cC^{I-\times}}$ is a cartesian $I$-symmetric monoidal $\infty$-category.
\end{corollary}

\section{\tI-commutative algebras}\label{Algebras section}
\begin{philremark}\label{Philosophical remark}
  On one hand, it follows from \cref{Conservativity corollary} that $I$-operads are determined conservatively by their theories of \emph{algebras in $I$-symmetric monoidal categories};
  indeed, it suffices to characterize their algebras in the universal case $\ucS_{\cT}^{I-\times}$.

  On the other hand, the right adjoint $\Cat_{I}^{\otimes} \rightarrow \Op_{I}$ is full on cores, since automorphisms in the slice category $\Cat_{/\Span_I(\FF_{\cT})}$ automatically preserve cocartesian morphisms.
    Hence the associated map of spaces
    \[ 
        \begin{tikzcd}
            \Cat_{I}^{\otimes, \simeq}
                \arrow[r] \arrow[d,"\in"ulabl, phantom]
            & \Op_{I}^{\simeq}
                \arrow[r]
            & \Fun(\Op_{I},\Cat)^{\simeq}
                \arrow[d,"\in"ulabl, phantom]\\
            \cC^{\otimes} 
                \arrow[rr,mapsto]                
            && \Alg_{(-)}(\cC)
        \end{tikzcd}
    \]
    is a summand inclusion.
    That is, an $I$-symmetric monoidal $\infty$-category is determined (functorially on equivalences) by its categories of $\cO$-algebras for each $I$-operad $\cO$.
\end{philremark}

Following along these lines and using \cref{Cartesian algebras}, we will generally characterize algebraic theories in \emph{arbitrary} settings by reducing to the universal case of $\ucS_{\cT}^{\cT-\times}$, which we study using category theoretic means. 
Indeed, in \cref{Cocartesian subsection,Smashing localization subsection} we use this to bootstrap $I$-semiadditivity of $\uCMon_I(\cC)$ to $I$-cocartesianness of $\uCAlg_I^{\otimes}(\cC)$ for $\cC^{\otimes}$ an arbitrary $I$-operad.
Using work from \cref{Cocartesian proof subsubsection}, we use this to conclude generalizations of \cref{Cocartesian main theorem,Absorption main corollary}, answering \cref{q: absorption to I,q: absorption to O}.

We take this to its logical extreme in \cref{Smashing localization subsection}, using this to completely characterize the smashing localizations associated with $\otimes$-idempotent weak $\cN_\infty$-operads.
As promised in the introduction, we use this classification to prove a generalization of \cref{Main theorem}, answering \cref{q: BH conjecture}.
Following this, in \cref{EH failure subsection} we show that our results are sharp; 
if $I$ is not almost essentially unital, then $\cN_{I \infty}^{\otimes} \obv \cN_{I \infty}^{\otimes}$ fails to be connected, so $\cN_{I \infty}^{\otimes}$ is (abstractly) idempotent under $\obv$ if and only if $I$ is almost essentially unital.

\subsection{Indexed tensor products of \tI-commutative algebras}\label{Cocartesian subsection}
Fix $I$ an almost-unital weak indexing system.
In \cref{Cocartesian algebras computation}, we showed that every object in a cocartesian $I$-symmetric monoidal structure bears a canonical $I$-commutative algebra algebra structure, i.e. $\uCAlg_{I}(\cC) \rightarrow \cC$ is an equivalence.
In this subsection, we demonstrate the converse, i.e. we show the following.
\begin{theorem}[Indexed tensor products of $I$-commutative algebras]\label{Cocartesian algebras theorem}
  The following are equivalent for $\cC^{\otimes} \in \Cat_I^{\otimes}$.
  \begin{enumerate}[label={(\alph*)}]
    \item $\cC^{\otimes}$ is cocartesian.
      \label{Cocartesian rigidity condition}
    \item For all unital $I$-operads $\cO^{\otimes}$, the forgetful functor $\Alg_{\cO}(\cC) \rightarrow \Fun_{\cT}(U\cO,\cC)$ is an equivalence.\label{Vanishing algs condition}
    \item The forgetful $\cT$-functor $\uCAlg_I(\cC) \rightarrow \cC$ is an equivalence.\label{Vanishing calgs condition}
    \item There exists an $I$-symmetric monoidal $\infty$-category $\cD^{\otimes}$ and an equivalence $\uCAlg^{\otimes}_I(\cD) \xrightarrow{\;\;\sim\;\;} \cC^{\otimes}$. \label{Abstractly calg condition}
  \end{enumerate}
\end{theorem}

The implications \ref{Cocartesian rigidity condition} $\implies$ \ref{Vanishing algs condition}, \ref{Vanishing calgs condition} are simply \cref{Cocartesian algebras computation}. 
For the implication \ref{Vanishing algs condition} $\implies$ \ref{Cocartesian rigidity condition}, note that \cref{Algebras are unital lemma} states that $\cC^{\otimes}$ is unital; 
hence Yoneda's lemma applied to $\Op_I^{\uni}$ constructs an $I$-operad equivalence $\cC^{\otimes} \simeq \cC^{I-\sqcup}$, which is an $I$-symmetric monoidal equivalence by \cref{Philosophical remark}.
The implication \ref{Vanishing calgs condition} $\implies$ \ref{Abstractly calg condition} follows by neglect of assumptions.
To summarize, we've arrived at the implications
\begin{equation}\label{Cocartesian implications equation}
  \begin{tikzcd}[ampersand replacement=\&, row sep=small]
    \& {\text{\ref{Abstractly calg condition}}} \\
	{\text{\ref{Vanishing algs condition}}} \& {\text{\ref{Cocartesian rigidity condition}}} \& {\text{\ref{Vanishing calgs condition}}}
	\arrow[Rightarrow, 2tail reversed, from=2-1, to=2-2]
	\arrow[Rightarrow, from=2-2, to=2-3]
	\arrow[Rightarrow, from=2-3, to=1-2]
  \end{tikzcd}
\end{equation}
Our workhorse lemma for closing the gap is the following.
\begin{lemma}\label{Cocartesian rigidity workhorse lemma}
  The following are equivalent for $\cP^{\otimes} \in \Op_{\cT}$:
  \begin{enumerate}[label={(\alph*)}]\setcounter{enumi}{4}
    \item The $\cT$-$\infty$-category $\uAlg_{\cP}(\ucS_{\cT})$ is $I$-semiadditive.\label{Mon semiadditive condition}
    \item For all $\cO^{\otimes} \in \Op_{I}^{\uni}$, the forgetful functor
      \[
        \Alg_{\cO \otimes \cP}(\ucS_{\cT}) \simeq \Alg_{\cO} \uAlg^{\otimes}_{\cP}(\ucS_{\cT}) \longrightarrow \Alg_{\triv(\U\cO)} \uAlg_{\cP}^{\otimes}(\ucS_{\cT}) \simeq  \Fun_{\cT}(U\cO,\uAlg_{\cP}(\ucS_{\cT}))
      \]
      is an equivalence.\label{Mon vanishing algs condition}
    \item For all $\cO^{\otimes} \in \Op_{I}^{\uni}$, the map $\triv\prn{\cO}^{\otimes} \obv \cP^{\otimes} \rightarrow \cO^{\otimes} \obv \cP^{\otimes}$ is an equivalence.\label{Tensor equivalence condition}
    \item 
      \label{Alg vanishing algs condition}
      For all $\cO^{\otimes} \in \Op_I^{\uni}$ and $\cQ^{\otimes} \in \Op_I$, the forgetful $\cT$-operad map
      \[
        \uAlg_{\cO \otimes \cP}^{\otimes}(\cQ) \simeq \Alg_{\cO} \uAlg_{\cP}^{\otimes}(\cQ) \longrightarrow \uAlg_{\triv(U\cO)}^{\otimes} \uAlg_{\cP}^{\otimes}(\cQ)
      \]
      is an equivalence
  \end{enumerate}
\end{lemma}
\begin{proof}
  Since \cref{Algebras are cartesian corollary} shows that $\Bor_I^{\cT} \uAlg^{\otimes}_{\cO}(\ucS_{\cT})$ is cartesian, \cref{Semiadditivity is cocartesianness} identifies the bi-implication \ref{Mon semiadditive condition} $\iff$ \ref{Mon vanishing algs condition} with \ref{Cocartesian rigidity condition} $\iff$ \ref{Vanishing algs condition} applied to $\Bor_I^{\cT} \uAlg^{\otimes}_{\cP}(\ucS_{\cT})$.
  \ref{Mon vanishing algs condition} $\implies$ \ref{Tensor equivalence condition} follows from \cref{Unital conservativity corollary}, and the implications \ref{Tensor equivalence condition} $\implies$ \ref{Alg vanishing algs condition} $\implies$ \ref{Mon vanishing algs condition} are obvious.
\end{proof}

\begin{proof}[Proof of \cref{Cocartesian algebras theorem}]
  After the implications illustrated in \cref{Cocartesian implications equation}, it suffices to show for all $\cD^{\otimes} \in \Cat_I^{\otimes}$ that $\uCAlg_I(\cD)$ satisfies \ref{Vanishing algs condition}, i.e. \ref{Abstractly calg condition} $\implies$ \ref{Vanishing algs condition};
  by \cref{Cocartesian rigidity workhorse lemma}, it suffices to prove that $\uCAlg_I(\ucS_{\cT})$ is $I$-semiadditive.
  But in fact, \cref{CMon is CAlg corollary} constructs an equivalence $\uCAlg_I(\ucS_{\cT}) \simeq \uCMon_I(\ucS_{\cT})$ and the latter is $I$-semiadditive by Cnossen-Lenz-Linsken's result, \cref{Semiadditivization theorem}.
\end{proof}

Rephrasing things somewhat, we've arrived at the following theorem.
\begin{repmaintheorem}{cooltheorem}{cooltheorem}{Cocartesian main theorem}
  Let $\cO^{\otimes}$ be an almost essentially reduced $\cT$-operad.
  Then, the following properties are equivalenent.
  \begin{enumerate}[label={(\alph*)}]
    \item \label[condition]{Cocartesian a} The $\cT$-$\infty$-category $\uAlg_{\cO} \ucS_{\cT}$ is $A\cO$-semiadditive.
    \item \label[condition]{Cocartesian b} The unique map $\cO^{\otimes} \rightarrow \cN^{\otimes}_{A\cO \infty}$ is an equivalence.
  \end{enumerate}
  Furthermore, for any almost essentially unital weak indexing system $I$ and $I$-symmetric monoidal $\infty$-category $\cC^{\otimes}$, the $I$-symmetric monoidal $\infty$-category $\uCAlg_{I}^{\otimes} \cC$ is cocartesian.
\end{repmaintheorem}
\begin{proof}
By \cref{Wirthmuller cocartesianness,Algebras are cartesian corollary,Cocartesian algebras theorem}, \cref{Cocartesian a} is equivalent to the forgetful $\cT$-functor
\[\begin{tikzcd}[ampersand replacement=\&, column sep=tiny]
  {\uCAlg_{A\cO}(\ucS_{\cT})} \& {\uAlg_{\cO} \uCAlg^{\otimes}_{A\cO}(\ucS_{\cT})} \& {\uCAlg_{A\cO} \uAlg^{\otimes}_{\cO}(\ucS_{\cT})} \&\&\&\& {\uAlg_{\cO}(\ucS_{\cT})}
	\arrow["\simeq"{description}, draw=none, from=1-2, to=1-1]
	\arrow["\simeq"{description}, draw=none, from=1-3, to=1-2]
	\arrow[from=1-3, to=1-7]
\end{tikzcd}\]
being an equivalence, which is equivalent to \cref{Cocartesian b} by \cref{Conservativity corollary}.
The remaining statement follows immediately from the implication \ref{Abstractly calg condition} $\implies$ \ref{Cocartesian rigidity condition} of \cref{Cocartesian algebras theorem}.
\end{proof}

This implies that $\cN_{I\infty}^{\otimes} \in \Op^{\red}_I$ is $\obv$-absorptive.
\begin{corollary}\label{Absorption corollary}
  Let $\cO^{\otimes}$ be an almost-reduced $I$-operad.
  The map $F\colon \cN^{\otimes}_{I \infty} \rightarrow \cN^{\otimes}_{I \infty} \obv \cO^{\otimes}$ is an equivalence.
\end{corollary}
\begin{proof}
  By \cref{Cocartesian algebras theorem}, the forgetful map 
  \[
    F^*\colon \Alg_{\cO \otimes \cN_{I \infty}}(\ucS_{\cT}) \simeq \Alg_{\cO} \uAlg_{\cN_{I\infty}}^{\otimes} (\ucS_{\cT}) \rightarrow \uAlg_{\cN_{I \infty}} (\ucS_{\cT})
  \]
  is an equivalence. 
  The statement then follows from \cref{Conservativity corollary}.
\end{proof}

\begin{remark}
  At this point, we may answer \cref{q: absorption to I} of the introduction;
  if $\cO^{\otimes}$ is almost-essentially unital and $\cN_{I\infty}^{\otimes} \rightarrow \cN^{\otimes}_{I \infty} \obv \cO^{\otimes}$ is an equivalence, then \cref{Support prop} implies that $\cO^{\otimes}$ is an almost-reduced $I$-operad;
  in particular, the assumptions of \cref{Absorption main corollary} are necessary and sufficient.
\end{remark}

\subsection{The smashing localization for \texorpdfstring{$\cN_{I \infty}^{\otimes}$}{N-infinity operads} and Blumberg-Hill's conjecture}\label{Smashing localization subsection}
In view of \cref{Absorption corollary}, when $I$ is an almost-unital weak indexing category the unique map $\triv_{\cT}^{\otimes} \rightarrow \cN_{I \infty}^{\otimes}$ induces an equivalence
\[
  \cN_{I \infty}^{\otimes} \simeq \cN_{I \infty}^{\otimes} \obv \triv_{\cT}^{\otimes} \xrightarrow{\;\; \sim \;\;} \cN_{I\infty}^{\otimes} \obv \cN_{I\infty}^{\otimes},
\]
i.e. it uniquely witnesses $\cN_{I\infty}^{\otimes}$ as an idempotent object in the sense of \ha{Def.}{4.8.2.1}.
To conclude \cref{Main theorem}, we will characterize the smashing localization classified by $\cN_{I\infty}^{\otimes}$-modules.\footnote{The identification between idempotent objects and smashing localizations is stated in \cite{HA,Carmeli} under the unnecessary specification that the constructions live in a given \emph{symmetric monoidal} $\infty$-category, as this leads to canonical lifts of idempotent objects to idempotent algebras.
  In fact, their arguments for the identification only make use of the underlying $\AA_2$-structure and the \emph{existence} of a braiding $A \otimes B \simeq B \otimes A$, separately for each pair.
Even before \cref{Envelopes section}, $(\Op_{\cT}, \obv, \triv_{\cT}^{\otimes})$ with the braiding determined by symmetry of the universal property for $\obv$ is certainly such a structure.}

\subsubsection{The smashing localization classified by $\cN_{I \infty}^{\otimes}$}
$(-) \obv \cN_{I\infty}^{\otimes}$ classifies algebraic Wirthm\"uller isomorphisms.
\begin{theorem}\label{The smashing localization theorem}
  Let $I$ be an almost essentially unital weak indexing system.
  Then, a $\cT$-operad $\cO^{\otimes}$ possesses an equivalence $\cP^{\otimes} \obv \cN_{I\infty}^{\otimes} \simeq \cP^{\otimes}$ if and only if the following conditions are satisfied:
  \begin{enumerate}[label={(\alph*)}]
    \item $c(\cP) = c(I)$, and
    \item $\uAlg_{\cP}(\ucS_{\cT})$ is $I$-semiadditive. 
  \end{enumerate}
\end{theorem}

By the arity support computation of \cref{Support prop}, \cref{The smashing localization theorem} is equivalent to the following.
\begin{proposition}\label{The smashing localization prop}
  Let $I$ be an almost-unital weak indexing system.
  Then, an at-least one color $\cT$-operad $\cP^{\otimes}$ satisfies $\cP^{\otimes} \obv \cN_{I\infty}^{\otimes} \simeq \cP^{\otimes}$ if and only if $\uAlg_{\cP}(\ucS_{\cT})$ is $I$-semiadditive.
\end{proposition}
\begin{proof}
  Just as in the proof of \cref{Cocartesian rigidity workhorse lemma}, note that \cref{Cocartesian algebras theorem} implies that the conditions of \cref{Cocartesian rigidity workhorse lemma} are equivalent to the additional condition
  \begin{enumerate}[label={(\alph*)}]\setcounter{enumi}{9}
    \item The forgetful $\cT$-functor $\uAlg_{\cN_{I\infty} \otimes \cP}(\ucS_{\cT}) \simeq \uCAlg_I \uAlg^{\otimes}_{\cP}(\ucS_{\cT}) \rightarrow \uAlg_{\cP}(\ucS_{\cT})$ is an equivalence.
  \end{enumerate}
  But this is equivalent to the desired equivalence by \cref{Conservativity corollary}
\end{proof}

\subsubsection{The proof of Blumberg-Hill's conjecture}
We start by answering \cref{q: BH conjecture}.
\begin{proposition}\label{Binary proposition}
  When $I$ and $J$ are almost-unital, there is an equivalence $\cN_{I \infty}^{\otimes} \obv \cN_{J \infty}^{\otimes} \simeq \cN_{I \vee J \infty}^{\otimes}$.
\end{proposition}
\begin{proof}
  By \csy{Prop.}{5.1.8}, $\cN_{I \infty}^{\otimes} \obv \cN_{J \infty}^{\otimes} \in \Op_{\cT} \simeq \Op_{\cT, \triv_{\cT}^{\otimes}/}$ is an idempotent object classifying the conjunction of the properties which are classified by $\cN_{I \infty}^{\otimes}$ and $\cN_{J\infty}^{\otimes}$;
  that is, $\cT$-operad $\cO^{\otimes}$ is fixed by $(-) \obv \cN_{I \infty}^{\otimes} \obv \cN_{J \infty}^{\otimes}$ if and only if $\uAlg_{\cO}(\ucS_{\cT})$ is $I$-semiadditive and $J$-semiadditive.
  \cref{Semiadditivity over joins}, identifies this with the property that $\uAlg_{\cO}(\ucS_{\cT})$ is $I \vee J$-semiadditive , i.e. $\cO^{\otimes}$ is fixed by $(-) \obv \cN_{I \vee J}^{\otimes}$.
  Thus, we have
  \[
    \cN_{I \vee J \infty}^{\otimes} \simeq \cN_{I \vee J \infty}^{\otimes} \obv \cN_{I \infty}^{\otimes} \obv \cN_{J \infty}^{\obv} \simeq \cN_{I \infty}^{\otimes} \obv \cN_{J \infty}.\qedhere
  \]
\end{proof}

We may now conclude the full theorem, which we restate in the atomic orbital case.
\begin{repmaintheorem}{cooltheorem}{cooltheorem}{Main theorem}
    \label{Main theorem generalized}
    $\cN_{(-)\infty}^{\otimes}\colon \wIndSys_{\cT} \rightarrow \Op_{\cT}$ restricts to a fully faithful symmetric monoidal $\cT$-right adjoint
    \[
        \begin{tikzcd}
        {\uwIndSys_{\cT}^{aE\uni, \otimes}} && {\uOp^{aE\uni, \otimes}_{\cT}.}
            \arrow[""{name=0, anchor=center, inner sep=0}, "{\cN_{(-)\infty}^{\otimes}}"', curve={height=20pt}, hook, from=1-1, to=1-3]
            \arrow[""{name=1, anchor=center, inner sep=0}, "A"', curve={height=20pt}, from=1-3, to=1-1]
            \arrow["\dashv"{anchor=center, rotate=-90}, draw=none, from=0, to=1]
        \end{tikzcd}
    \]
    Furthermore, the resulting tensor product of weak $\cN_\infty$-operads is computed by the Borelified join
    \begin{equation}\label{Ninfty equation}
      \cN_{I \infty}^{\otimes} \obv \cN_{J\infty}^{\otimes} \simeq \cN_{\Bor_{cI \cap cJ}^{\cT}(I \vee J) \infty}^{\otimes}.
    \end{equation}
    Hence whenever $I,J$ are almost-unital weak indexing categories and $\cC^{\otimes}$ is an $I \vee J$-symmetric monoidal $\infty$-category, there is a canonical equivalence of $I \vee J$-symmetric monoidal $\infty$-categories
    \begin{equation}\label{Main calg equation generalized}
        \uCAlg^{\otimes}_I \uCAlg^{\otimes}_J(\cC) \simeq \uCAlg_{I \vee J}(\cC).
    \end{equation}
\end{repmaintheorem}

\def\cNI{\cN_{I\infty}^{\otimes}}
\def\cNJ{\cN_{J\infty}^{\otimes}}
\begin{proof}[Proof of \cref{Main theorem generalized}]
  The $\cT$-adjunction is precisely \cref{Windex T-adjunction corollary}, and \cref{Ninfty equation,Main calg equation generalized} will follow from symmetric monoidality of $\cN_{(-)\infty}^{\otimes}$ and the support computation of \cref{Support prop}.

  We're left with proving that almost essentially unital weak $\cN_\infty$-operads are closed under $\obv$, i.e.
  the unique map $\varphi\colon \cN_{I \infty}^{\otimes} \obv \cN_{J \infty}^{\otimes} \rightarrow \cN_{I \vee J\infty}^{\otimes}$ is an equivalence.
 By \cref{One-color tensor product observation}, it suffices to prove that $\Bor_{cI \cap cJ}^{\cT}(\varphi)$ is an equivalence, i.e. we may assume that $I$ and $J$ are almost-unital;
 this is \cref{Binary proposition}.
\end{proof}

\subsubsection{$\cN_{I\infty}^{\otimes}$ classifies $I$-cocartesianness}
We now study a variant of \cref{The smashing localization theorem}, motivated by the following.
\begin{observation}
  The computation of \cite[\S~\href{https://people.math.harvard.edu/~lurie/papers/HA.pdf\#subsection.2.3.1}{2.3.1}]{HA} and resulting theory may be stated simply:
  the operad $\EE^{\otimes}_0$ is an idempotent object in $\Op^{\otimes}$ under the unique map $\triv^{\otimes} \rightarrow \EE_0^{\otimes}$, and the corresponding smashing localization classifies unitality.
  In particular, a symmetric monoidal $\infty$-category $\cC^{\otimes}$ is a $\EE_0^{\otimes}$-module with respect to $\obv$ if and only if the unit object $1_{\cC} \in \cC^{\otimes}$ is initial, i.e. for all $X \in \cC$, the unique map $X^{\sqcup \emptyset} \rightarrow X^{\otimes \emptyset}$ is an equivalence;
  that is, $(-) \obv \EE_0^{\otimes}$ classifies \emph{$I_0$-cocartesianness} in the symmetric monoidal case.
\end{observation}
We say that a $\cT$-operad $\cO^{\otimes}$ is \emph{$I$-cocartesian} if it is $\upsilon(I)$-unital and the identity on $U\cO$ classifies an equivalence $\Bor_I^{\cT} \cO^{\otimes} \xrightarrow{\;\;\;\; \sim \;\;\;\;} U\cO^{I-\sqcup}$.
We begin with an $I$-operadic variant of \cref{Cocartesian algebras theorem,The smashing localization prop}.
\begin{proposition}\label{Workhorse operadic case}
  Given $I$ an almost-unital weak indexing category, the following are equivalent for $\cP^{\otimes} \in \Op_{\cT}$.
  \begin{enumerate}[label={(\alph*)}]
    \item \label[condition]{Ca} $\cP^{\otimes}$ is $I$-cocartesian.
    \item \label[condition]{Cb}  For all unital $I$-operads $\cO^{\otimes}$, the forgetful functor $\Alg_{\cO}(\cP) \rightarrow \uFun_{\cT}(I\cO,\cP)$ is an equivalence.
    \item \label[condition]{Cc} The forgetful $\cT$-operad map $\uCAlg^{\otimes}_I(\cP) \rightarrow \cP^{\otimes}$ is an equivalence.
    \item \label[condition]{Cd} There exists a $\cT$-operad $\cQ^{\otimes}$ and an equivalence  of $\cT$-operads $\uCAlg_I^{\otimes}(\cQ) \xrightarrow{\;\;\sim\;\;} \cP^{\otimes}$.\setcounter{enumi}{10}
    \item The canonical map $\cP^{\otimes} \rightarrow \cP^{\otimes} \obv \cN_{I\infty}^{\otimes}$ is an equivalence.
  \end{enumerate}
\end{proposition}
\begin{proof}
  The proof of the equivalence between \cref{Ca,Cb,Cc,Cd} is identical to \cref{Cocartesian algebras theorem}, so we omit it.
  Now, by a standard two-out-of-three argument, $\cP^{\otimes}$ is local for the smashing localization associated with $\cN_{I\infty}^{\otimes}$ if and only if pullback
  along the localization map $\cO^{\otimes} \rightarrow \cO^{\otimes} \obv \cN_{I \infty}^{\otimes}$ induces an equivalence $\Alg_{\cO} \uCAlg^{\otimes}_I(\cP)^{\simeq} \simeq \Alg_{\cO \otimes \cN_{I \infty}}(\cP)^{\simeq} \xrightarrow{\;\;\; \sim \;\;\;} \Alg_{\cO}(\cP)^{\simeq}$ for all $\cO^{\otimes} \in \Op_{\cT}$, which is equivalent to the condition that $\CAlg^{\otimes}_I(\cP) \rightarrow \cP^{\otimes}$ is an equivalence by Yoneda's lemma.
\end{proof} 

Applying \cref{Support prop} to \cref{Workhorse operadic case} yields a variant of \cref{The smashing localization theorem}, answering \cref{q: absorption to O}.
\begin{proposition}\label{Generalized smashing localization theorem}
  Let $I$ be an almost essentially unital weak indexing category and $\cO^{\otimes}$ a $\cT$-operad.
Then, $\cO^{\otimes}$ admits an (essentially unique) $\cN_{I \infty}^{\otimes}$-module structure if and only if the following conditions hold:
  \begin{enumerate}[label={(\alph*)}]
    \item $c(\cO) = c(I)$, and
    \item $\cO^{\otimes}$ is $I$-cocartesian.
  \end{enumerate}
\end{proposition}

\begin{remark}
  If $J \subset I$ and $\cO^{\otimes}$ is $I$-cocartesian, then $\cO^{\otimes}$ is $J$-cocartesian.
  In particular, applying this to $I \cap \EE_{0,\upsilon(I)}^{\otimes} \subset I$ shows when $I$ is unital that the conditions of \cref{Generalized smashing localization theorem} implies that $\cO^{\otimes}$ is unital. 
\end{remark}

Given a related pair of weak indexing categories $I \subset J$, let $\Op_J^{I-\cocart} \subset \Op_J$ be the full subcategory of $I$-cocartesian $J$-operads.
We find that $\Op_J^{I-\cocart}$ is absorptive under $\obv$ and the internal hom.
\begin{corollary}\label{J-cocartesian algebras are cocartesion}
  Suppose $\cO^{\otimes},\cP^{\otimes}$ are at-least-one-color $J$-operads such that either $\cO^{\otimes}$ or $\cP^{\otimes}$ are $I$-cocartesian.
  Then, $\cO^{\otimes} \obv \cP^{\otimes}$ and $\uAlg_{\cO}^{\otimes}(\cP)$ are $I$-cocartesian.
\end{corollary}
\begin{proof}
  $I$-cocartesianness of $\cO^{\otimes} \obv \cP^{\otimes}$ follows from \cref{Generalized smashing localization theorem}.
  If $\cP^{\otimes}$ is $I$-cocartesian, then \cref{Workhorse operadic case} constructs equivalences
    \[
      \CAlg_{I} \uAlg^{\otimes}_{\cO}(\cP)
        \simeq \Alg_{\cO} \uCAlg_{I}^{\otimes}(\cP)
        \simeq \Alg_{\cO}(\cP),
    \]
  so the result follows from another application of \cref{Workhorse operadic case}.
  If $\cO^{\otimes}$ is $I$-cocartesian, the result follows from two more applications of \cref{Workhorse operadic case}. 
  as we acquire equivalences
  \begin{equation}\label{J cocart equation}
    \uCAlg_I \uAlg_{\cO}^{\otimes}(\cP) \simeq \uAlg_{\cN_{I\infty} \otimes \cO}(\cP) \simeq \uAlg_{\cO}(\cP).\qedhere
  \end{equation}
\end{proof}

Finally, we vastly extend the results of \cite[\S~\href{https://people.math.harvard.edu/~lurie/papers/HA.pdf\#subsection.2.3.1}{2.3.1}]{HA}.
\begin{corollary}\label{J-cocartesian double adjunction}
  Suppose $I \subset J$ is a related pair of almost-unital weak indexing categories.
  Then, $\Op_J^{I-\cocart} \subset \Op_J$ is a smashing localization and a cosmashing colocalization:
  \[\begin{tikzcd}[ampersand replacement=\&, column sep=large]
	{\Op_J^{I-\cocart}} \& {\Op_J}
	\arrow[""{name=0, anchor=center, inner sep=0}, hook, from=1-1, to=1-2]
\arrow[""{name=1, anchor=center, inner sep=0}, "{\cN_{I\infty}^{\otimes} \obv (-)}"', curve={height=18pt}, from=1-2, to=1-1]
	\arrow[""{name=2, anchor=center, inner sep=0}, "{\uCAlg^{\otimes}_{I}(-)}", curve={height=-18pt}, from=1-2, to=1-1]
	\arrow["\dashv"{anchor=center, rotate=-90}, draw=none, from=0, to=2]
	\arrow["\dashv"{anchor=center, rotate=-90}, draw=none, from=1, to=0]
\end{tikzcd}\]
\end{corollary} 
\begin{proof}
  \cref{Generalized smashing localization theorem} exhibits the top adjunction and \cref{Cocartesian algebras theorem} shows that $\uCAlg_I^{\otimes}\colon \Op_J \rightarrow \Op_J$ factors through $\Op_J^{I-\cocart} \subset \Op_J$.
  Moreover, applying $\Gamma^{\cT} (-)^{\simeq}$ to \cref{J cocart equation} yields a natural equivalence
  \[
    \Map_{\Op_J}\prn{\cO^{\otimes},\cP^{\otimes}} \simeq \Map_{\Op_J^{I-\cocart}}\prn{\cO^{\otimes}, \uCAlg_I^{\otimes}(\cP)} 
  \]
  for all $J$-operads $\cO^{\otimes},\cP^{\otimes}$ such that $\cO^{\otimes}$ is $I$-cocartesian, yielding the bottom adjunction.
\end{proof}

\subsubsection{(Co)localization to unital \tI-operads}\label{Unitalization subsection}
The specialization of \cref{Cocartesian algebras theorem,Generalized smashing localization theorem} to $I \cap I_0$ is quite useful, so we state it explicitly here.
\begin{corollary}\label{Maximally powerful unitality corollary}
  Given $I$ an almost-unital weak indexing system and $\cO^{\otimes} \in \Op_I$, the following are equivalent:
  \begin{enumerate}[label={(\alph*)}]
    \item \label[condition]{cond: Trivial P-algebras} For all unital $I_0$-operads $\cP^{\otimes}$, the forgetful $I$-operad map $\uAlg^{\otimes}_{\cP}(\cO) \rightarrow \uAlg_{\triv(U\cP)}^{\otimes}(\cO^{\otimes})$ is an equivalence.
    \item \label[condition]{cond: Trivial I_0-algebras} The forgetful $I$-operad map $\uAlg^{\otimes}_{\EE_{0,\upsilon(I)}}(\cO) \rightarrow \cO^{\otimes}$ is an equivalence.
    \item \label[condition]{cond: unital} $\cO^{\otimes}$ is unital.
    \item \label[condition]{cond: E0 module} There exists an equivalence $\cO^{\otimes} \simeq \EE_{0,\upsilon(I)}^{\otimes} \obv \cO^{\otimes}$.
    \item \label[condition]{cond: algebras E0 module} For all $I$-operads $\cC^{\otimes}$, the forgetful $I$-operad map $\uAlg^{\otimes}_{\cO} \uAlg^{\otimes}_{\EE_{0,\upsilon(I)}}(\cC) \rightarrow \uAlg^{\otimes}_{\cO}(\cC)$ is an equivalence.
    \item \label[condition]{cond: unital algebras} For all $I$-operads $\cC^{\otimes}$, the $I$-operad $\uAlg^{\otimes}_{\cO}(\cC)$ is unital.
    \item \label[condition]{cond: pointed algebras} The $\cT$-$\infty$-category $\uMon_{\cO}(\cS)$ is $\upsilon(I)$-pointed.
  \end{enumerate}
  In particular, 
  $\Op_I^{\uni} \subset \Op_I$ is a smashing localization and cosmashing colocalization:
  \[\begin{tikzcd}[ampersand replacement=\&]
	{\Op_I^{\uni}} \& {\Op_I}
	\arrow[""{name=0, anchor=center, inner sep=0}, hook, from=1-1, to=1-2]
	\arrow[""{name=1, anchor=center, inner sep=0}, "{\EE_{0,\upsilon(I)}^{\otimes} \obv (-)}"', curve={height=18pt}, from=1-2, to=1-1]
	\arrow[""{name=2, anchor=center, inner sep=0}, "{\uAlg^{\otimes}_{\EE_{0, \upsilon(I)}}(-)}", curve={height=-18pt}, from=1-2, to=1-1]
	\arrow["\dashv"{anchor=center, rotate=-90}, draw=none, from=0, to=2]
	\arrow["\dashv"{anchor=center, rotate=-90}, draw=none, from=1, to=0]
\end{tikzcd}\]
\end{corollary}
The double adjunction is \cref{J-cocartesian double adjunction} and the following diagram shows how to recover the corollary.
\[\begin{tikzcd}[ampersand replacement=\&, column sep=huge]
  {\text{\ref{cond: unital}}} \& {\text{\ref{cond: E0 module}}} \& {\text{\ref{cond: algebras E0 module}}} \& {\text{\ref{cond: unital algebras}}} \\
      {\text{\ref{cond: Trivial P-algebras}}} \& {\text{\ref{cond: Trivial I_0-algebras}}} \&\& {\text{\ref{cond: pointed algebras}}}
  \arrow[Rightarrow, from=1-1, to=2-1, "\ref{Workhorse operadic case}"']
	\arrow[Rightarrow, 2tail reversed, from=1-2, to=1-1, "\ref{Generalized smashing localization theorem}"']
  \arrow[Rightarrow, 2tail reversed, from=1-3, to=1-2, "\text{Yoneda}"']
  \arrow[Rightarrow, from=1-3, to=1-4, "\ref{Workhorse operadic case}"]
  \arrow[Rightarrow, from=1-4, to=2-4, "\ref{Unital SMC observation}"]
  \arrow[Rightarrow, from=2-1, to=2-2, "\text{obvious}"']
  \arrow[Rightarrow, from=2-2, to=1-1, "\ref{F-unitality lemma}" description]
  \arrow[Rightarrow, from=2-4, to=1-2, "\ref{Cocartesian rigidity workhorse lemma}"{description}]
\end{tikzcd}\]

\subsubsection{The underlying $\cT$-$\infty$-category.}%
We get an immediate corollary from \cref{Cocartesian rigidity main theorem generalized,Cocartesian algebras theorem}.
\begin{corollary}\label{Left adjoint to U corollary}
  Suppose $I$ is almost-unital.
  Then, $U_{\uni}\colon \uOp_{I}^{\uni} \rightarrow \uCat_{\cT}$ is left $\cT$-adjoint to $(-)^{I-\sqcup}$.
\end{corollary}

\begin{warning}
    \cref{Left adjoint to U corollary} shows that no nontrivial $\cT$-colimit of one-color $\cT$-operads has one color;
    in particular, no one-color $\cT$-operads are the result of a nontrivial induction.
\end{warning}

\def\cOo{\cO^{\otimes}}
\def\cPo{\cP^{\otimes}}
We use this to compute the $\cT$-$\infty$-category underlying Boardman-Vogt tensor products.
\begin{proposition}\label{Underlying category of BV tensor product}
    The underlying category functor $U|_{\uni}:\Op_{I}^{\uni} \rightarrow \Cat_{\cT}$ sends
    \[
      U\prn{\cO^{\otimes} \obv \cP^{\otimes}} \simeq U(\cO^{\otimes}) \times U(\cP^{\otimes}).
    \]
\end{proposition}
\begin{proof}
\cref{Left adjoint to U corollary,J-cocartesian algebras are cocartesion}
yield a string of natural equivalences
    \begin{align*}
      \Fun_{\cT}\prn{U\prn{\cO^{\otimes} \obv \cP^{\otimes}},\cC}
        &\simeq \Alg_{\cO \otimes \cP}\prn{\cC^{I-\sqcup}}\\ 
        &\simeq \Alg_{\cO} \uAlg_{\cP}^{\otimes}\prn{\cC^{I-\sqcup}}\\
        &\simeq \Alg_{\cO} \uFun_{\cT}\prn{U(\cP^{\otimes}),\cC}^{I-\sqcup}\\
        &\simeq \Fun_{\cT}\prn{U(\cO^{\otimes}),\uFun_{\cT}\prn{U(\cP^{\otimes}),\cC}}\\
        &\simeq \Fun_{\cT}\prn{(U(\cO^{\otimes}) \times U(\cP^{\otimes}),\cC},
    \end{align*}
    so the result follows by Yoneda's lemma.
\end{proof}
Applying \cref{One-color tensor product observation,Support prop,Underlying category of BV tensor product}, we acquire the following.
\begin{corollary}\label{aE-reduced tensor products}
  The full subcategories $\Op_{\cT}^{\red} \subset \Op_{\cT}^{a\red} \subset  \Op_{\cT}^{aE\red} \subset \Op_{\cT}$ are closed under $\obv$.
\end{corollary}

Finally, this gives us a formula for the $I$-Borelification of tensor products with $\cN_{I\infty}^{\otimes}$.
\begin{corollary}
  Given a $\cT$-operad $\cC^{\otimes}$ and an almost essentially unital weak indexing system $I$, there exists an equivalence $\Bor_I^{\cT} \prn{\cC^{\otimes} \obv \cN_{I \infty}^{\otimes}} \simeq U\prn{\cO^{\otimes} \obv \EE_{0,\upsilon(I)}^{\otimes}}^{I-\sqcup}$.
  In particular, if $\Bor_I^{\cT} \cO^{\otimes}$ is a unital $I$-operad, then
  \[
    \Bor_I^{\cT} \prn{\cC^{\otimes} \obv \cN_{I\infty}} \simeq U\cC^{I-\sqcup}.
  \]
\end{corollary}
\begin{proof}
  In sight of \cref{One-color tensor product observation} we may assume $I$ is almost-unital.
  \cref{Underlying category of BV tensor product} shows that these have the same underlying $\cT$-category and \cref{J-cocartesian double adjunction} shows that the left-hand size is cocartesian, so \cref{Cocartesian rigidity main theorem generalized} yields the corollary.
\end{proof}
We see in the next subsection that this behaves poorly outside of the almost essentially unital setting.

\subsection{Failure of the nonunital equivariant Eckmann-Hilton argument}\label{EH failure subsection}
We say that a $\cT$-operad $\cO^{\otimes}$ with at most one color is \emph{$n$-connected} if the nonempty structure spaces $\cO(S)$ are each $n$-connected.
We write the full subcategory of $n$-connected $\cT$-operads as
\[
  \Op_{\cT, \geq n+1}^{\leq \oc} \subset \Op_{\cT}^{\leq \oc}.
\]
By \cref{Conservativity corollary}, this is equivalent to the condition that the forgetful functor $\CAlg_{A\cO}(\cC) \rightarrow \Alg_{\cO}(\cC)$ is an equivalence for all $\cT$-symmetric monoidal $(n+1)$-categories , which itself is equivalent to the same condition in the case $\cC \simeq \cS_{\leq n}$.
We first observe compatibility with $\obv$ in the almost-essentially reduced setting.
\begin{corollary}\label{n-connected tensor products}
  $\Op^{aE\red}_{\cT, \geq n+1}$ is closed under $\obv$.
\end{corollary}
\begin{proof}
  By \cref{Support prop,Main theorem generalized}, $\cO^{\otimes}, \cP^{\otimes} \in \Op^{aE\red}_{\cT, \geq n+1}$ participate in a string of natural equivalences
  \begin{align*}
    \Alg_{\cO \otimes \cP}(\ucS_{\cT, \leq n})
    &\simeq \Alg_{\cO}\uAlg_{\cP}(\ucS_{\cT, \leq n})\\
    &\simeq \CAlg_{A\cO}\uCAlg_{A\cP}(\ucS_{\cT, \leq n})\\
    &\simeq \CAlg_{A\cO \vee A\cP}(\ucS_{\cT, \leq n})\\
    &\simeq \CAlg_{A(\cO \otimes \cP)}(\ucS_{\cT, \leq n}),
  \end{align*}
    induced by the unique map $\cO^{\otimes} \obv \cP^{\otimes} \rightarrow \cN_{A\prn{\cO \otimes \cP}}^{\otimes}$, so \cref{n-connected prop}, implies that $\cO^{\otimes} \obv \cP^{\otimes}$ is $n$-connected.
    \cref{aE-reduced tensor products} implies that $\cO^{\otimes} \obv \cP^{\otimes}$ is almost essentially reduced.
\end{proof}
\begin{remark}
  The unit object $\triv_{\cT}^{\otimes} \in \Op_{\cT}$ is $n$-connected for all $n$, so $n$-connected $\cT$-operads are closed under $k$-fold tensor products for all $k \in \NN$.
\end{remark}
The example $\triv_{\cT}^{\otimes} \obv \cO^{\otimes} \simeq \cO^{\otimes}$ demonstrates that this is the best we can say without further assumptions on the $\cT$-operads in question;
the author hopes to return to this question in forthcoming work, constructing analogues to \cite{Schlank}.
For the time being, we demonstrate that \cref{n-connected tensor products} dramatically fails without the \emph{almost essentially unital} assumption, exhibiting a failure of the nonunital Eckmann-Hilton argument.
\begin{observation}\label{Split codiagonal maps always exist observation}
    Fix $I$ a weak indexing system.
    By \cref{Support prop}, there is a contractible space of diagrams of the following form:
    \[
        \cN_{I \infty}^{\otimes} \simeq \cN_{I \infty}^{\otimes} \obv \triv_{c(I)}^{\otimes} \xrightarrow{\id \otimes !} \cN_{I \infty}^{\otimes} \obv \cN_{I \infty}^{\otimes} \rightarrow \cN_{I \infty}^{\otimes};
    \]
    furthermore, the composite $\cN_{I \infty}^{\otimes} \rightarrow \cN_{I \infty}^{\otimes}$ is homotopic to the identity since $\cN_{I \infty}^{\otimes}$ has contractible endomorphism space.
    In particular, this implies that there is a unique natural \emph{split diagonal} diagram
     \[
        \begin{tikzcd}[row sep=tiny]
            & {\CAlg_{I}\uCAlg^{\otimes}_{I}(-)} \\
            {\CAlg_{I}(-)} && {\CAlg_{I}(-)}
            \arrow[from=1-2, to=2-3,"U"]
            \arrow[from=2-1, to=1-2, "\delta"]
            \arrow[Rightarrow, no head, from=2-1, to=2-3]
        \end{tikzcd}
    \]
  $\delta$ takes a structure to two interchanging copies of itself, and $U$ simply forgets one of the structures.
\end{observation}
A weak $\infty$-categorical form of the \emph{Eckmann-Hilton argument} for $I$-commutative algebras would state that the functor $U$ is an equivalence, or equivalently, $\delta$ is an equivalence, i.e. $\cN_{I\infty}^{\otimes} \obv \cN_{I\infty}^{\otimes}$ is $\infty$-connected;
the specialization to $(n+1)$-categories is that $\cN_{I \infty}^{\otimes} \obv \cN_{I \infty}^{\otimes}$ is $n$-connected.
Unfortunately, this does not hold for all $I \in \wIndSys_{\cT}$.
The following simple counterexample was pointed out to the author by Piotr Pstr\k{a}gowski.
\begin{example}\label{Nonunital nonequivariant example}
  Let $R$ be a nonzero commutative ring and let $\Comm_{nu}^{\otimes}$ be the weak $\cN_\infty$-$*$-operad associated with the $*$-weak indexing system $\FF^{nu} = \FF - \cbr{\emptyset}$.
    Then, the Abelian group underlying $R$ supports a $\Comm_{nu}^{\otimes} \obv \Comm_{nu}^{\otimes}$ structure given by the two multiplications
    \[
        \mu(r,s) = rs, \hspace{50pt} \mu_0(r,s) = 0,
    \]
    which are easily seen to satisfy interchange but be distinct.
    This lies outside of the essential image of 
    \[
        \delta\colon \Alg_{\Comm_{nu}}(\Ab) \longrightarrow \Alg_{\Comm_{nu}} \uAlg_{\Comm_{nu}}(\Ab),
    \]
    so $\delta$ is not an equivalence;
    by \cref{n-connected prop}, this implies that $\Comm_{nu}^{\otimes} \obv \Comm_{nu}^{\otimes}$ is not connected.
\end{example}

In the positive direction, \cite{Schlank} yields a classification of $\obv$-idempotent algebras in \emph{reduced} $\infty$-operads.
In fact, \cref{Nonunital nonequivariant example} shows that the associated unitality assumption only misses one example among nonequivariant one-color weak $\cN_\infty$-operads.
\begin{corollary}\label{Nonunital nonequivariant corollary}
    A weak $\cN_{\infty}$-$*$-operad $\cO^{\otimes}$ possesses a map $\triv^{\otimes} \rightarrow \cO^{\otimes}$ inducing an equivalence
    \[
      \cO^{\otimes} \xrightarrow{\;\;\;\; \sim \;\;\;\;} \cO^{\otimes} \obv \cO^{\otimes}
    \]
    if and only if $\cO^{\otimes}$ is equivalent to $\triv^{\otimes}$, $\EE_0^{\otimes}$, or $\EE_\infty^{\otimes}$.
\end{corollary}
\begin{proof}
    \SY{Cor.}{5.3.4}{39} covers the reduced case, so it suffices to assume that $\cO(\emptyset) = \emptyset$ and show that $\cO^{\otimes} \simeq \triv^{\otimes}$.
    Note that $\Comm^{\otimes}_{nu}$ is the terminal nonunital $\cN_\infty$-$*$-operad, i.e. there exists a map $\cO^{\otimes} \rightarrow \Comm^{\otimes}_{nu}$, yielding a diagram
    \[
        \begin{tikzcd}
            \cO^{\otimes} \otimes \cO^{\otimes} \arrow[r,"\varphi"]
            & \Comm_{nu}^{\otimes} \otimes \Comm_{nu}^{\otimes}\\
            \cO^{\otimes} \arrow[r] \arrow[u]
            & \Comm_{nu}^{\otimes} \arrow[u]
        \end{tikzcd}
    \]
    Pulling back \cref{Nonunital nonequivariant example}, we find that if $\cO(n) = *$ for any $n\neq 1$, then $\varphi^* R \in \Alg_{\cO \otimes \cO}(\Ab)$ is not in the image of the diagonal;
    contrapositively, $\cO(n) = \emptyset$ when $n \neq 1$, i.e. it's equivalent to $\triv^{\otimes}$.
\end{proof}

We saw in \cite{Windex} that $\cbr{\triv^{\otimes} \rightarrow \EE_0^{\otimes} \rightarrow \EE_\infty^{\otimes}} = \Op^{a\uni, \text{ weak-}\cN_\infty}_*$.
In this section, we introduce an equivariant analogue to this argument in order to prove the following proposition.
\begin{proposition}\label{Nonunital EH proposition}
    Suppose $\cN_{I \infty}^{\otimes} \obv \cN_{I \infty}^{\otimes}$ is connected.
    Then, $I$ almost essentially unital.
\end{proposition}

By combining \cref{Nonunital EH proposition,Absorption corollary}, we conclude the remaining part of \cref{Absorption main corollary}.
\begin{corollary}\label{Absorption main corollary generalized}
  $\cN_{I\infty}^{\otimes} \otimes \cN_{I \infty}^{\otimes}$ is a weak $\cN_\infty$-operad if and only if $I$ is almost essentially unital;
  in particular, there exists a (necessarily unique) map $\triv^{\otimes} \rightarrow \cN_{I \infty}^{\otimes}$ inducing an equivalence $\cN_{I \infty}^{\otimes} \xrightarrow{\;\;\;\; \sim\;\;\;\;} \cN_{I \infty}^{\otimes} \obv \cN_{I \infty}^{\otimes}$ if and only if $I$ is almost-unital.
\end{corollary}

To show \cref{Nonunital EH proposition}, we pass to a universal case.
First, the weak indexing system.
\begin{recollection}
  In \cite{Windex}, we computed the terminal weak indexing system with unit family $\cF$ to be
  \[
    \FF_{\cF^{\perp}-nu,V} = \begin{cases}
      \FF_V & V \in \cF; \\ 
      \FF_V - \cbr{S \mid \forall U \in \Orb(S), \, U \in \cF} & V \not\in \cF; 
    \end{cases}
  \]
  in particular, $\uFF_I$ fails to be almost essentially unital if and only if there is some non-contractible $W$-set in $\FF_{\upsilon(I)^{\perp}-nu,W} \cap \FF_{I,W}$ for some $W \in \upsilon(I)^{\perp}$.
  We refer to the associated weak indexing category as $I_{\cF^{\perp}-nu}$;
  note that $I_{\cF^{\perp}-nu} \subset \FF_{\cT}$ is the wide subcategory of maps $T \rightarrow S$ such that either $S,T \in \FF_{\cF}$ or $S,T \in \FF_{\cF^{\perp}}$.
\end{recollection}
Now, we construct a family of problematic $\cN_{I_{\cF^{\perp}}-nu}^{\otimes} \obv \cN_{I_{\cF^{\perp}}-nu}^{\otimes}$-algebras.
\begin{construction}
  Let $M$ be a $\cT$-commutative monoid in pointed sets.
  We define a new functor
  \[
    M^0\colon h_1\Span_{I_{\cF^{\perp}-nu}}(\FF_{\cT}) \rightarrow \Set_*
  \]
  which agrees with $M$ on objects, backwards maps, forwards maps lying in $\FF_{\cF}$, but whose forward maps lying in $\FF_{\cF^{\perp}}$ are zero.
  This is evidently functorial on backwards and forward maps, and the restriction to backwards maps is product-preserving.
  We're left with verifying the double coset formula that, given a cartesian square as on the left such that $S,T,R \in \cT$ and $f,f' \in I_{\cF^{\perp}-nu}$, the right square commutes, where $(-)_*$ denotes covariant functoriality and $(-)^*$ contravariant.
  \[
    \begin{tikzcd}[ampersand replacement=\&, row sep=small]
	\& {R \times_S T} \&\&\&\& {M^0(R \times_ST)} \\
	R \&\& T \&\& {M^0(R)} \&\& {\widetilde M^0(T)} \\
	\& S \&\&\&\& {M^0(R)}
	\arrow["{g'}"{description}, from=1-2, to=2-1]
	\arrow["{f'}"{description}, from=1-2, to=2-3]
	\arrow["\lrcorner"{anchor=center, pos=0.125, rotate=-45}, draw=none, from=1-2, to=3-2]
	\arrow["{f'_*}"{description}, from=1-6, to=2-7]
	\arrow["f"{description}, from=2-1, to=3-2]
	\arrow["g"{description}, from=2-3, to=3-2]
	\arrow["{g^{\prime *}}"{description}, from=2-5, to=1-6]
	\arrow["{f_*}"{description}, from=2-5, to=3-6]
	\arrow["{g^*}"{description}, from=3-6, to=2-7]
\end{tikzcd}
  \]
  The assertions that $f,f' \in I_{\cF^{\perp}-nu}$ and that $\cF$ is a family together imply that $T \in \cF$ if and only if the entire diagram lives in $\FF_{\cF}$, and $T \in \cF^{\perp}$ if and only if the entire diagram lives in $\FF_{\cF^{\perp}}$. 
  In the former case, the right diagram commutes by the double coset formula for $M$, and in the latter case it commutes as each composite map is zero.
\end{construction}
\begin{lemma}\label{Interchanging algebras}
For all $\cT$-commutative algebras $M$, the $I_{\cF^{\perp}-nu}$-commutative algebras $M$ and $M_0$ interchange.
\end{lemma}
\begin{proof}
  Note a diagram of $\cT$-coefficient systems in a 1-category commutes if and only if the $V$-fixed point diagram commutes for all $V \in \cT$;
  the $V$-fixed points of the diagram in \cref{1-categorical recollection} in our case correspond with the diagram
  \[
    \begin{tikzcd}[ampersand replacement=\&, column sep = large]
      {\prn{X^T}^S} \& {X^{S \times T} } \& {\prn{X^S}^T} \& {X^T} \\
      {X^S} \&\&\& X
      \arrow["\simeq"{description}, draw=none, from=1-1, to=1-2]
      \arrow["{(\tr_T)_S}", from=1-1, to=2-1]
      \arrow[draw=none, from=1-2, to=1-3]
      \arrow["\simeq"{description}, draw=none, from=1-2, to=1-3]
      \arrow["{\prn{\tr^0_S}_T}", from=1-3, to=1-4]
      \arrow["{\tr_T}", from=1-4, to=2-4]
      \arrow["{\tr_S^0}", from=2-1, to=2-4]
    \end{tikzcd}
  \]
  where $\tr_*$ is the indexed multiplication in $M$ and $\tr_*^0$ is the indexed multiplication in $M^0$;
  when $V \in \cF^{\perp}$, this commutes as each of the composites factor through a zero map.

  Moreover, note that $\Bor_{\cF}^{\cT} I_{\cF^{\perp}-nu}$ is unital, so \cref{Absorption corollary} implies that $\Bor_{\cF}^{\cT} I_{\cF^{\perp}-nu}$-algebras interchange with themselves;
  in particular, the interchange relation of \cref{1-categorical recollection} for $M$ with itself implies the same relation for $M$ and $M_0$ whenever $V \in \cF$.
\end{proof}
We are left with constructing a highly nontrivial $\cT$-commutative algebra;
we choose a universal one.
\begin{construction}
  Since the ``isomorphism classes of objects'' functor $\pi_0\colon \Cat \rightarrow \Set$ preserves limits, pushforward along it lifts to a functor
  \[
    \upi_0\colon \Cat_{\cT}^{\otimes} \simeq \CMon_{\cT}(\Cat) \rightarrow \CMon_{\cT}(\Set);
  \]
  the \emph{effective Burnside $\cT$-commutative monoid} is $\uA_{\cT} \deq \upi_0 \uFF_{\cT}^{\cT-\sqcup}$.
  We denote its image under the maps 
  \[
    \CMon_{\cT}(\Set) \simeq \CMon_{\cT}(\Set_*) \rightarrow \CMon_{I_{\cF^{\perp}-nu}}(\Set_*)
  \]
  implied by \cref{CMon is CAlg corollary,Absorption corollary} by $\widetilde \uA_{\cT}$.
\end{construction}
\begin{lemma}\label{Burnside lemma}
  The $S$-indexed multiplication in $\widetilde \uA_{\cT}$ and $\widetilde \uA^0_{\cT}$ are distinct for all $S \in \FF_{\cF^{\perp}-nu,V} - \cbr{*_V}$ and $V \in \cF^{\perp}$.
\end{lemma}
\begin{proof}
  It suffices to prove that, for all $S \neq \emptyset_{V} \in \FF_{V}$, the $S$-ary multiplication of $\uA_{\cT}$ takes some element to another element other than the unit;
  unraveling definitions, this is equivalent to the property that some nonempty $V$-set can be expressed as an $S$-indexed coproduct.
  $S$ provides such an example.
\end{proof}
We now exhibit failure of the non-almost-essentially-unital 1-categorical Eckmann-Hilton argument.
\begin{proof}[Proof of \cref{Nonunital EH proposition}]
Note that
\begin{align*}
  \cN_{I\infty}^{\otimes} \obv \cN_{I \infty}^{\otimes} \; \text{ is connected} 
  & \iff h_{1} \cN_{I \infty}^{\otimes} \obv \cN_{I \infty}^{\otimes} \simeq \cN^{\otimes}_{A \prn{\cN_{I \infty} \otimes \cN_{I \infty}}} \simeq \cN_{I \infty}^{\otimes}\\ 
  &\implies \CAlg_I(\Set_*) \longrightarrow \Alg_{\cN_{I\infty} \otimes \cN_{I \infty}}(\Set_*) \; \text{ is essentially surjective}.
\end{align*}
  Furthermore, \cref{Interchanging algebras,Burnside lemma} construct an $\cN_{\upsilon(I)^{\perp}-nu \infty}^{\otimes} \obv \cN_{\upsilon(I)^{\perp}-nu \infty}^{\otimes}$-algebra $A$ satisfying the condition that its two individual structure maps $A(S) \rightarrow A(*_V)$ differ whenever $V \in \upsilon(I)^{\perp}$ and $S \neq *_V$.
  Since $I$ is not almost essentially unital, it must admit some noncontractible $S \in \FF_{I,V}$ for $V \in \upsilon(I)^{\perp}$, so the pullback $\cN_{I \infty}^{\otimes} \obv \cN_{I \infty}^{\otimes}$ structure on $A$ has two distinct underlying $I$-algebra structures, implying it is outside of this essential image.
  The contrapositive shows that $\cN_{I \infty}^{\otimes} \obv \cN_{I \infty}^{\otimes}$ is not connected.
\end{proof}

\begin{remark}\label{Idempotents remark}
  Using the above argument, one can show that if $\cO^{\otimes}$ is a idempotent object in $\cT$-operads, then its nullary spaces $\cO(\emptyset_V)$ are nonempty.
  If additionally $\cO(\emptyset_V)$ are assumed to be contractible (i.e. $\cO^{\otimes}$ is almost-unital), then \cref{Underlying category of BV tensor product} shows that the underlying fixed point catgeories $\cO_V$ are all idempotent algebras, i.e. they are contractible.
  Hence $\cO^{\otimes}$ will be shown to be almost-reduced.
  In forthcoming work, we will develop an equivariant lift of \cite{Schlank}, which would imply that every idempotent almost-unital $\cT$-operad is a weak $\cN_\infty$-operad.
\end{remark}

\section{Corollaries in higher algebra}
We now indulge in a number of corollaries.
We begin in \cref{Envelopes section} by making the equivariant equifibered perspective symmetric monoidal, establishing that $\obv$ is $\cT$-bifunctorial and distributive, fulfilling \cref{Bifunctor assumption} in full generality.
We make use of the distributivity in \cref{sec:tensor disintegration}, where we apply the disintegration and assembly procedure of \cref{sec:disinte gration} to compute tensor products of $\cT$-space colored $\cT$-operads.
Then, in \cref{Right modules subsection} we show how \cref{Absorption main corollary} constructs an $I$-symmetric monoidal structure on right-modules over an $I$-commutative algebra.

Moving on from this, in \cref{Infinitary Dunn additivity} we spell out the cases of \emph{equivariant Dunn additivity} that follows from \cref{Absorption main corollary,Main theorem};
in \cref{THR subsection} we show that equivariant factorization homology is $G$-symmetric monoidal and spell out how to infinitely iterate equivariant factorization homology and $\TCR$.

\subsection{Coherences and restrictions of equivariant Boardman-Vogt tensor products}\label{Envelopes section}
We would like to construct coherences for $\obv$ using the argument of \cite{Barkan_Segal}, but it is currently not known whether $\Env_I$ is monic in $\Cat$, so we must modify their argument:
we use that the \emph{sliced} envelope fully-faithful \Barkan{Prop.}{4.2.1}.

Now, we proved in \cref{Absorption main corollary generalized} that $\Comm_{\cT}^{\otimes}$ bears a unique structure as an idempotent object;
in particular, \cref{Env prop} shows that $\Env$ is compatible with finitary tensor products, so it induces a $\circledast$-idempotent object structure on 
$\uFF_{\cT}^{\cT-\sqcup} \in \Cat_{\cT}^{\otimes}$.
By \ha{Prop.}{4.8.2.9}, this underlies a unique $\EE_\infty$-algebra under the mode structure.
Then, \cref{I-symmetric monoidal slice corollary} constructs a symmetric monoidal $\cT$-$\infty$-category structure on the $\cT$-overcategory $\uCat_{\cT, /\uFF_{\cT}^{\cT-\sqcup}}^{\otimes}$ whose underlying tensor functor has value
\[
  \cC \circledast \cD \xrightarrow{\pi_{\cC} \circledast \pi_{\cD}} \uFF_{\cT}^{\cT-\sqcup} \circledast \uFF_{\cT}^{\cT-\sqcup} \xrightarrow{\;\;\;\; \sim \;\;\;\;} \uFF_{\cT}^{\cT-\sqcup}.
\] 
and whose unit is
\[
  \Env\prn{\triv_{\cT}^{\otimes}} \xrightarrow{\eta} \uFF_{\cT}^{\cT-\sqcup}.
\]
We acquire existence and uniqueness of the Boardman-Vogt symmetric monoidal structure.
\begin{repmaintheorem}{coolcorollary}{cooltheorem}{Operad main theorem}\label{Env coherence corollary}
  $\uOp_{\cT} \subset \uCat_{\cT,/\uFF_{\cT}^{\cT-\sqcup}}^{\otimes}$ is a symmetric monoidal subcategory under $\circledast$, with unit corresponding with $\triv_{\cT}^{\otimes}$ and tensor bifunctor corresponding with $\obv$.
  Hence there exists a unique symmetric monoidal $\cT$-$\infty$-category $\uOp_{\cT}^{\otimes}$ and symmetric monoidal $\cT$-functor
    \[
      \uOp_{\cT}^{\otimes} \rightarrow \uCat_{\cT,/\uFF_{\cT}^{\cT-\sqcup}}^{\otimes, \circledast} 
    \]
    lifting the sliced $\cT$-symmetric monoidal envelope.
\end{repmaintheorem}
\begin{proof}
  We're tasked with proving that the image of $\Env^{/\uFF_{\cT}^{\cT-\sqcup}}(-)$ contains the unit and is closed under binary tensor products.
  The unit is \cref{Env prop}, and by construction we have a commutative diagram
  \[
    \begin{tikzcd}[column sep=huge]
      {\Env\prn{\cO^{\otimes} \obv \cP^{\otimes}}} & {\Env\prn{\Comm_{\cT}^{\otimes} \obv \Comm_{\cT}^{\otimes}} } & {\Env\prn{\Comm_{\cT}^{\otimes}} } \\
      {\Env\prn{\cO^{\otimes}} \circledast \Env\prn{ \cP^{\otimes}}} & {\uFF_{\cT}^{\cT-\sqcup} \circledast \uFF_{\cT}^{\cT-\sqcup}} & {\uFF_{\cT}^{\cT-\sqcup}}
      \arrow["{\Env(\pi_{\cO} \obv \pi_{\cP})}", from=1-1, to=1-2]
      \arrow["\simeq"{marking, allow upside down}, draw=none, from=1-1, to=2-1]
      \arrow["\simeq"{marking, allow upside down}, draw=none, from=1-2, to=2-2]
      \arrow["{\Env(\id \otimes \eta)}"', "\sim", from=1-3, to=1-2]
      \arrow["\simeq"{marking, allow upside down}, draw=none, from=1-3, to=2-3]
      \arrow["{\pi_{\Env(\cO^{\otimes})} \circledast \pi_{\Env(\cP^{\otimes})}}"', from=2-1, to=2-2]
      \arrow["{\id \circledast \eta}", "\sim" above, from=2-3, to=2-2]
    \end{tikzcd}  
  \]
  Inverting $\Env(\id \otimes \eta)$ and $\id \circledast \eta$ yields the desired equivalence
  \[
    \Env^{/ \uFF_{\cT}^{\cT-\sqcup}}\prn{\cO^{\otimes} \obv \cP^{\otimes}} 
    \simeq \Env^{/ \uFF_{\cT}^{\cT-\sqcup}}\prn{\cO^{\otimes}} 
    \circledast \Env^{/ \uFF_{\cT}^{\cT-\sqcup}}\prn{\cP^{\otimes}}. \qedhere
  \] 
\end{proof}

\begin{corollary}\label{Assumption (b) is attained}
  There exists a $G$-bifunctor $\obv\colon \uOp_{\cT} \times \uOp_{\cT} \longrightarrow \uOp_{\cT}$ 
  whose $V$-value $\Op_V \times \Op_V \rightarrow \Op_V$ is
  \[
    \cO^{\otimes} \obv \cP^{\otimes} \simeq L_{\Op_V} \prn{\cO^{\otimes} \times \cP^{\otimes} \longrightarrow \Span(\FF_V) \times \Span(\FF_V) \xrightarrow{\;\;\; \wedge \;\;\;} \Span(\FF_V)}.
  \]
\end{corollary}

This recovers the correct construction in the nonequivariant case.

\begin{corollary}
  When $\cT = *$, there is an equivalence of symmetric monoidal $\infty$-categories
  \[
    \Op_*^{\otimes} \simeq \Op^{\otimes},
  \]
  where the latter is the Boardman-Vogt symmetric monoidal $\infty$-category of \cite{Barkan_Segal}.
\end{corollary}
\begin{proof}
  In \cite{EBV} we supplied an equivalence $\Op_* \simeq \Op$, so it suffices to upgrade this to a symmetric monoidal equivalence.
  In fact, the unslicing functor $\Cat^{\otimes}_{\infty, /\FF^{\sqcup}} \rightarrow \Cat^{\otimes}_{\infty}$ bears a symmetric monoidal structure (see \cref{Unslicing is symmetric monoidal}), so \cref{Env coherence corollary} constructs a symmetric monoidal structure on the composite $\Op_*^{\otimes} \rightarrow \Cat_\infty^{\otimes}$.
  Thus \cite[Thm~\href{https://arxiv.org/pdf/2301.08650v3\#thmA.5}{E}]{Barkan_Segal} constructs a symmetric monoidal equivalence extending the equivalence $\Op_* \simeq \Op$.
\end{proof}

\begin{corollary}\label{Eveerything is closed under tensor products}
  Let $I$ be a one color weak indexing system and $n \in \NN \cup \cbr{\infty}$.
  Then, $\Op_I \subset \Op_{\cT}$ is a symmetric monoidal subcategory, $\Op_I^{\uni} \subset \Op_{I}$ is a smashing localization, and the following are symmetric monoidal full subcategory inclusions:
  \[
    \begin{tikzcd}[row sep=0em, column sep=tiny]
      {\Op^{aE\red}_{I,  \geq n}} & {\Op^{aE\red}_{I}} & {\Op^{aE\uni}_{I}} & {\Op_I}\\
      {\Op^{\red}_{I,  \geq n}} & {\Op^{\red}_{I}} & {\Op^{\uni}_{I}}
      \arrow["\subset"{description}, draw=none, from=1-1, to=1-2]
      \arrow["\subset"{marking, allow upside down}, draw=none, from=1-2, to=1-3]
      \arrow["\subset"{marking, allow upside down}, draw=none, from=1-3, to=1-4]
      \arrow["\subset"{marking, allow upside down}, draw=none, from=2-1, to=2-2]
      \arrow["\subset"{marking, allow upside down}, draw=none, from=2-2, to=2-3]
    \end{tikzcd}
  \]
\end{corollary}
\begin{proof}
  The first statement follows from \cref{Support prop} and the second from \cref{Maximally powerful unitality corollary}.
  $\triv_\cT^{\otimes}$ and $\EE_{0,\upsilon(I)}^{\otimes}$ are $\infty$-connected, so in particular, the symmetric monoidal units are compatible with each of the above subcategory inclusions.
  We're left with verifying that each subcategory inclusion is closed under tensor products;
  the lefthand inclusions both follow from \cref{n-connected tensor products},
  the middle inclusions both follow from \cref{Underlying category of BV tensor product}, and
  the righthand inclusion is \cref{aE-reduced tensor products}.
\end{proof}

We finish the subsection by confirming a convenient structural result.
\begin{corollary}\label{Distributtivity of BV tensor products}
  $\uOp^{\otimes}_I$, and $\uOp^{J-\cocart,\otimes}_I$ are presentably symmetric monoidal $\cT$-$\infty$-categories.
\end{corollary}

To see this, $\uOp^{\otimes}_I$ is presentable by the localizing inclusion $\uOp_I \subset \uCat_{\cT, /\uFF_{I}^{I-\sqcup}}$, and it is distributive by the tensor-hom $\cT$-adjunction $(-) \obv \cO^{\otimes} \dashv \uAlg_{\cO}^{\otimes}(-)$.
The remaining case follows from the following easy lemma.
\begin{lemma}
  If $L\colon \cD^{\otimes} \rightarrow \cC^{\otimes}$ is a smashing $\cT$-localization and $\cD^{\otimes}$ is a presentably symmetric monoidal $\cT$-$\infty$-category, then $\cC^{\otimes}$ is a presentably symmetric monoidal $\cT$-$\infty$-category.
\end{lemma}
\begin{proof}
  It's clear that $\cC$ is a presentable $\cT$-$\infty$-category, so we're left with verifying that $- \otimes C\colon \cC \rightarrow \cC$ possesses a right $\cT$-adjoint;
  by the usual argument, it suffices to show this on fixed points, so we may assume $\cT = *$. 
  
  We claim that $- \otimes C \dashv \hom_{\cD}(C,-)$, the latter denoting the $\cD$-internal hom.
  It suffices to verify that $\hom_{\cD}(C,D)$ is $L$-local for all $D \in \cD$.
  We apply the standard argument: if $f\colon X \rightarrow Y$ is an $L$-equivalence, then $C \otimes f \sim C \otimes Lf\colon C \otimes X \rightarrow C \otimes Y$ is an equivalence, so the horizontal arrows in the following are equivalences.
    \[\begin{tikzcd}[ampersand replacement=\&, row sep=small]
	{\Map(Y,\hom(C,D))} \& {\Map(X,\hom(C,D))} \\
	{\Map(Y \otimes C, D)} \& {\Map(X \otimes C, D)}
	\arrow["{f^*}", from=1-1, to=1-2]
	\arrow["\simeq"{marking, allow upside down}, draw=none, from=1-1, to=2-1]
	\arrow["\simeq"{marking, allow upside down}, draw=none, from=1-2, to=2-2]
	\arrow["{(C \otimes f)^*}", from=2-1, to=2-2]
\end{tikzcd}\]
  The fact that the top arrow is an equivalence is the desired locality.
\end{proof}
 
\subsection{Disintegration and equivariant Boardman-Vogt tensor products}\label{sec:tensor disintegration}
We show the following generalization of the main results of \cite[\S~2.3.3-2.3.4]{HA} in \cref{sec:disinte gration}.
\begin{theorem}[Disintegration and assembly]\label{thm:general disintegration and assembly}
  Let $X$ be a $\cT$-space.
  Taking fibers yields an equivalence
  \[
    \uOp_{I, /X^{I-\sqcup}} \simeq \uFun_{\cT}(X, \uOp_I).
  \]
  The counit of this specifies a natural equivalence
  \[
    \ucolim_{x \in X} \prn{\Res_{\stab(x)}^{\cT} \cO^{\otimes} \times_{\Res_{\stab(x)}^{\cT} X^{I-\sqcup}} \cN_{I_V\infty}^{\otimes}}
    \xrightarrow{\;\;\; \sim \;\;\;} \cO^{\otimes}.
  \]
\end{theorem}
Here, $\ucolim_{x \in X}$ refers to a $\cT$-colimit of an $X$-indexed diagram.
Given $x \in X^V$ we've written $\stab(x) \deq V$.
Given $\cO^{\otimes}$ a $\cT$-operad, $I$ a one-color $\cT$-weak indexing category, and $x \in \cO_V$ a $V$ object, we define the \emph{reduced endomorphism $I_V$-operad of $x$} to be the pullback
\[
  \begin{tikzcd}[ampersand replacement=\&]
    {\End_x^{I,\red}(\cO)} \& {\Res_V^{\cT} \cO^{\otimes}} \\
    {\cN_{I_V \infty}^{\otimes}} \& {\Res_V^{\cT} U\cO^{I-\sqcup}}
    \arrow["{\iota_x}", from=1-1, to=1-2]
    \arrow["{!}"', from=1-1, to=2-1]
    \arrow["\lrcorner"{anchor=center, pos=0.125}, draw=none, from=1-1, to=2-2]
    \arrow["\eta", from=1-2, to=2-2]
    \arrow["{\cbr{x}}"', from=2-1, to=2-2]
  \end{tikzcd}
\]
In the case $I = \cT$, we simply write $\End_x^{\red}(\cO) \deq \End_x^{\cT,\red}(\cO)$.
\begin{remark}
  If $\cO^{\otimes}$ is unital (resp. almost-unital) then $\End_x^{\red}(\cO)$ is reduced (almost-reduced).
\end{remark}
We acquire the following from \cref{thm:general disintegration and assembly}.
\begin{corollary}\label{thm:disintegration and assembly}
  Suppose $\cO^{\otimes}$ is a $\cT$-operad whose underlying $\cT$-$\infty$-category $U\cO$ is a $\cT$-space and $I$ is a one-color weak indexing system.
  Then, the inclusion maps $\iota_x$ assemble to a $\cT$-colimit diagram in $I$-operads:
  \[
    \ucolim_{x \in U\cO} \End_x^{I,\red}(\cO) \xrightarrow{\;\;\;\; \sim \;\;\;\;} \Bor_I^{\cT} \cO.
  \]
\end{corollary}

In essence, this says that an at-least one color $\cT$-operad $\cO^{\otimes}$ whose underlying $\cT$-$\infty$-category is a $\cT$-space disintegrates into a $U\cO$-local system of one color $\cT$-operads, and the $U\cO$-indexed colimit $\cT$-operad (i.e. the Grothendieck construction) assembles $\cO^{\otimes}$ from this local system.
In particular, $\cO$-algebras are $U\cO$-indexed systems of $\End_x^{\red}(\cO)$-algebras:
\begin{align*}
  \uAlg_{\cO}(\cC) 
  &\simeq \uAlg_{\ucolim_{x \in U\cO} \End_x^{\red}(\cO)}(\cC)\\
  &\simeq \ulim_{x \in U\cO} \uAlg_{\End_x^{\red}(\cO)}(\cC).
\end{align*}
The corresponding picture for $\cO^{\otimes} \obv \cP^{\otimes}$-algebras is $U\cO \times U\cP$-local systems of $\End_x^{\red}(\cO) \obv \End_y^{\red}(\cP)$-algebras: 
that is, we can compute tensor products of at-least one color $\cT$-operads in terms of one color $\cT$-operads, as long as they are $\cT$-space colored.
\begin{corollary}[Disintegration of tensor products] \label{Disintegration of tensor products}
  Suppose $\cO^{\otimes},\cP^{\otimes},\cQ^{\otimes}$ are at-least one colored $\cT$-operads whose underlying $\cT$-$\infty$-categories are $\cT$-spaces and $\varphi\colon \cO^{\otimes} \obv \cP^{\otimes} \rightarrow \cQ^{\otimes}$
  is a map such that
  \begin{enumerate}[label={(\alph*)}]
    \item the underlying map of $\cT$-spaces $U\varphi\colon U\cO \times U\cP \rightarrow U\cQ$ is an equivalence, and
    \item for all pairs $(x,y) \in U\cO \times U\cP$, $\varphi$ pulls back to an equivalence
      \[
        \varphi_{(x,y)}\colon \End_x^{\red}(\cO) \obv \End_y^{\red}(\cP) \rightarrow \End_{(x,y)}^{\red}(\cQ).
      \]
  \end{enumerate}
  Then $\varphi$ is an equivalence.
\end{corollary}
\begin{proof}
    \cref{Distributtivity of BV tensor products,thm:disintegration and assembly} construct equivalences of arrows
    \[\begin{tikzcd}[ampersand replacement=\&, column sep=small]
	{\cO^{\otimes} \obv \cP^{\otimes}} \& {\ucolim_{x \in U\cO} \End_x^{\red}(\cO) \obv \ucolim_{y \in U\cP} \End_y^{\red}(\cP)} \& {\ucolim_{(x,y) \in U\cO \times U\cP} \End_x^{\red}(\cO) \obv \End_y^{\red}(\cP)} \\
  {\cQ^{\otimes}} \& {\ucolim_{(x,y) \in U\cO \times U\cP} \End_{(x,y)}^{\red}(\cQ)} \& {\ucolim_{(x,y) \in U\cO \times U\cP} \End_{(x,y)}^{\red}(\cQ)}
	\arrow["\simeq"{description}, draw=none, from=1-1, to=1-2]
	\arrow[from=1-1, to=2-1, "\varphi"]
	\arrow["\simeq"{description}, draw=none, from=1-2, to=1-3]
	\arrow[from=1-2, to=2-2]
	\arrow["\sim", from=1-3, to=2-3]
	\arrow["\simeq"{description}, draw=none, from=2-1, to=2-2]
	\arrow[equals, from=2-2, to=2-3]
\end{tikzcd}\]
  The right vertical arrow is an equivalence by assumption, so $\varphi$ is an equivalence by two out of three.
\end{proof}
We will make crucial use of this in forthcoming work concerning variants of $\EE^{\otimes}_V$ with tangential structure.

\subsection{Norms of right-modules over \tI-commutative algebras}\label{Right modules subsection}
Let $I$ be an indexing category, $\cC^{\otimes}$ an $I$-symmetric monoidal $\infty$-category, $t\colon W \rightarrow V$ an $I$-admissible transfer, $A$ an $I_V$-commutative algebra, and $M$ a right module over the associative algebra underlying $\Res_W^V A$.
Then, we may define the \emph{$A$-module norm of $M$} by the base-changed $A$-module
\[
    \,_A N_W^V M \deq A \otimes_{N_W^V \Res_W^V A} N_W^V M;
\]
that is, the normed multiplication recognizes $A$ as an $N_W^V \Res_W^V A$-module, and the $A$-module norm of $M$ is the free $A$-module on the normed $N_W^V \Res_W^V A$-algebra of $M$;
see \cite{Yang} for a detailed account in the $C_p$-equivariant case.

In this subsection, we use the equivalence $\cN_{I \infty}^{\otimes} \simeq \EE_1^{\otimes} \obv \cN_{I \infty}^{\otimes}$ to lift this to an $I$-symmetric monoidal structure, yielding coherent functoriality and a coherent double coset formula for $A$-module norms.
To do this, we begin by bootstrapping $G$-symmetric monoidality of the right module construction from the non-equivariant case.
\begin{observation}
  Fix $\cO^{\otimes}$ a $\cT$-operad.
  By \ha{Rmk.}{4.8.3.8}, functors $F\colon \Tot \Tot_{\cT} \cO^{\otimes} \rightarrow \Cat^{\Alg}$ are data
  \[
    \begin{tikzcd}[ampersand replacement=\&]
    \& {\cC^{\otimes}} \\
    {\EE_1^{\otimes} \times \Tot \Tot_{\cT} \cO^{\otimes}} \& {\EE_1^{\otimes} \times \Tot\Tot_{\cT} \cO^{\otimes}}
    \arrow[from=1-2, to=2-2, "\pi_F"]
    \arrow[from=2-1, to=1-2, "A_F"]
    \arrow[equals, from=2-1, to=2-2]
  \end{tikzcd}
  \]
  such that $\pi$ is a cocartesian fibration whose fibers $\cC_V^{\otimes} \rightarrow \EE_1^{\otimes}$ are the unstraightenings of small monoidal $\infty$-categories and such that the composite arrows $\EE_1^{\otimes} \times \cbr{O} \hookrightarrow \cC^{\otimes}$ are associative algebras.
  Moreover, unwinding definitions and applying \cref{O-monoid lemma}, the condition that $F$ corresponds with an $\cO$-monoid $\Tot_{\cT} \cO^{\otimes} \rightarrow \uCoFr^{\cT} \Cat^{\Alg}$ corresponds with the condition that each of the fibers $\cC_n^{\otimes} \rightarrow \Tot \Tot_{\cT} \cO^{\otimes}$ is an $\cO$-monoidal $\infty$-category.
\end{observation}
Given $\cC^{\otimes} \in \Cat_{\EE_1 \otimes \cO}^{\otimes}$ and $A \in \Alg_{\EE_1 \otimes \cO}(\cC)$, we acquire a functorial diagram
  \[
    \begin{tikzcd}[ampersand replacement=\&, sep = small]
      {\Tot_{\cT} \cO^{\otimes}} \\
      \& {\uCoFr^{\cT} \Cat^{\Alg} \times_{\uCoFr^{\cT}\Alg (\Cat)} \Tot_{\cT} \cO^{\otimes}} \& {\uCoFr^{\cT} \Cat^{\Alg}} \& {\uCoFr^{\cT} \Cat^{\Mod}} \& {\uCoFr^{\cT} \Cat} \\
      \&\& {\uCoFr^{\cT} \Alg(\Cat)} \\
      \& {\Tot_{\cT} \cO^{\otimes}} \& {\uAlg(\uCoFr^{\cT} \Cat)}
      \arrow[dashed, from=1-1, to=2-2]
      \arrow["{\prn{A,\cC^{\otimes}}}"{description}, curve={height=-8pt}, from=1-1, to=2-3]
      \arrow["{\RMod_A(\cC)^{\otimes}}" description, curve={height=-24pt}, from=1-1, to=2-5]
      \arrow[curve={height=24pt}, equals, from=1-1, to=4-2]
      \arrow[from=2-2, to=2-3]
      \arrow[from=2-2, to=4-2]
      \arrow["\lrcorner"{anchor=center, pos=0.125}, draw=none, from=2-2, to=4-3]
      \arrow["\Theta"', from=2-3, to=2-4]
      \arrow[from=2-3, to=3-3, "U"]
      \arrow["Y"', from=2-4, to=2-5]
      \arrow["\simeq"{marking, allow upside down}, draw=none, from=3-3, to=4-3]
      \arrow[from=4-2, to=4-3, "\cC^{\otimes}"]
    \end{tikzcd}
  \]
  $\theta$ and $Y$ are product preserving functors, so $\RMod_A(\cC)^{\otimes}$ is an $\cO$-monoid in $\Cat$, i.e. an $\cO$-monoidal $\infty$-category. 
Unwinding definitions, this proves the following proposition.
  \begin{proposition}
  Let $\cO^{\otimes}$ be a $\cT$-operad and let $\cC^{\otimes}$ an $\EE_1 \otimes \cO$-monoidal $\infty$-category.
  There is a lift
  \[
    \begin{tikzcd}[column sep=huge, ampersand replacement=\&]
      \&\& {\Cat_{\cO}^{\otimes}} \\
      {\Alg_{\cO \otimes \EE_1}(\cC)} \& {\Alg_{\EE_1}(\cC)} \&  \Cat,
      \arrow[from=2-1, to=2-2]
      \arrow["{\Gamma^{\cT}}", from=1-3, to=2-3]
      \arrow["{\RMod_{(-)}^{\otimes}}(\cC)", curve={height=-12pt}, from=2-1, to=1-3]
      \arrow[from=2-2, to=2-3, "\RMod_{(-)}(\cC)"]
    \end{tikzcd}   
  \]
  natural separately in $\cO^{\otimes}$ and $\cC^{\otimes}$;
  that is, left modules over $\EE_1 \otimes \cO$-algebras bear a natural $\cO$-algebra structure.
\end{proposition}
We immediately acquire the following corollary, confirming a hypothesis of 
\cite[Rmk.~\href{https://arxiv.org/pdf/1708.03017v1\#theorem.3.15}{3.15}]{Hill_chromatic}.
\begin{corollary}
  Let $\cO^{\otimes}$ be a $\cT$-operad whose underlying $I^\infty$-operad is $\EE_\infty$ and $\cC^{\otimes}$ an $\cO$-monoidal $\infty$-category. 
  There is a lift
  \[
    \begin{tikzcd}[ampersand replacement=\&, column sep=huge]
      \&\& {\Cat_{\cO}^{\otimes}} \\
      {\Alg_{\cO}(\cC)} \& {\Alg_{\EE_1}(\cC)} \& \Cat
      \arrow[from=1-3, to=2-3]
      \arrow[from=2-1, to=1-3, "\RMod_{(-)}^{\otimes}(\cC)", curve={height=-12pt}]
      \arrow[from=2-1, to=2-2]
      \arrow[from=2-2, to=2-3, "\RMod_{(-)}(\cC)"]
    \end{tikzcd}
  \]
  natural separately in $\cO^{\otimes}$ and $\cC^{\otimes}$.
  In particular, if $I$ is an indexing category, the $\infty$-category of right-modules over an $I$-commutative algebra admits a natural $I$-symmetric monoidal structure.
\end{corollary}

\subsection{Equivariant infinitary Dunn additivity}\label{Infinitary Dunn additivity}
\def\uRep{\underline{\Rep}}
\def\orth{\mathrm{orth}}
\def\Emb{\mathrm{Emb}}
In \cite{Bonventre-nerve}, a \emph{genuine operadic nerve} 1-categorical functor was constructed between a model of graph-$G$ operads and a model of $G$-operads.
In \cite{EBV}, we lifted this to a conservative functor of $\infty$-categories $N^{\otimes} \cln g\Op_G \rightarrow \Op_G$.
Given $V$ an orthogonal $G$-representation, we define 
\[
  \EE^{\otimes}_V \deq N^{\otimes}D_V,
\]
where $D_V$ is the \emph{little $V$-disks graph $G$-operad} of \cite{Guillou-May}, whose $n$-ary $G \times \Sigma_n$ space is the configuration space
\[
    D_V(n) \simeq \mathrm{Conf}_n(V)
  \]
by \cite[Lem~1.2]{Guillou-May}.
The resulting unital $G$-operad $\EE_V$ was studied in \cite{Horev}, who showed for instance that
\[
    \EE_V(S) \simeq \mathrm{Conf}_{S}^H(V) \deq \dcolim_{\stackrel{W \subset V}{\mathrm{fin. dim}}} \mathrm{Conf}_S^H(W),
\]
in view of the fact that the assignment $\cO \mapsto \cO(S)$ preserves sifted colimits \cite[\S~2.3]{EBV};
here, $\Conf_S^H(W)$ is the space of $H$-equivariant configurations of $S$ into $W$ under the compact open topology.
\def\Conf{\mathrm{Conf}}
\def\stab{\mathrm{stab}}

A weak form of the following easy claim appears to be folklore.
\begin{proposition}\label{EV swindle}
  Let $G$ be a topological group, $H \subset G$ a closed subgroup, $S \in \FF_H$ a finite $H$-set admitting an configuration $\iota:S \hookrightarrow W$, and $V,W$ orthogonal $G$-representations whose associated map
  \[
    \Conf_S^H(V) \hookrightarrow \Conf_S^H(V \oplus W)
  \]
  is an equivalence.
  Then, $\Conf_S^H(V)$ is contractible.
\end{proposition}
\begin{proof}
    Linear interpolation to $\iota$ yields a deformation of $\Map^H(S,V \oplus W)$ onto $\cbr{\iota}$.
    The path of a point beginning in the subspace $\Conf_S^H(V) \subset \Conf_S^H(V \oplus W)$ consisting of configurations with zero projection to $W$ lands within $\Conf^H_S(V \oplus W)$ at all times;
    composing this deformation after the deformation retract $\Conf^H_S(V \oplus W) \xrightarrow\sim \Conf_H^S(V)$ yields a deformation retract of $\Conf_S^H(V \oplus W)$ onto $\cbr{\iota}$, so it is contractible.\footnote{
      Said explicitly, let $h\colon \brk{0,1} \rightarrow \Conf_S^H(V \oplus W)$ be the deformation retract onto those configurations with zero projection to $W$.
      Then, our deformation retract $h'$ onto $\iota(w)$ is computed by
      \[
        h'(t) = \begin{cases}
          h(2t)                                               & t \leq \frac{1}{2},\\
          \prn{2 - 2t} \cdot h(1) + \prn{2t - 1} \iota        & t \geq \frac{1}{2}.
        \end{cases}
      \]
    The second is an \emph{isotopy} since $h(1)$ and $\iota$ are pointwise-linearly independent embeddings.}
  By the equivalence $\Conf_S^H(V) \simeq \Conf_S^H(V \oplus W)$, the space $\Conf_S^H(V)$ is contractible as well.
\end{proof}

\begin{remark}
  This argument only produces \emph{contractibility}, whereas the nonequivariant argument using Fadell and Neuwirth's fibration \cite{Fadell} sharply characterizes \emph{$n$-connectivity} of $\Conf_k(\RR^n)$, and hence of $\EE_k^{\otimes}$;
  the author will equivariantize this in forthcoming work.
\end{remark}

We say that $V$ is a \emph{weak universe} if it is a direct sum of infinitely many copies of a collection of irreducible orthogonal $G$-representations;
equivalently, there is an equivalence $V \simeq V \oplus V$.
Given $V$ an orthogonal $G$-representation, we let $AV \deq A\EE_V$, i.e. $AV$ corresponds with the weak indexing system $\uFF^V = \uFF_{AV}$ of finite $H$-sets admitting an embedding into $V$.
The following corollary follows immediately from \cref{EV swindle}.

\begin{corollary}\label{EV is ninfty}
  If there exists an equivalence $\EE_V^{\otimes} \simeq \EE_{V \oplus W}^{\otimes}$, then the canonical map $\Bor_{AW}^{G} \EE_V^{\otimes} \rightarrow \cN_{AW}^{\otimes}$ is an equivalence;
  in particular, if $V$ is a weak universe, then there is a unique equivalence
  \[
    \EE_V^{\otimes} \xrightarrow{\;\;\; \sim \;\;\;} \cN_{AV}^{\otimes}.
  \]
\end{corollary}
\begin{observation}
  If $V$ is a \emph{universe} (i.e. it is a weak universe admitting a positive-dimensional fixed point locus), then it admits embeddings of all finite sets with trivial $G$-action;
  in this case, $\EE_V^{\otimes}$ is not just a weak $\cN_\infty$-operad, but an $\cN_\infty$-operad.
\end{observation}
Much study has been dedicated to the less general setting of \emph{universes};
for instance, Rubin has given a complete and simple characterization of those indexing systems (equivalently, transfer systems) occurring as the arity-support of an $\EE_V$-operad in \cite{Rubin_steiner} for $G$ abelian.

An inclusion $V \subset W$ yields a map of graph $G$-operads $D_V \rightarrow D_W$, hence a map $\EE^{\otimes}_V \rightarrow \EE_W^{\otimes}$.
This yields a map of weak indexing systems $\uFF^V \rightarrow \uFF^W$;
in \cite{Windex} we showed that this is additive, i.e.
\begin{equation}\label{Additivity equation}
  \uFF^V \vee \uFF^W = \uFF^{V \oplus W}.
\end{equation}

\begin{corollary}[Equivariant infinitary Dunn additivity]\label{Equivariant infinitary Dunn additivity}
  Let $G$ be a finite group and $V,W$ real orthogonal $G$-representations satisfying at least one of the following conditions:
  \begin{enumerate}[label={(\alph*)}]
    \item $V,W$ are weak $G$-universes, or
    \item the functoriality map $\EE_{V}^{\otimes} \rightarrow \EE_{V \oplus W}^{\otimes}$ is an equivalence.
  \end{enumerate}
  Then there is a canonical equivalence
  \[
    \EE_{V}^{\otimes} \obv \EE_W^{\otimes} \rightarrow \EE_{V \oplus W}^{\otimes};
  \]
  equivalently, for any $G$-symmetric monoidal category $\cC$, there are canonical equivalences\footnote{What we mean by ``canonical'' depends on the case;
  for case (a), there is a contractible space of equivalences, and for case (b), this equivalence comes from inverting arrows of the zigzag $\EE_{V \oplus W}^{\otimes} \simeq \EE_{V \oplus W}^{\otimes} \obv \triv_G^{\otimes} \xrightarrow{\id \otimes !} \EE_{V \oplus W}^{\otimes} \obv \EE_W^{\otimes} \xleftarrow{\iota \otimes \id} \EE_V^{\otimes} \obv \EE_W^{\otimes}$.}
  \[
    \Alg_{\EE_V} \uAlg_{\EE_W}^{\otimes}(\cC) \leftarrow \Alg_{\EE_{V \oplus W}}(\cC) \rightarrow \Alg_{\EE_W} \uAlg_{\EE_V}^{\otimes}(\cC).
  \]
\end{corollary}
\begin{proof}
  Given \cref{EV is ninfty}, case (a) follows from \cref{Main theorem,Additivity equation} and case (b) follows from \cref{Absorption main corollary}.
\end{proof}

\begin{remark}
  In \cite{Szczesny}, an ostensibly similar result to \cref{Equivariant infinitary Dunn additivity} is proved:
  given $D_{V}$ the \emph{little Disks graph $G$-operad}, Szczesny constructs a non-homotopical Boardman-Vogt tensor product $\otimes$ and a canonical map $D_V \otimes D_W \rightarrow D_{V \oplus W}$, which he shows to be a weak equivalence of graph $G$-operads in \cite[Thm~7.1]{Szczesny}.
  Neither this result nor \cref{Equivariant infinitary Dunn additivity} imply each other.
  
  On one hand, Szczesny's result concerns a tensor product with no known homotopical properties, so it is incomparable with results concerning $\infty$-categories of algebras defined by homotopy-coherent universal properties. 
  On the other hand, while \cref{Equivariant infinitary Dunn additivity} is homotopical, it only concerns cases where at least one of the representations induces $I$-symmetric monoidal $\infty$-categories of algebras whose indexed tensor products are indexed coproducts;
  this property will not be satisfied for any nontrivial indexed tensor products in the finite-dimensional case, so the range of representations in Szczesny's result is significantly larger.
  The author will address the general case in forthcoming work.
\end{remark}

\subsection{Norms on Real topological Hochschild and cyclic homology}\label{THR subsection}
\subsubsection{Factorization homology in general}
In classical algebra, there is a well-known tensor products of functors $F,G:\cC \rightarrow \cD$ using monoidal structure of $\cD$:
the \emph{pointwise tensor product} sets $F \otimes G(C) \deq F(C) \otimes G(C)$.
We will use a lift of this due to Nardin-Shah.
\begin{theorem}[{\cite[Thm~\href{https://arxiv.org/pdf/2203.00072v1\#nul.3.3.1}{3.3.1},Thm.~\href{https://arxiv.org/pdf/2203.00072v1\#nul.3.3.3}{3.3.3}]{Nardin}}]
  Let $\cK$ be a $\cT$-$\infty$-category and $\cC^{\otimes}$ a $\cT$-operad.
  Then, there exists a unique (functorial) $I$-operad structure $\Fun_{\cT}(\cK, \cC)^{\otimes-\ptw}$ on $\Fun_{\cT}(\cK,\cC)$ satisfying the universal property
  \[
    \Alg_{\cO}(\uFun_{\cT}(\cK,\cC)^{\otimes-\ptw}) \simeq \uFun_{\cT}(\cK, \uAlg_{\cO}(\cC))
  \]
  for $\cO \in \Op_I$.
  Furthermore, when $\cC^{\otimes}$ is $I$-symmetric monoidal, $\uFun_{\cT}(\cK,\cC)^{\otimes-\ptw}$ is $I$-symmetric monoidal and satisfies the universal property
  \[
    \Fun^{I-\otimes}_{\cT}\prn{\cD,\uFun_{\cT}(\cK,\cC)^{\otimes-\ptw}}
    \simeq 
    \Fun_{\cT}\prn{\cK, \uFun^{I-\otimes}_{\cT}(\cD,\cC)}.
  \]
  If additionally, $S$ is $I$-admissible, then the $S$-indexed tensor product of $(F_U) \in \uFun_{\cT}(\cK,\cC)^{\otimes-\ptw}_S$ has value
\[\begin{tikzcd}
	{\cD_V} & {\cD_S} & {\cC_S} & {\cC_V}
	\arrow["{\Delta^S}", from=1-1, to=1-2]
  \arrow["{\bigotimes\limits^S_U F_U}"{description}, curve={height=25pt}, from=1-1, to=1-4]
	\arrow["{(F_U)}", from=1-2, to=1-3]
  \arrow["{\bigotimes^S}", from=1-3, to=1-4]
\end{tikzcd}\]
\end{theorem}
The following proposition is easy, so we omit its proof.
\begin{proposition}
  There exists a natural equivalence $\uFun_{\cT}(\cK,\cC)^{\otimes-\ptw} \simeq \uAlg_{\triv(\cK)}^{\otimes}(\cC)$.
\end{proposition}
To use this, we need a number of functors involving these terms.
\begin{lemma}\label{U sym mon}
  Let $\cC^{\otimes}$ be a $\cT$-symmetric monoidal category, $\cO^{\otimes}$ a $\cT$-operad, and $\cK$ a $\cT$-category.
  \begin{enumerate}
    \item There is a natural $\cT$-symmetric monoidal functor
      $U\colon \uAlg_{\cO}^{\otimes}(\cC) \rightarrow \uFun_{\cT}(\Env(\cO),\cC)^{\otimes-\ptws}.$
    \item There is a lift of coevaluation of $\cT$-functors,
      \[\begin{tikzcd}[ampersand replacement=\&]
        \& { \uFun_G^{\otimes}\prn{\uFun_G\prn{\cK,\cC}^{\otimes-\ptws}, \cC}.} \\
        {\cK } \& { \uFun_G\prn{\uFun_G\prn{\cK,\cC}, \cC}.}
        \arrow["U", from=1-2, to=2-2]
        \arrow["\coev", from=2-1, to=1-2]
        \arrow["\coev"', from=2-1, to=2-2]
      \end{tikzcd}\]
      natural in the sense that varying $\cK$ forms a $\cT$-functor $\ucS_{\cT} \rightarrow \uFun_{\cT}(\Infl_e^{\cT} \partial \Delta^1, \uCat_{\cT})$. 
  \end{enumerate}
\end{lemma}
\begin{proof}
  (1) is pullback along the $\cT$-operad map $\triv(\Env(\cO))^{\otimes} \rightarrow \cO^{\otimes}$ adjunct to the structure functor $\Env(\cO) \rightarrow \cO$.
  For (2), this comes down to fussing with with adjunctions;
  the equivalence 
  \begin{align*}
    \Alg_{\triv(\cK)} \uAlg_{\uFun_G(\cK,\cC)^{\otimes-\ptws}}(\cC) 
    \simeq 
    \Alg_{\uFun_G(\cK,\cC)^{\otimes-\ptws}}\uFun_G(\cK,\cC)^{\otimes-\ptws} 
  \end{align*}
  carries $\Alg_{\triv(\cK)} \uFun_G^{\otimes}\prn{\uFun_G(\cK,\cC)^{\otimes-\ptws}, \cC}$ onto the $G$-symmetric monoidal functors.
  In particular, the identity pulls back to a $\cT$-operad map 
  \[
    \triv(\cK)^{\otimes} \rightarrow \uFun_G^{\otimes}\prn{\uFun_G(\cK,\cC)^{\otimes-\ptws},\cC}^{\otimes-\ptws},
  \]
  Now, note that we have an additional natural equivalence
  \[
    \Alg_{\triv(\cK)} \uFun_G\prn{\uFun_G(\cK,\cC),\cC}^{\otimes-\ptws} 
    \simeq 
    \Alg_{\triv\prn{\uFun_G(\cK,\cC)}}\uFun_G(\cK,\cC)^{\otimes-\ptws} 
  \]
  and precomposition with the counit $\triv\prn{\uFun_G(\cK,\cC)}^{\otimes} \rightarrow \uFun_G(\cK,\cC)^{\otimes-\ptws}$ yields a $\cT$-operad diagram 
  \[
    \begin{tikzcd}[ampersand replacement=\&]
	\& { \uFun_G^{\otimes}\prn{\uFun_G\prn{\cK,\cC}^{\otimes-\ptws}, \cC}^{\otimes-\ptws}} \\
	{\triv(\cK) } \& { \uFun_G\prn{\uFun_G\prn{\cK,\cC}^{\otimes-\ptws}, \cC}^{\otimes-\ptws}}
	\arrow["U", from=1-2, to=2-2]
	\arrow["\coev", from=2-1, to=1-2]
	\arrow["\coev"', from=2-1, to=2-2]
\end{tikzcd}
  \]
  The desired diagram is constructed from this using the adjunction $\triv \dashv U$.
\end{proof}
We use this to construct a $G$-symmetric monoidal lift for \emph{genuine equivariant factorization homology}.
\begin{corollary}\label{Symmetric monoidal factorization}
  Given $\cC^{\otimes}$ a distributive $G$-symmetric monoidal $\infty$-category, genuine equivariant factorization homology assembles to a $G$-functor
  \[
    \int\colon \uMfld^{V-fr,G} \rightarrow \uFun_G^{\otimes}\prn{\uAlg_{\EE_V}^{\otimes}(\cC), \cC},
  \]
  natural with respect to $G$-sifted $G$-colimit preserving $G$-symmetric monoidal functors in $\cC$.
  In particular, evaluation at $M$ yields a commutative diagram of $G$-symmetric monoidal functors
  \[
    \begin{tikzcd}[ampersand replacement=\&, column sep=large]
	{\uCAlg_{AV}^{\otimes}(\cC)} \& {\uCAlg_{AV}^{\otimes}(\cC)} \\
	{\uAlg^{\otimes}_{\EE_V}(\cC)} \& {\cC^{\otimes}}
	\arrow["{\int_M}", from=1-1, to=1-2]
	\arrow["U", from=1-1, to=2-1]
	\arrow["U", from=1-2, to=2-2]
	\arrow["{\int_M}", from=2-1, to=2-2]
\end{tikzcd}
\]
\end{corollary}
\begin{proof}
  In the notation of \cite{Horev}, let $\iota^{\otimes}:\uDisk^{G,V-fr,\sqcup} \rightarrow \uMfld^{G,V-fr,\sqcup}$ be the $G$-symmetric monoidal inclusion of $V$-framed $G$-disks into $V$-framed $G$-manifolds.
  By \hor{Prop.}{4.1.4}, $\int_M$ may be presented as the $G$-value of a composition
  \[
    \int_M\colon \uAlg_{\EE_V}(\cC) \simeq \uFun_{G}^{\otimes}\prn{\uDisk^{G,V-fr}, \cC} \xrightarrow{U} \uFun_G\prn{\uDisk^{G,V-fr},\cC} \xrightarrow{\iota_!} \uFun_G\prn{\uMfld^{G,V-fr}, \cC} \xrightarrow{\ev_M} \cC.
  \]
  To construct the desired $G$-functor $\int_M$, it suffices to construct 
  $G$-symmetric monoidal lifts of $U$ and $\iota_!$;
  then, we have a $G$-symmetric monoidal functor
  \[
    \begin{tikzcd}[ampersand replacement=\&, column sep=2em]
	{\uMfld^{V-fr,G}} \&\& {\uFun_G^{\otimes}\prn{\uAlg_{\EE_V}^{\otimes}(\cC),\cC}} \\
	{\uFun_G^{\otimes}\prn{\uFun_G\prn{\uMfld^{G,V-fr}, \cC},\cC}} \& {\uFun_G^{\otimes}\prn{\uFun_G\prn{\uDisk^{G,V-fr}, \cC},\cC}} \& {\uFun_G^{\otimes}\prn{\uFun_G^{\otimes}\prn{\uDisk^{G,V-fr}, \cC},\cC}}
	\arrow["\int", from=1-1, to=1-3]
  \arrow["\mathrm{coev}", from=1-1, to=2-1]
	\arrow["\simeq"{marking, allow upside down}, draw=none, from=1-3, to=2-3]
	\arrow["{U^*}", from=2-1, to=2-2]
	\arrow["{\prn{\iota_!}^*}", from=2-2, to=2-3]
  \end{tikzcd}
\]
  whose value on $M$ is $\int_M$.
  Indeed, $U$ is given by \cref{U sym mon}.

  For $\iota_!$, we use the $G$-symmetric monoidality of $G$-operadic left Kan extension argued in \cref{cor:T-sym-mon-kan}, noting that $G$-siftedness of the relevant slice $G$-category $\uDisk^{G,V-fr}_{/M}$ follows from \hor{Lem.}{5.2.7}.
\end{proof}

\subsubsection{Multiplication and norms on $\THR$}
We specialize \cref{Symmetric monoidal factorization} to $(V,M) = (\sigma,S^\sigma)$.
\begin{corollary}\label{THR corollary}
  Real topological Hochschild homology lifts to a $C_2$-symmetric monoidal functor
  \[
    \THR:\uAlg^{\otimes}_{\EE_\sigma}(\cC) \rightarrow \cC;
  \]
  in particular, $\THR$ lifts to a $C_2$-symmetric monoidal endofunctor
  \[
    \THR:\uAlg_{\EE_{V  + \infty\sigma}}^{\otimes}(\cC) \rightarrow \uAlg_{\EE_{V + \infty\sigma}}^{\otimes}(\cC).
  \]
  Given $A \in \Alg_{\EE_{V + \infty \sigma}}(\cC)$, there is an equivalence
  \[
    \THR(A) \simeq \ucolim_{S^{\sigma}} A,
  \]
  with colimit taken in $\Alg_{\EE_{V + \infty \sigma}}(\cC)$, naturally in $A$.
\end{corollary}
\begin{proof}
  The last sentence is the only part which does not follow immediately from combining Horev's facorization homology formula \hor{Rmk.}{7.1.2} with \cref{Equivariant infinitary Dunn additivity,Symmetric monoidal factorization}.
  It suffices to show the colimit property for $\cO \simeq \cO \otimes \EE_\sigma$-algebras whenever $\Bor_{A\sigma}^{\cT} \cO^{\otimes} \simeq \cN_{A\sigma \infty}^{\otimes}$, which holds for $\EE_{V + \infty\sigma}$ by \cref{EV swindle}.
  In any case, naturality of the dihedral bar construction together with the the Wirthm\"uller maps of \cref{Wirthmuller construction} yields a diagram
  \[
    \begin{tikzcd}[ampersand replacement=\&]
      \vdots \& {A^{\sqcup \mu_3}} \& {A^{\sqcup \mu_2}} \& A \& {\colim_{S^{\sigma}} A} \\
      \vdots \& {A^{\otimes \mu_3}} \& {A^{\otimes \mu_2}} \& A \& {\THR(A);}
      \arrow[shift left=3, from=1-1, to=1-2]
      \arrow[shift right=3, from=1-1, to=1-2]
      \arrow[shift left, from=1-1, to=1-2]
      \arrow[shift right, from=1-1, to=1-2]
      \arrow[shift left=2, from=1-2, to=1-3]
      \arrow[shift right=2, from=1-2, to=1-3]
      \arrow[from=1-2, to=1-3]
      \arrow[from=1-2, to=2-2]
      \arrow[shift left, from=1-3, to=1-4]
      \arrow[shift right, from=1-3, to=1-4]
      \arrow[from=1-3, to=2-3]
      \arrow[from=1-4, to=1-5]
      \arrow[from=1-4, to=2-4]
      \arrow[from=1-5, to=2-5, "\varphi"]
      \arrow[shift left, from=2-1, to=2-2]
      \arrow[shift right, from=2-1, to=2-2]
      \arrow[shift left=3, from=2-1, to=2-2]
      \arrow[shift right=3, from=2-1, to=2-2]
      \arrow[shift left=2, from=2-2, to=2-3]
      \arrow[shift right=2, from=2-2, to=2-3]
      \arrow[from=2-2, to=2-3]
      \arrow[shift left, from=2-3, to=2-4]
      \arrow[shift right, from=2-3, to=2-4]
      \arrow[from=2-4, to=2-5]
    \end{tikzcd}
  \]
  where $\mu_n$ is the $n$-element ``dihedral'' $C_2$-set, i.e. the unique $\sigma$-admissible $C_2$-set of size $n$, and each row is a geometric realization diagram.
  When the domain category is $A\sigma$-semiadditive, the vertical maps between the bar constructions are equivalences, so $\varphi$ is an equivalence.
  The result then follows by $A\sigma$-semiadditivity of $\cO$-algebras, as in \cref{Main theorem generalized}.
\end{proof}

\begin{remark}
  The computation $\THR(A) = \colim_{S^{\sigma}} A$ when $A$ is pulled back from a $C_2$-commutative algebra is not new;
  indeed, it appears as \qs{Rmk.}{5.4}.
  In fact, the ambiguity induced by the potential discrepancy between our construction $\uAlg_{\cO}^{\otimes}(\cC)$ and that of \ns{Thm.}{5.3.4} vanishes for the $I$-symmetric monoidal structure on $\CAlg_I(\cC)$ by applying \cref{Cocartesian rigidity main theorem generalized} in view of the fact that each are cocartesian \ns{Thm.}{5.3.9}.
  The new element of this identification is that the operation on $C_2$-commutative algebras is induced canonically from the operation on $\EE_\sigma$-algebras and that the colimit formulas need only an $\EE_{\infty\sigma}$-algebra structure.
\end{remark}

\begin{remark}
  In the above, we only needed $\cC$ to be an $\EE_\sigma$-monoidal $\infty$-category;
  however, to easily understand $\cO$-algebras, one ought to assume that $\cC$ is an $\cO$-monoidal $\infty$.category.
\end{remark}
Now, we can construct a circle action on $\THR$.
\def\Diff{\mathrm{Diff}}
\def\uDiff{\underline{\Diff}}
\def\uEmb{\underline{\Emb}}
\def\uMap{\underline{\Map}}
\begin{construction}\label{Action construction}
  Define $\uDiff^{V-fr}(M) \subset \uEmb^{V-fr}(M,M)$ to be the topological subspace of diffeomorphisms and embeddings of $M$ with conjugation $G$-action, considered as a (grouplike) $\EE_1$-$G$-space.
  Precomposing the functoriality of \cref{Symmetric monoidal factorization} along the action $B\uDiff^{V-fr} \rightarrow B\uAut_{\uMfld^{V-fr}}(M) \subset \uMfld^{G, V-fr}$ yields a $\uDiff^{V-fr}(M)$-action on $\int_M(-)$ through $G$-symmetric monoidal natural transformations, 
  where $B\uDiff^{V-fr}(M)$ is the unique connected $G$-space with $\Omega B \uDiff^{V-fr}(M) \simeq \uDiff^{V-fr}(M)$ as $\EE_1$-$G$-spaces (see \ha{\S}{5.2.6}).
  In particular, this yields a natural lift
\[\begin{tikzcd}[ampersand replacement=\&]
  \& {\uAlg_{\cO}(\cC)^{B\uDiff^{V-fr}(M)}} \\
	{\uAlg_{\cO \otimes \EE_V}(\cC)} \& {\cC.}
	\arrow["U", from=1-2, to=2-2]
	\arrow[dashed, from=2-1, to=1-2]
	\arrow["{\int_M}"', from=2-1, to=2-2]
\end{tikzcd}\]
  Applying the left-action $S^\sigma \rightarrow \uDiff^{\sigma-fr}(S^\sigma)$ yields a lift
  \[\begin{tikzcd}[ampersand replacement=\&]
    \& {\uAlg_{\cO}(\cC)^{BS^{\sigma}}} \\
	{\uAlg_{\cO \otimes \EE_\sigma}(\cC)} \& \cC
	\arrow["U", from=1-2, to=2-2]
	\arrow[dashed, from=2-1, to=1-2]
	\arrow["\THR"', from=2-1, to=2-2]
\end{tikzcd}\]
  We refer to this as the \emph{$S^\sigma$ action on $\THR$}.
\end{construction}
In fact, this type of action has been seen in previous work.
\begin{remark}
  For $\psi\colon G \rightarrow G/N$ a surjective topological group homomorphism, the \emph{$N$-free $G$-family} is
  \[
    B_{G/N}^{\psi} N = \cbr{[G / H] \;\;\middle|\;\; \Res_N^G([G/H]) \text{ is a free } N \text{-set}} \subset \cO^{\op}_G;
  \]
  \qs{Lem.}{4.14} recognizes the ``quotient by $N$'' forgetful functor
  \[
    B_{G/N}^{\psi} N \rightarrow \cO^{\op}_{G/N},
  \]
  as the unstraightening of a $G/N$-space, which we also refer to as $B_{G/N}^{\psi} N$.
  In particular, for an arbitrary Abelian group $A$, they define $B_{C_2}^t A$ with respect to the semidirect product extension $A \rightarrow A \rtimes C_2 \rightarrow C_2$ for the $C_2$-action on $A$ by inversion.
  In the case that $A = \TT$ is the circle group, this is given by the usual extension $\TT = \SO(2) \rightarrow O(2) \rightarrow C_2$.

  We may explicitly identify $B_{C_2}^t \TT$ by hand;
  it's well-known that there are exactly two conjugacy classes of subgroups $(H) \subset O(2)$ whose homogeneous spaces $[O(2)/H]$ are $\TT$-free:
  the reflections and the trivial subgroup.
  The normalizer of a subgroup generated by a reflection is a dihedral group of order 4, so
  \[
    \Aut_{O(2)}([O(2) / C_2]) \simeq W_{O(2)}(C_2) \simeq D_4/C_2 \simeq C_2.
  \]
  In particular, we may picture the unstraightening of $B_{C_2}^t \TT$ via the diagram
  \[\begin{tikzcd}[ampersand replacement=\&]
	{B^t_{C_2}\TT} \&\&\&\&\& {[O(2)/C_2]} \& {[O(2)/e]} \\
	{*_{C_2}} \&\&\&\&\& {[C_2/C_2]} \& {[C_2/e]}
  \arrow["{\mathrm{Unstraightening}}"{pos=0.45}, shorten <=15pt, shorten >=30pt, squiggly={pre length = 15pt, post length=30pt}, tail reversed, from=1-1, to=1-6]
	\arrow[from=1-1, to=2-1]
	\arrow["{C_2}"', from=1-6, to=1-6, loop, in=215, out=145, distance=10mm]
	\arrow[""{name=0, anchor=center, inner sep=0}, from=1-7, to=1-6]
	\arrow["{O(2)}", from=1-7, to=1-7, loop, in=325, out=35, distance=10mm]
  \arrow[shorten <=15pt, shorten >=30pt, squiggly={pre length = 15pt, post length=30pt}, tail reversed, from=2-1, to=2-6]
	\arrow[""{name=1, anchor=center, inner sep=0}, from=2-7, to=2-6]
	\arrow["{C_2}", from=2-7, to=2-7, loop, in=325, out=35, distance=10mm]
	\arrow[shorten <=9pt, shorten >=6pt, Rightarrow, from=0, to=1]
\end{tikzcd}\]
In particular, taking fibers, we find that $\prn{\Omega B^t_{C_2} \TT}^e = \TT$, $\prn{\Omega B^t_{C_2} S^1}^{C_2} = C_2$, and the restriction map $C_2 \rightarrow \TT$ is the unique nontrivial such homomorphism.
This evidently agrees with the circle group structure on $S^\sigma$, inducing an equivalence $BS^\sigma \simeq B_{C_2}^t \TT$, so we really have acquired a natural lift
  \[\begin{tikzcd}[ampersand replacement=\&]
    \& {\uAlg_{\cO}(\cC)^{h_{C_2} \TT}} \\
	{\uAlg_{\cO \otimes \EE_\sigma}(\cC)} \& \cC
	\arrow["U", from=1-2, to=2-2]
	\arrow[dashed, from=2-1, to=1-2]
	\arrow["\THR"', from=2-1, to=2-2]
\end{tikzcd}\]
where we write $\cC^{h_{C_2} \TT} \simeq \uFun_{C_2}\prn{B^t_{C_2} \TT, \cC}$.
\end{remark}
The following reformulation of the equivalence $\THR(A) \simeq \colim_{S^\sigma} A$ may be familiar.
\begin{observation}\label{Free action obs}
  Let $\Gamma$ be a grouplike $\EE_1$-$\cT$-space.
  Evaluation and left Kan extension yield an adjunction
  \[
    - \otimes \Gamma\colon \Gamma^{\cT} \cC \rightleftarrows \Fun_\cT(B \Gamma, \cC)\colon U;
  \]
  unwinding the left Kan extension formula along $*_{\cT} \rightarrow B \Gamma$ shows that the $\cT$-object underlying $- \otimes \Gamma$ is the constant $\Gamma = \Omega B\Gamma$-indexed $\cT$-colimit functor;
  in particular, \cref{THR corollary} understands $\THR(A)$ to be the free $\EE_{V + \infty \sigma}$-object on $A$ with $C_2$-equivariant $S^{\sigma}$-action.
  This free action agrees with that of \cref{Action construction}, and in particular this identifies our action with that of \qstwosec{5}.
\end{observation}
 
\subsubsection{Multiplication and norms on $\TCR$}
Having produced a $C_2$-symmetric monoidal construction 
\[
\THR\colon \Alg_{\EE_\sigma}\prn{\Sp_{C_2}} \rightarrow \uSp_{B_{C_2} \TT}^{\otimes} \deq \uFun_{C_2}\prn{B^t_{C_2} \TT, \uSp_{C_2}^{\otimes}}^{\otimes-\ptws}
\]
which lifts to Quigley-Shah's construction, we're poised to become the ``future work'' indicated in \qs{Warning}{0.12}
by constructing a lax symmetric monoidal $p$-typical (Borel) Real topological cyclic homology functor which lifts to Quigley-Shah's construction.
Now, consider the $C_2$-space map $i\colon B_{C_2}^t \mu_{p^\infty} \rightarrow B_{C_2}^t \TT$ defined by the following colimit
\[\begin{tikzcd}[ampersand replacement=\&, row sep=small]
  	{B^t_{C_2} e} \& \cdots \& {B_{C_2}^t \mu_{p^n}} \& {B_{C_2}^t \mu_{p^{n+1}}} \& \cdots \& {B_{C_2}^{t} \mu_{p^\infty}} \\
	\&\&\&\&\& {B_{C_2}^t \TT}
	\arrow[hook', from=1-1, to=1-2]
	\arrow[curve={height=24pt}, from=1-1, to=2-6]
	\arrow[hook', from=1-2, to=1-3]
	\arrow[hook', from=1-3, to=1-4]
	\arrow[curve={height=12pt}, from=1-3, to=2-6]
	\arrow[hook', from=1-4, to=1-5]
	\arrow[curve={height=6pt}, from=1-4, to=2-6]
	\arrow[hook', from=1-5, to=1-6]
	\arrow[from=1-6, to=2-6, "i"]
\end{tikzcd}\]
Postcomposing with $i^*$ yields a $C_2$-symmetric monoidal functor $\uAlg_{\EE_\sigma}^{\otimes}(\Sp_{C_2}) \rightarrow \uSp_{B_{C_2} \mu_{p^\infty}}^\otimes$.

Following Nikoulaus and Scholze \cite[Cons~\href{https://arxiv.org/pdf/1707.01799v2\#thm.4.2.1}{IV.2.1}]{Nikolaus} Quigley-Shah defined the \emph{$C_2$-operad of (Borel) Real $p$-cyclotomic spectra} by the pullback $C_2$-operad 
  \[\begin{tikzcd}[ampersand replacement=\&]
    \uRCyc^{\otimes}_p \& {\uFun_{C_2}\prn{\Infl_e^{C_2} \Delta^1, \uSp_{C_2}^{\otimes}}^{\otimes-\ptws}} \\
  {\uSp_{B_{C_2}^t \mu_{p^\infty}}^{\otimes}} \& {\uSp_{C_2}^{\otimes} \times \uSp_{C_2}^{\otimes}}
	\arrow[from=1-1, to=1-2]
	\arrow[from=1-1, to=2-1]
	\arrow["\lrcorner"{anchor=center, pos=0.125}, draw=none, from=1-1, to=2-2]
	\arrow["{(\ev_0, \ev_1)}", from=1-2, to=2-2]
	\arrow["{\prn{U, (-)^{t_{C_2} \mu_p}}}"', from=2-1, to=2-2]
\end{tikzcd}\]
where $(-)^{t_{C_2}\mu_p}$ has the lax $C_2$-symmetric monoidal structure of \qs{Setup}{2.1};
this is a \emph{lax equalizer of $C_2$-operads}.
Pullback-stability of cocartesian fibrations guarantees that $\uRCyc^{\otimes} \rightarrow \uSp_{C_2}^{\otimes}$ is a cocartesian fibration of $C_2$-operads, and in particular, $\uRCyc^{\otimes}$ is a $C_2$-symmetric monoidal $\infty$-category.

In \qstwo{Const.}{5.5}, Quigley-Shah define a lax $C_2$-symmetric monoidally natural \emph{dihedral Tate diagonal} functor $\Delta_p\colon \uSp_{C_2}^{\otimes} \rightarrow \uFun_{C_2}\prn{\Infl_e^{C_2} \Delta^1, \uSp_{C_2}^{\otimes}}^{\otimes-\ptws}$
whose composites are $\ev_1 \circ \Delta_p A \sim \prn{A^{\otimes \mu_p}}^{t_{C_2} \mu_p}$ and $\ev_0 \sim \id$.
\begin{remark}
  Note that the full subcategory $\FF^{\sigma}_{C_2} = \cbr{\mu_{n} \mid n \in \NN} \subset \FF_{C_2}$ of $\sigma$-admissible $C_2$-sets participates in a weak indexing system with $\FF^\sigma_e = \FF$;
  in particular these together are closed under indexed coproducts.
  One can easily identify the underlying set of output by the \emph{double coset formula} for indexed coproducts, which in this case simply shows that $\abs{\coprod_U^S X_U} = \sum_{U \in \Orb(S)} \abs{U} \cdot \abs{X_U}$;
  since $\Res_e^{C_2}\colon \FF^{\sigma}_{C_2} \rightarrow \FF^{\sigma}_e = \FF$ is injective, the dihedral Tate diagonal of an indexed tensor power has signature $X^{\otimes \mu_n} \longrightarrow  {\prn{\prn{X^{\otimes \mu_n}}^{\otimes \mu_p}}^{t_{C_2}\mu_p}} \simeq \prn{X^{\otimes \mu_{np}}}^{t_{C_2} \mu_p}$.
\end{remark}

\begin{construction}[Dihedral Frobenius]
  Naturality of the dihedral Tate diagonal, commutativity of $C_2$-colimits with $C_2$-colimits, and the coassembly map for $C_2$-limits yields a natural diagram
  \[\begin{tikzcd}[ampersand replacement=\&, column sep=0.9em]
	\cdots \& {X^{\otimes \mu_3}} \&\& {X^{\otimes \mu_2}} \&\& X \&\&\& {\THR(X)} \\
  \cdots \& {\prn{X^{\otimes \mu_{3p}}}^{t_{C_2}{mu_p}}} \&\& {\prn{X^{\otimes \mu_{2p}}}^{t_{C_2}{\mu_p}}} \&\& {\prn{X^{\otimes \mu_{p}}}^{t_{C_2}{\mu_p}}} \&\&\& {\colim \prn{X^{\otimes \mu_p}}^{t_{C_2}\mu_p}} \\
  {\Bigg(\cdots} \& {X^{\otimes \mu_{3p}}} \&\& {X^{\otimes \mu_{2p}}} \&\& {X^{\otimes \mu_{p}}\Bigg)^{t_{C_2} \mu_p}} \&\&\& {\prn{\colim \prn{X^{\otimes \mu_{p}}}}^{t_{C_2}\mu_p}} \\
  {\Bigg(\cdots} \& {X^{\otimes \mu_{3p}}} \& \cdots \& {X^{\otimes \mu_{2p}}} \& \cdots \& {X^{\otimes \mu_{p}}} \& \cdots \& {X\Bigg)^{t_{C_2}\mu_p}} \& {\THR(X)^{t_{C_2}\mu_p}}
	\arrow[shift left=3, from=1-1, to=1-2]
	\arrow[shift right=3, from=1-1, to=1-2]
	\arrow[shift left, from=1-1, to=1-2]
	\arrow[shift right, from=1-1, to=1-2]
	\arrow[from=1-2, to=1-4]
	\arrow[shift left=2, from=1-2, to=1-4]
	\arrow[shift right=2, from=1-2, to=1-4]
	\arrow["{{\Delta_p(X^{\otimes 3})}}", from=1-2, to=2-2]
	\arrow[shift left, from=1-4, to=1-6]
	\arrow[shift right, from=1-4, to=1-6]
	\arrow["{{\Delta_p(X^{\otimes 2})}}", from=1-4, to=2-4]
	\arrow[from=1-6, to=1-9]
	\arrow["{{\Delta_p(X)}}", from=1-6, to=2-6]
	\arrow[from=1-9, to=2-9]
	\arrow[shift left, from=2-1, to=2-2]
	\arrow[shift right, from=2-1, to=2-2]
	\arrow[shift left=3, from=2-1, to=2-2]
	\arrow[shift right=3, from=2-1, to=2-2]
	\arrow[from=2-2, to=2-4]
	\arrow[shift left=2, from=2-2, to=2-4]
	\arrow[shift right=2, from=2-2, to=2-4]
	\arrow[shift right, from=2-4, to=2-6]
	\arrow[shift left, from=2-4, to=2-6]
	\arrow[from=2-6, to=2-9]
	\arrow[from=2-9, to=3-9]
	\arrow[shift left, from=3-1, to=3-2]
	\arrow[shift right, from=3-1, to=3-2]
	\arrow[shift left=3, from=3-1, to=3-2]
	\arrow[shift right=3, from=3-1, to=3-2]
	\arrow[shift left=2, from=3-2, to=3-4]
	\arrow[shift right=2, from=3-2, to=3-4]
	\arrow[from=3-2, to=3-4]
	\arrow[from=3-2, to=4-2]
	\arrow[shift left, from=3-4, to=3-6]
	\arrow[shift right, from=3-4, to=3-6]
	\arrow[from=3-4, to=4-4]
	\arrow[from=3-6, to=3-9]
	\arrow[from=3-6, to=4-6]
	\arrow[from=3-9, to=4-9]
	\arrow[""{name=0, anchor=center, inner sep=0}, shift left=3, from=4-1, to=4-2]
	\arrow[""{name=1, anchor=center, inner sep=0}, shift right=3, from=4-1, to=4-2]
	\arrow[""{name=2, anchor=center, inner sep=0}, shift left=3, from=4-2, to=4-3]
	\arrow[""{name=3, anchor=center, inner sep=0}, shift right=3, from=4-2, to=4-3]
	\arrow[""{name=4, anchor=center, inner sep=0}, shift left=3, from=4-3, to=4-4]
	\arrow[""{name=5, anchor=center, inner sep=0}, shift right=3, from=4-3, to=4-4]
	\arrow[""{name=6, anchor=center, inner sep=0}, shift left=3, from=4-4, to=4-5]
	\arrow[""{name=7, anchor=center, inner sep=0}, shift right=3, from=4-4, to=4-5]
	\arrow[""{name=8, anchor=center, inner sep=0}, shift left=3, from=4-5, to=4-6]
	\arrow[""{name=9, anchor=center, inner sep=0}, shift right=3, from=4-5, to=4-6]
	\arrow[""{name=10, anchor=center, inner sep=0}, shift left=3, from=4-6, to=4-7]
	\arrow[""{name=11, anchor=center, inner sep=0}, shift right=3, from=4-6, to=4-7]
	\arrow[shift right, from=4-7, to=4-8]
	\arrow[shift left, from=4-7, to=4-8]
	\arrow[from=4-8, to=4-9]
	\arrow["\cdots"{marking, allow upside down}, draw=none, from=0, to=1]
	\arrow["\cdots"{marking, allow upside down}, draw=none, from=2, to=3]
	\arrow["\cdots"{marking, allow upside down}, draw=none, from=4, to=5]
	\arrow["\cdots"{marking, allow upside down}, draw=none, from=6, to=7]
	\arrow["\cdots"{marking, allow upside down}, draw=none, from=8, to=9]
	\arrow["\cdots"{marking, allow upside down}, draw=none, from=10, to=11]
\end{tikzcd}\]
We acquire a (lax $C_2$-symmetric monoidally-) natural composite map $\varphi_p\colon \THR(X) \rightarrow \THR(X)^{t_{C_2} \mu_p}$, 
yielding a lax $C_2$-symmetric monoidal functor
\[\begin{tikzcd}[ampersand replacement=\&, column sep=huge]
	{\uAlg_{\EE_\sigma}^{\otimes}\prn{\uSp_{C_2}}} \\
	\& {\uRCyc^{\otimes}_p} \& {\uFun_{C_2}\prn{\Infl_e^{C_2} \Delta^1, \uSp_{C_2}^{\otimes}}^{\otimes-\ptws}} \\
  \& {\uSp_{B_{C_2}^t \mu_{p^\infty}}^{\otimes}} \& {\uSp_{C_2}^{\otimes} \times \uSp_{C_2}^{\otimes}}
	\arrow["{\widetilde\THR}"{description}, dashed, from=1-1, to=2-2]
	\arrow["{\varphi_p}", curve={height=-12pt}, from=1-1, to=2-3]
	\arrow["\THR"', curve={height=12pt}, from=1-1, to=3-2]
	\arrow[from=2-2, to=2-3]
	\arrow[from=2-2, to=3-2]
	\arrow["\lrcorner"{anchor=center, pos=0.125}, draw=none, from=2-2, to=3-3]
	\arrow["{(\ev_0, \ev_1)}", from=2-3, to=3-3]
	\arrow["{\prn{U, (-)^{t_{C_2} \mu_p}}}"', from=3-2, to=3-3]
\end{tikzcd}\qedhere\]
\end{construction}

Now, in \qssec{2.3} they also define a lax $C_2$-symmetric monoidal functor $\TCR(-,p)\colon \uRCyc_p^{\otimes} \rightarrow \uSp_{C_2}^{\otimes}$.
From this, we may conclude the following.

\begin{corollary}\label{TCR corollary}
  Given $\cO^{\otimes}$ a $C_2$-operad, there is a commutative diagram of lax $C_2$-symmetric monoidal functors
\[\begin{tikzcd}[ampersand replacement=\&, row sep=small, column sep=huge]
	{\uCAlg_{C_2}^{\otimes}\prn{\uSp_{C_2}}} \& {\uCAlg^{\otimes}_{C_2} \prn{\uRCyc_{p}}} \& {\uCAlg_{C_2}\prn{\uSp_{C_2}}} \\
	{\uAlg_{\EE_\sigma \otimes \cO}^{\otimes}\prn{\uSp_{C_2}}} \& {\uAlg^{\otimes}_{\cO} \prn{\uRCyc_{p}}} \& {\uAlg^{\otimes}_{\cO} \prn{\uSp_{C_2}}} \\
	{\uAlg_{\EE_\sigma}^{\otimes}\prn{\uSp_{C_2}}} \& {\uRCyc_{p}^{\otimes}} \& {\uSp_{C_2}^{\otimes}} \\
	\& {\uSp_{C_2}^{\otimes}}
	\arrow["{{\THR}}", from=1-1, to=1-2]
	\arrow[from=1-1, to=2-1]
  \arrow[from=1-2, to=1-3, "{\TCR(-,p)}"]
	\arrow[from=1-2, to=2-2]
	\arrow[from=1-3, to=2-3]
	\arrow[dashed, from=2-1, to=2-2]
	\arrow[from=2-1, to=3-1]
	\arrow[from=2-2, to=2-3]
	\arrow[from=2-2, to=3-2]
	\arrow[from=2-3, to=3-3]
	\arrow[dashed, from=3-1, to=3-2]
	\arrow["\THR"', from=3-1, to=4-2]
  \arrow["{\TCR(-,p)}", from=3-2, to=3-3]
	\arrow[from=3-2, to=4-2]
\end{tikzcd}\]
  whose top row recovers the constructions of \cite{Quigley}.
  In particular, if $I$ is a $C_2$-weak indexing system containing the dihedral $C_2$-sets,
  then Quigley-Shah's $p$-typical Real topological cyclic homology is lifted from a lax $C_2$-symmetric monoidal endofunctor of $I$-commutative ring spectra.
\end{corollary}

\begin{remark}
  The conditions of \cref{THR corollary,TCR corollary} were varied for the sake of diversity of examples;
  they can be interchanged, and each can be weakened to construct (lax) $C_2$-symmetric monoidal endofunctors of $\cO$-algebras whenever $\Bor_{A\sigma}^{C_2} \cO^{\otimes} \simeq \EE_{\sigma \infty}^{\otimes} \simeq \cN_{A\sigma \infty}^{\otimes}$.
\end{remark}

\begin{remark}
  For simplicity, we made the \emph{$p$-typical} construction;
  however, Real cyclotomic spectra and their Real topological cyclic homology were constructed integrally in \cite{Quigley_cyclotomic}, and we may lift our constructions to a lax $C_2$-symmetric monoidal Real cyclotomic structure as follows.
\[\begin{tikzcd}[ampersand replacement=\&, column sep=huge]
	{\uAlg_{\EE_\sigma}^{\otimes}\prn{\uSp_{C_2}}} \\
  \& {\uRCyc^{\otimes}} \& {\prod\limits_{p \text{ prime}}\uFun_{C_2}\prn{\Infl_e^{C_2} \Delta^1, \uSp_{B_{C_2}^t \TT}^{\otimes}}^{\otimes-\ptws}} \\
  \& {\uSp_{B_{C_2}^t \TT}^{\otimes}} \& {\prod\limits_{p \text{ prime}} \uSp_{B_{C_2}^t \TT} \times \uSp_{B_{C_2}^t \TT}}
	\arrow["{{\widetilde\THR}}"{description}, dashed, from=1-1, to=2-2]
	\arrow["{\prn{\varphi_p}}", curve={height=-12pt}, from=1-1, to=2-3]
	\arrow["\THR"', curve={height=12pt}, from=1-1, to=3-2]
	\arrow[from=2-2, to=2-3]
	\arrow[from=2-2, to=3-2]
	\arrow["\lrcorner"{anchor=center, pos=0.125}, draw=none, from=2-2, to=3-3]
	\arrow["{{(\ev_0, \ev_1)}}", from=2-3, to=3-3]
	\arrow["{\prn{U, (-)^{t_{C_2} \mu_p}}_p}"', from=3-2, to=3-3]
\end{tikzcd}\]
The lax $C_2$-symmetric monoidal \emph{Real topological Cyclic homology} functor may be defined to be the composite
\[
  \uAlg_{\EE_\sigma}^{\otimes}\prn{\uSp_{C_2}} \xrightarrow{\;\;\;\; \widetilde \THR \;\;\;\;} \uRCyc^{\otimes} \xrightarrow{\;\;\;\; \TCR \;\;\;\;} \uSp_{C_2}^{\otimes}.\qedhere
\]
\end{remark}
\subsubsection{Speculations on genuine equivariant $\THH$ for other groups}
Merling posed the following question.
\begin{question}[{\cite[Prob~1.6]{AIM_SHT}}]\label{Merling question}
  Is it possible to build a version of $\THH$ for $G$-ring spectra which is $G$-symmetric monoidal?
\end{question}
On its face, \cref{Merling question} receives a positive answer by setting $\cC = \uSp_G^{\otimes}$ in the following.
\begin{corollary}
  Given $\cC^{\otimes}$ a $G$-symmetric monoidal $\infty$-category, there is a $G$-symmetric monoidal functor
  \[
    \int_{S^1}\colon \uAlg^{\otimes}_{\EE_1}(\cC) \rightarrow \cC^{\otimes}
  \]
  whose $H$-value functor $\Alg_{\EE_1}(\cC_H) \rightarrow \cC_H$ is $\THH$.
\end{corollary}
\begin{proof}
  Specialize \cref{Symmetric monoidal factorization} to $M  = S^1$.
\end{proof}

This is possibly unsatisfying;
it certainly recovers the notion of \cite{Mehrle_Koszul}, but it recovers neither twisted nor Real $\THH$.
We can give a more general result recovering the latter as follows.
\begin{corollary}
  Let $\cC^{\otimes}$ be a $G$-symmetric monoidal $\infty$-category, $V$ be an orthogonal $G$-representation, and $S \subset V$ an embedded $G$-set.
  Then, there exists a $G$-symmetric monoidal functor
  \[
    \THH_{S,V}\colon \uAlg_{\EE_V}^{\otimes}(\cC) \rightarrow \cC^{\otimes},
  \]
  natural in $G$-symmetric monoidal strongly $G$-colimit preserving $G$-functors in $\cC$, satisfying
  \[
    \THH_{S,V}(A) \simeq 
    A \otimes_{A^{\otimes S}} A.
  \]
\end{corollary}
\begin{proof}
  Define the invariant closed topological subspace 
  \[
    S^{\spoke(S,V)} \deq \RR S \cup \cbr{\infty} \subset S^V.
  \]
  Giving $S^V$ the standard metric, pick some small $\varepsilon > 0$ and let $\tau S^{\spoke\prn{S,V}} \subset S^V$ be the open ball around $S^{\spoke(S,V)}$ of radius $\varepsilon$;
  this is canonically a $V$-framed open $G$-submanifold of $S^V$.
  We define 
  \[
    \THH_{S,V}(A) \deq \int_{\tau S^{\spoke(S,V)}} A\colon \uAlg_{\EE_V}^{\otimes}(\cC) \rightarrow \cC^{\otimes}.
  \]
  We're left with verifying the tensor product formula, which itself will be an application of Horev's $G$-$\otimes$-excision result \hor{Prop.}{5.2.3}.
  We need a collar decomposition;
  define the function $f\colon S^V \rightarrow [-1,1]$ by
  \[
    f(v) = \begin{cases}
      -1 & \abs{v} < 2\varepsilon; \\ 
      \frac{\abs{V} - 3\varepsilon}{\varepsilon} & 2\varepsilon \leq \abs{v} \leq 4\varepsilon;\\
      1 & \abs{v}  > 4 \varepsilon.
    \end{cases}
  \]
  The restriction of this to $\tau S^{\spoke(S,V)}$ is a collar decomposition with positive and negative part framed-diffeomorphic to $D(V)$ and with interior $V$-framed diffeomorphic to the indexed disjoint union $S \cdot D(V)$.
  $\otimes$-excision yields
  \begin{align*}
    \THH_{S,V}(A)
    &\simeq 
    \int_{\tau S^{\spoke(S,V)}_+} A
    \otimes_{\int_{\tau S^{\spoke(S,V)}_+ \cap \tau S^{\spoke(S,V)}} A}
    \int_{\tau S^{\spoke(S,V)}_-} A;\\
    &\simeq 
    \int_{D(V)} A
    \otimes_{\int_{S \cdot D(V)} A}
    \int_{D(V)^{\op}} A\\
    &\simeq A \otimes_{A^{\otimes S}} A^{\op}.\qedhere
  \end{align*}
\end{proof}

\begin{appendix}
\section{Technicalities on (co)cartesian \tI-symmetric monoidal \tinfty-categories}\label{Cocartesian proof subsubsection}
\stoptocwriting
For the duration of this appendix, we assume the notation of \cite[\S~A]{EBV}. 
Fix $\cP \subset \cT$ an atomic orbital subcategory and $I \subset \FF_{\cT}^{\cP}$ an almost-unital weak indexing category.
We define the \emph{$\infty$-category of $\Gamma$-$I$-preoperads} 
\[
  \PreOp_I^{\Gamma} \deq \Cat_{\cT, /\uFF_{I,*}}^{\Int-\cocart},
\]
so that the results of \cite{Barkan} recognize $\Op_I \subset \PreOp_I^{\Gamma}$ as a localizing subcategory.\footnote{Here, $\Gamma$ is a reference to Segal's category $\Gamma$, whereas the undecorated version centers the \emph{effective Burnside category}.}

This appendix can be understood as a lift of \cite[\S~2.4.1-2.4.3]{HA} to the setting of (co)cartesian $I$-symmetric monoidal $\infty$-categories, working in the specific model of $\Gamma$-$I$-preoperads;
we proceed by an essentially similar strategy, complicated only by less convenient combinatorics.
We suggest only readers in need of the minutiae inquire within, and the shunt remaining readers to the summaries contained in \cref{Cartesian subsection}.

First, define the $\cT$-1-category $\uGamma_I^*$ to have $V$-values
\[
  \Gamma_{I,V}^* \deq \cbr{U_+ \xrightarrow{s.i.} S_+ \;\; \middle| \;\; U \in \cT_{/V}} \subset \Ar(\uFF_{I,*})_V;
\]
that is, the objects of $\Gamma^*_{I,V}$ are pointed $I$-admissible $V$-sets with a distinguished orbit, and the morphisms of $\Gamma_{I,V}$ preserve distinguished orbits.
This possesses a \emph{target} forgetful $\cT$-functor $t\colon \uGamma_{I}^* \rightarrow \uFF_{I,*}$.
We use this to construct an $\infty$-category $\Tot \Tot_{\cT} \cC$ over $\Tot \uFF_{I,*}$ in \cref{Quasicat subsection} satisfying the following universal property.
\begin{proposition}\label{UFI exists prop}
  Given $\cC$ a $\cT$-$\infty$-category, there exists an $\infty$-category $\Tot \Tot_{\cT} \cC^{I-\sqcup}$ over $\Tot \uFF_{I,*}$ satisfying the universal property that there is a natural equivalence
  \[
    \Fun_{/\Tot \uFF_{I,*}}(\cD,\Tot \Tot_{\cT} \cC^{I-\sqcup}) \simeq \Fun_{/\cT^{\op}}(\cD \times_{\Tot \uFF_{I,*}} \Tot \uGamma_{I}^{*},\Tot \cC);
  \]
\end{proposition}

Second, define the (non-full) $\cT$-subcategory $\Gamma^\times_{I} \subset \Ar(\uFF_{I,*})$ to have $V$-objects given by summand inclusions of pointed $V$-sets $\overline{S}_+ \hookrightarrow S_+$ and morphisms of $V$-objects given by maps $\alpha\colon S_+ \rightarrow T_+$ with the property that $\alpha^{-1}\prn{\overline{T}_+} \subset \overline{S}_+$. 
In \cref{Quasicat subsection} we prove the following.
\begin{proposition}\label{UFI cart exists prop}
  Given $\cC$ a $\cT$-$\infty$-category, there exists an $\infty$-category $\Tot \Tot_{\cT} \widetilde \cC^{I-\times}$ over $\uFF_{I,*}$ satisfying the universal property that there is a natural equivalence
  \[
    \Fun_{/\Tot \uFF_{I,*}}(K,\Tot \Tot_{\cT} \widetilde \cC^{I-\times}) \simeq \Fun_{/\cT^{\op}}(K \times_{\Tot \uFF_{I,*}}\Tot \uGamma_I^\times, \Tot \cC).
  \]
\end{proposition}
Note that there is an equivalence
\[
  \cbr{S_+} \times_{\uFF_{I,*}} \uGamma_{I}^{\times} \simeq \sP_{\uV}(S),
\]
where $\sP_{\uV}(S)$ is the $V$-poset with $U$-value given by subsets of $\Res_U^V S$ ordered under inclusion.
In particular, for $S_+ \in \uFF_{I,*}$, we view objects in $\widetilde \cC^{I-\times}_{S_+}$ as $V$-functors $\sP_{\uV}(S)^{\op} \rightarrow \cC_{V}$.
Let $\Tot \Tot_{\cT} \cC^{I-\times} \subset \Tot \Tot_{\cT} \widetilde \cC^{I-\times}$ be the full subcategory whose objects over $V$ are spanned by those functors $F\colon \sP_{\uV}(S)^{\op} \rightarrow \cC_{\uV}$ satisfying the property that, for all $U \rightarrow V$ and $T \subset \Res_U^V S$, the maps $\Res_W^V F(T) \rightarrow F(W)$ exhibit $F(T)$ as the $T$-indexed product $F(T) \simeq \prod_W^T F(U)$ in $\cC$.

Following \cref{Quasicat subsection}, we construct cocartesian lifts and characterize algebras and $I$-symmetric monoidal functors into $\cC^{I-\sqcup}$ and $\cC^{I-\times}$ in \cref{Comonoids subsection,O-monoids subsection}.
We spell out a corollary in \cref{sec:morita} relating $L_{\Op_I}$-equivalences to the Morita theory of algebraic patterns.
\subsection{Quasicategories modeling \texorpdfstring{$\cC^{I-\sqcup}$ and $\cC^{I-\times}$}{C cocart and C cart}}\label{Quasicat subsection}
\subsubsection{The cocartesian case}
Let $\cT^{\op}$ be a quasicategory and $\Tot \cC \in \sSet^{\cocart}_{/\cT}$ a cocartesian fibration to $\cT$. 
There exists a simplicial set $\Tot \Tot_{\cT} \cC^{I-\sqcup}$ satisfying the universal property
\begin{equation}\label{sseq universal prop}
  \Hom_{/\Tot \uFF_{I,*}}(K,\Tot \Tot_{\cT} \cC^{I-\sqcup}) \simeq \Hom_{/\cT^{\op}}(K \times_{\Tot \uFF_{I,*}} \Tot \uGamma_{I}^*, \Tot \cC), \hspace{40pt} \forall K \in \sSet_{/\Tot \uFF_{I,*}}^{\cocart}.
\end{equation}
\vspace{-\baselineskip}
\begin{lemma}
  The map $\Tot \Tot_{\cT} \cC^{I-\sqcup} \rightarrow \Tot \uFF_{I,*}$ is an inner fibration;
  hence $\Tot \cC^{I-\sqcup}$ is a quasicategory.
\end{lemma}
\begin{proof}
  The proof is exactly analogous to \cite[Prop~2.4.3.3]{HA}:
  apply the universal property
\[\begin{tikzcd}[ampersand replacement=\&, column sep=large]
  {{\Lambda_i^n}} \&\& {\Tot \Tot_{\cT} \cC^{I-\sqcup}} \&\& {\Lambda_i^n \times_{\uFF_{I,*}} \uGamma^*_I} \& {\dcoprod_{\stackrel{U \in \Orb(S)}{f(U) \in S_{n,+}^{\circ}}} \Lambda_i^n} \& \Tot \cC \\
	{{\Delta^n}} \&\& {\Tot \uFF_{I,*}} \&\& {\Delta^n \times_{\uFF_{I,*}} \uGamma^*_I} \& {\dcoprod_{\stackrel{U \in \Orb(S)}{f(U) \in S_{n,+}^{\circ}}} \Delta^n} \& {\cT^{\op}}
	\arrow["{{{f_0}}}", from=1-1, to=1-3]
	\arrow[from=1-1, to=2-1]
	\arrow[""{name=0, anchor=center, inner sep=0}, from=1-3, to=2-3]
	\arrow["\simeq"{description}, draw=none, from=1-5, to=1-6]
	\arrow[""{name=1, anchor=center, inner sep=0}, from=1-5, to=2-5]
	\arrow[from=1-6, to=1-7]
	\arrow[from=1-6, to=2-6]
	\arrow[from=1-7, to=2-7]
	\arrow[dashed, from=2-1, to=1-3]
	\arrow["{\prn{S_{0,+} \rightarrow \cdots \rightarrow S_{n,+}}}"', from=2-1, to=2-3]
	\arrow["\simeq"{description}, draw=none, from=2-5, to=2-6]
	\arrow[dashed, from=2-6, to=1-7]
	\arrow[from=2-6, to=2-7]
  \arrow[shorten <=30pt, shorten >=30pt, squiggly={pre length=30pt, post length=30pt}, tail reversed, from=0, to=1]
\end{tikzcd}\]
  to note that inner horn lifts of $\Tot \cC^{I-\sqcup} \rightarrow \Tot \uFF_{I,*}$ correspond with tuples of inner horn lifts along $\Tot \cC \rightarrow \cT^{\op}$, which exist by assumption that it is a cocartesian fibration (hence an inner fibration).
  The remaining claim follows by noting that $\Tot \uFF_{I,*}$ is a quasicategory, so the composite map $\Tot \Tot_{\cT} \cC^{I-\sqcup} \rightarrow \Tot \uFF_{I,*} \rightarrow *$ is an inner fibration. 
\end{proof}

\begin{proof}[Proof of \cref{UFI exists prop}]
  We've verified that $\Tot \Tot_{\cT} \cC^{I -\sqcup}$ is a quasicategory over $\Tot \uFF_{I,*}$.
  Fixing some quasicategory $\cD$ over $\uFF_{I,*}$ and applying \cref{sseq universal prop} for $K \deq \cD \times \Delta^n$, we find that $\Fun(\cD,\Tot \Tot_{\cT} \cC^{I-\sqcup}) \simeq \Fun_{/\cT^{\op}}(\cD \times_{\Tot \uFF_{I,*}} \Tot \uGamma_I^*, \Tot \cC)$.
  The result then follows by replacing ``quasicategory'' with ``$\infty$-category.''
\end{proof}

\subsubsection{The cartesian case} 
Now, we define $\Tot \Tot_{\cT} \widetilde \cC^{I-\times} \in \sSet_{/\Tot \uFF_{I,*}}$ by the universal property
\begin{equation}\label{cartesian universal prop}
  \Hom_{/\Tot \uFF_{I,*}}(K,\Tot \Tot_{\cT} \widetilde \cC^{I-\times}) \simeq \Hom_{/\cT^{\op}}(K \times_{\Tot \uFF_{I,*}} \Tot \uGamma_{I}^\times, \Tot \cC), \hspace{40pt} \forall K \in \sSet_{/\Tot \uFF_{I,*}}^{\cocart}.
\end{equation}\vspace{-\baselineskip}
\begin{recollection}[{\cite[Def~2.1.2]{Nardin}}]
  A morphism $f$ in $\Tot \uFF_{I,*}$ from $S_+ \in \FF_{I,*,U}$ to $T_+ \in \FF_{I,*,V}$ may be modelled as a morphism of spans
  \[\begin{tikzcd}[sep = small, row sep = tiny]
	S && {f^{-1}(T)} && T \\
	& {\Res_U^V S} \\
	U && V && V
	\arrow[from=1-1, to=3-1]
	\arrow[from=1-3, to=1-1]
  \arrow[from=1-3, to=1-5, "f^{\circ}"]
  \arrow[dashed, hook, from=1-3, to=2-2, "\iota_f"]
	\arrow[from=1-3, to=3-3]
	\arrow[from=1-5, to=3-5]
	\arrow[from=2-2, to=1-1]
	\arrow[from=2-2, to=3-3]
	\arrow[from=3-3, to=3-1]
	\arrow[Rightarrow, no head, from=3-3, to=3-5]
\end{tikzcd}\]
  such that $f^{\circ} \in I$ (c.f. \cref{FI construction}).
  Such a morphism is $\pi_{\uFF_{I,*}}$-cocartesian if $f^{\circ}$ and $\iota_f$ are both equivalences, i.e. it witnesses an equivalence $\Res_U^V S_+ \xrightarrow\sim T_+$.
\end{recollection}
Let $f\colon T_+ \rightarrow S_+$ be a map in $\Tot \uFF_{I,*}$ lying over an orbit map $U \rightarrow V$ and let $\overline{S} \subset S$ be an element of $\uGamma^{\times}_I$ lying over $S_+$.
We would like to construct a Cartesian edge landing on $\overline{S} \subset S$;
we do so by setting $\overline{T} \deq f^{-1}(\Res_U^V \overline{S}) \subset f^{-1}(\Res_U^V S) \subset T$, and letting the associated map $t\colon \prn{f^{-1}(\Res_U^V \overline{S}) \subset T} \rightarrow \prn{\overline{S} \subset S}$ be the canonical one.
The following lemma then follows by unwinding definitions, where $U\colon \uGamma_{I}^\times \rightarrow \uFF_{I,*}$ denotes the forgetful functor.
\begin{lemma}
  $t$ is a $U$-cartesian arrow;
  in particular, $U$ is a cartesian fibration. 
\end{lemma}

The following lemma then follows from \cite[Cor~3.2.2.12]{HTT}.
\begin{lemma}\label{Cart cocart morphisms}
  Let $\widetilde p\colon \widetilde \cC^{I-\times} \rightarrow \Tot \uFF_{I,*}$ be the projection and let $\tilde \alpha\colon F \rightarrow G$ be a $\widetilde \cC^{I-\times}$-morphism lying over a $\Tot \uFF_{I,*}$-morphism $\alpha:T_+ \rightarrow S_+$ lying over an orbit map $U \rightarrow V$. 
  Then, $\tilde \alpha$ is $\widetilde p$-cocartesian if and only if, for all $T' \subset T$, the induced map 
  \[
    F(\alpha^{-1}(\Res_U^V T')) \rightarrow \Res_U^V G(T')
  \]
  is an equivalence;
  in particular, $\widetilde p$ is a cocartesian fibration, so $\widetilde \cC^{I-\times}$ is a quasicategory.
\end{lemma}

\begin{proof}[Proof of \cref{UFI cart exists prop}]
  We concluded in \cref{Cart cocart morphisms} that $\Tot \Tot_{\cT} \widetilde \cC^{I-\times}$ is a quasicategory satisfying \cref{cartesian universal prop}, so it models an $\infty$-category satisfying our universal property by the same argument as \cref{UFI exists prop}. 
\end{proof}

\subsection{\texorpdfstring{$\cO$}{O}-comonoids and (co)cartesian rigidity}\label{Comonoids subsection}
An object of $\Tot \Tot_{\cT} \cC^{I-\sqcup}$ may be viewed as $S_+$ a pointed $V$-set and $\bC = (C_W) \in \cC_S$ an $S$-tuple of elements of $\cC$;
a morphism $f\colon \bC \rightarrow \bD$ may be viewed as a $\Tot \uFF_{I,*}$-map $(S_+ \rightarrow V_{S,+}) \xrightarrow{f} (T_+ \rightarrow V_{T,+})$ together with a collection of maps
\[
  \cbr{f_{W}\colon \Ind_W^U C_W \rightarrow D_U \mid W \in f^{-1}(U)} 
\]
for all $U \in \mathrm{Orb}(T)$.
Unwinding the universal property for cocartesian arrows, we find the following.\footnote{It is here that we use almost-unitality for the cocartesian setting;
if $I$ was not almost essentially unital, then there would exist some $S$ whose $I$-admissible orbits do not together cover $S$, so $\cC^{I-\sqcup} \rightarrow \Tot \uFF_{I,*}$ would not be an inert-cocartesian fibration.}
\begin{proposition}\label{Cocart cocart prop}
  \label{Cocartesian morphisms in cocartesian SMC}
  A map $f\colon \bC \rightarrow \bD$ is $\pi$-cocartesian if and only if $\cbr{f_{W}}$ witness $D_U$ as the indexed coproduct
  \[
    \dcoprod_{W}^{f^{-1}(U)} C_W \xrightarrow{\;\;\; \sim \;\;\;}  D_U
  \]
  for all $U \in \mathrm{Orb}(T)$.
  In particular, $f$ is inert if and only if the following conditions are satisfied:
  \begin{enumerate}[label={(\alph*)}]
    \item The projected morphism $\pi(f):S \rightarrow T$ is inert.
    \item The associated map $C_{f^{-1}(U)} \rightarrow D_{U}$ is an equivalence for all $U \in \mathrm{Orb}(T)$.
  \end{enumerate}
  Hence $\pi\colon \Tot \Tot_{\cT} \cC^{I-\sqcup} \rightarrow \Tot \uFF_{I,*}$ presents a $\Gamma$-$I$-preoperad $\cC^{I-\sqcup}$.
\end{proposition}
\begin{corollary}
  $\cC^{I-\sqcup}$ is an $I$-operad which is an $I$-symmetric monoidal $\infty$-category if and only if $\cC$ admits $I$-indexed coproducts.
\end{corollary}
\begin{proof}
  It follows from \cref{Cocart cocart prop} that $\Tot \Tot_{\cT} \cC^{I-\sqcup} \rightarrow \Tot \uFF_{I,*}$ is a cocartesian fibration if and only if $\cC$ admits $I$-indexed coproducts, so it suffices to verify the following conditions:
  \begin{enumerate}[label={(\alph*)}]\setcounter{enumi}{1}
    \item cocartesian transport yields an equivalence
      \[
        \cC_{S} \simeq \prod_{U \in \Orb(S)} \cC_{U};
      \]  
    \item cocartesian transport yields an equivalence
      \[
        \Map^{T \rightarrow S}_{\Tot \cO^{\otimes}}(\bC,\bD) \simeq \prod_{U \in \Orb(S)} \Map^{T_U \rightarrow U}_{\Tot \cO^{\otimes}}\prn{\bC_{\uU}, D}.
      \]
  \end{enumerate}
  In fact, each condition follows from \cref{Cocart cocart prop}.
\end{proof}

\begin{observation}
  It follows from the above discussion that $U(\cC) = \cC$.
  Moreover, it follows from \cref{Cocart cocart prop} that the indexed tensor product functor $\otimes^S\colon \cC_S \rightarrow \cC_V$ for $\cC^{I-\sqcup}$ is left adjoint to $\Delta^S$, i.e. indexed tensor products in $\cC^{I-\sqcup}$ are indexed coproducts.
\end{observation}

Given $\cO^{\otimes}$ a unital $I$-operad, define a diagram of Cartesian squares in $\Cat_{\cT}$.
\[
  \begin{tikzcd}
    \cO \arrow[d] \arrow[r,hook',"\iota"] \arrow[rd,phantom,"\lrcorner" very near start]
    & \cO^{\otimes}_\Gamma \arrow[d] \arrow[r] \arrow[rd,phantom,"\lrcorner" very near start]
    & \cO^{\otimes} \arrow[d]\\
    *_{\cT} \arrow[r,hook]
    & \uGamma_I^* \arrow[r]
    & \uFF_{I,*}
  \end{tikzcd}
\]
Note that the objects of $\cO_{\Gamma,V}^{\otimes}$ consist of triples $(S_+,U,X)$ where $U \in \mathrm{Orb}(S)$ and $X \in \cO_S$, and the image of $\iota$ is equivalent to the triples where $S \in \cT_{/V}$, hence $U = S$.

Further note that cocartesian transport along the inert morphism $U_+ \hookrightarrow S_+$ induces an equivalence
\[
  \Map_{\cO^{\otimes}_{\Gamma,V}}(\iota Y,(S_+,U,X)))
  \simeq \Map_{\cO^{\otimes}_{\Gamma,V}}(\iota Y,(U_+,U,X_U)))
\]
for all $Y \in \cO$.\footnote{This utilizes unitality of $\cO^{\otimes}$, as we implicitly use that, for each orbit $U' \in \Orb(S)$ other than $U$, the space $\cO(\emptyset_{U'};X_{U'})$ is contractible.}
In particular, $\iota$ witnesses $\cO$ as a \emph{colocalizing} $\cT$-subcategory, with colocalization $\cT$-functor
\[
  R(S_+,U,X) \simeq (U_+,U,X_U).
\]
This interacts with Kan extensions via the following lemmas.
\begin{lemma}\label{Adjoint Kan extension}
  Suppose $F\colon \cC \rightarrow \cD$ is a $\cT$-functor and $L\colon \cC \rightarrow \cE$ is $\cT$-left adjoint to $R\colon \cE \rightarrow \cC$.
  Then, $FR$ is the $\cT$-left Kan extension of $F$ along $L$.
\end{lemma}
\begin{proof}
  Using \cite[Thm~10.5]{Shah2} we simply repeat the nonequivariant proof:
  Yoneda's lemma yields
  \begin{align*}
    \Nat_{\cT}(L_!F, G)
    &\simeq \Nat_{\cT}(F, GL)\\
    &\simeq \Nat_{\cT}(FR, G),
  \end{align*}
  so another application of Yoneda's lemma constructs a natural equivalence $L_! F \sim FR$.
\end{proof}

\begin{lemma}\label{Kan extension lemma}
  Fix a $\cT$-functor $A\colon \cO^{\otimes}_\Gamma \rightarrow \cC$.
  Then, the following are equivalent
  \begin{enumerate}[label={(\alph*)}]
      \item The corresponding map $\cO^{\otimes} \rightarrow \cC^{I-\sqcup}$ is a functor of $I$-operads.
      \item For all morphisms $\alpha$ in $\cO^{\otimes}_\Gamma$ whose image in $\cO^{\otimes}$ is inert, $A(\alpha)$ is an equivalence in $\cC$.
      \item If $f\colon (S_+,U,X) \rightarrow (U_+,U,X_U)$ is a cocartesian lift of the corresponding inert morphism, then $A(f)$ is an equivalence.
      \item $A$ is $\cT$-left Kan extended from $\cO$.
    \end{enumerate}
    Furthermore, every functor $F:\cO \rightarrow \cC$ admits a left Kan extension along $\cO \hookrightarrow \cO^{\otimes}_\Gamma$;
    in particular, the forgetful functor $\uAlg_{\cO}(\cC) \rightarrow \uFun_G(\cO,\cC)$ is an equivalence.
\end{lemma}
\begin{proof}
  (a) $\iff$ (b) follows immediately from \cref{Cocartesian morphisms in cocartesian SMC}.
  (b) $\iff$ (c) is immediate by definition.
  (c) $\iff$ (d) and the remaining statement both follow by \cref{Adjoint Kan extension}.
\end{proof}

\subsection{\texorpdfstring{$\cO$}{O}-monoids}\label{O-monoids subsection}
We start with the following.
\begin{proposition}\label{CItimes is a preoperad}
  $\Tot \Tot_{\cT} \cC^{I-\times} \rightarrow \Tot \uFF_{I,*}$ is a cocartesian fibration, so in particular, it presents a $\Gamma$-$I$-preoperad $\cC^{I-\times}$;
  moreover, its straightening is an $I$-symmetric monoidal $\infty$-category if and only if $\cC$ admits $I$-indexed products, in which case the indexed tensor functors in $\cC$ are indexed products.
\end{proposition}
\begin{proof}
  For the first statement, it suffices to observe that the cocartesian arrows described in \cref{Cart cocart morphisms} lie in $\Tot \Tot_{\cT} \cC^{I-\times}$.
  For the second, note by unwinding definitions that cocartesian transport induces a fully faithful functor
  \[
    \cC_S \rightarrow \prod_{U \in \Orb(S)} \cC_U
  \]
  Moreover, this is essentially surjective if and only if $\cC$ admits $I$-indexed products, as desired.
\end{proof}

We organize ourselves around the following observations.
\begin{observation}
  The projection $\Tot \Tot_{\cT}\cO^{\otimes} \times_{\Tot\uFF_{I,*}} \Tot \uGamma_{I}^{\times} \rightarrow \Tot \Tot_{\cT}\cO^{\otimes}$ admits a left adjoint $L$ sending $X \in \cO^{\otimes}_{S_+}$ to $(X,S \subset S)$;
  the unit map of this adjunction is evidently an equivalence, so $L\colon \Tot\Tot_{\cT} \cO^{\otimes} \rightarrow \Tot \Tot_{\cT}\cO^{\otimes} \times_{\uFF_{I,*}} \uGamma_{I}^{\times}$ is fully faithful.
\end{observation}
\begin{observation}
  The section $L\colon \Tot \uFF_{I,*} \rightarrow \Tot \uGamma_{I}^\times$ induces a natural transformation $K \simeq K \times_{\Tot \uFF_{I,*}} \Tot \uFF_{I,*} \rightarrow K \times_{\Tot \uFF_{I,*}} \Tot \uGamma_I^\times$, which induces a natural transformation $\Tot \Tot_{\cT} \widetilde \cC^{I-\times} \rightarrow \Tot \cC$ under Yoneda's lemma.
  Unwinding definitions using \cref{Cart cocart morphisms}, this presents a $\cT$-functor $\Tot_{\cT} \widetilde \cC^{I-\times} \rightarrow \cC$.
\end{observation}

Given a $\cT$-functor $\Tot_{\cT} \cO^{\otimes} \xrightarrow{\;\;\;\; \varphi\;\;\;\;} \Tot_{\cT}\widetilde \cC^{I-\times}$, we acquire a corresponding functor
\[
  \Tot \Tot_{\cT} \cO^{\otimes} \xrightarrow{\;\;\;\;L\;\;\;\;} \Tot_{\cT} \cO^{\otimes} \times_{\Tot \uFF_{I,*}} \Tot \uGamma_I^\times \xrightarrow{\;\;\;\; \varphi'\;\;\;\;} \Tot \cC
\]
over $\cT^{\op}$.
Now, the following observation is important.
\begin{observation}
  The description of cocartesian arrows of \cref{Cart cocart morphisms,CItimes is a preoperad} implies that $\varphi'$ and $L$ are unstraightened from $\cT$-functors.
\end{observation}

Now, given a $\Gamma$-$I$-preoperad $\cO^{\otimes}$, we say that an \emph{$\cO$-monoid in $\cC$} is a $\cT$-functor $\Tot_{\cT} \cO^{\otimes} \rightarrow \cC$ satisfying the condition that, for all $X \in \cC_S$, the maps $\Res_U^V F(X) \rightarrow F(X_U)$ induced by cocartesian transport witness $F(X)$ as the indexed product
\[
  F(X) \simeq \dprod_U^S F(X_U).
\]\vspace{-\baselineskip}
\begin{proposition}\label{Cartesian algebras generalized}
  Fix $\cC$ a $\cT$-$\infty$-category and $\cO^{\otimes}$ a $\Gamma$-$I$-preoperad.
  Then, the postcomposition functor $\Alg_{\cO}(\cC^{I-\times}) \rightarrow \Fun_{\cT}(\Tot_{\cT}\cO^{\otimes}, \cC)$ is fully faithful with image spanned by the $\cO$-monoids. 
\end{proposition}

\begin{lemma}\label{O-monoid pre-lemma}
  The following conditions are equivalent:
  \begin{enumerate}[label={(\alph*)}]
    \item $\varphi$ is a map of $\Gamma$-$I$-preoperads.
    \item For all morphisms $\alpha$ in $\Tot \Tot_{\cT} \cO^{\otimes} \times_{\Tot \uFF_{I,*}} \Tot \uGamma_{I}^\times$ whose image in $\cO^{\otimes}$ is inert $\varphi'(\alpha)$ is an equivalence..
    \item If $f\colon (\overline{S}_+ \rightarrow V_+, \overline{S} , F, X) \rightarrow (S_+ \rightarrow V_+,\overline{S},F,X)$ is a cocartesian lift of the corresponding inert morphism, then $\varphi(f)$ is an equivalence.
    \item The composite $\Tot_{\cT} \cO^{\otimes} \xrightarrow{\;\;\; \varphi\;\;\;} \Tot_{\cT} \widetilde \cC^{I-\times} \rightarrow \cC$ is homotopic to $\varphi'$.
    \item The composite $\Tot_{\cT} \cO^{\otimes} \xrightarrow{\;\;\; \varphi\;\;\;} \Tot_{\cT} \widetilde \cC^{I-\times} \rightarrow \cC$ is $\cT$-right Kan extended from $\varphi'$ along $L$. 
  \end{enumerate}
\end{lemma}
\begin{proof}
  \cref{Cart cocart morphisms} immediately implies that (a) $\iff$ (b) $\iff$ (c) $\iff$ (d).
  (d) $\iff$ (e) follows  from \cref{Adjoint Kan extension}.
\end{proof}

The following lemma follows by unwinding definitions.
\begin{lemma}\label{O-monoid lemma}
  Suppose $\varphi$ is a functor of $\Gamma$-$I$-preoperads.
  Then, the following conditions are equivalent:
  \begin{enumerate}[label={(\alph*)}]
    \item $\varphi$ factors through a $\Gamma$-$I$-preoperad map $\overline{\varphi}\colon \cO^{\otimes} \rightarrow \cC^{I-\times}$.
    \item $\varphi'$ is an $\cO$-monoid.
  \end{enumerate}
\end{lemma}

\begin{proof}[Proof of \cref{Cartesian algebras generalized}]
  Since $\Tot_{\cT}\cC^{I-\times} \hookrightarrow \Tot_{\cT}\tilde \cC^{I-\times}$ is fully faithful,
  the first of the following functors is fully faithful
  \[
    \Alg_{\cO}\prn{\cC^{I-\times}} \hookrightarrow \Alg_{\cO}\prn{\widetilde \cC^{I-\times}} \hookrightarrow \Fun_{\cT}\prn{\Tot_{\cT}\cO^{\otimes}, \cC}.
  \]
  By \cref{O-monoid pre-lemma}, the second is fully faithful, and by \cref{O-monoid pre-lemma} the image of the composite functor consists of the $\cO$-mononids.
\end{proof}

We may additionally characterize $I$-symmetric monoidal functors via a lift of \cite[Prop~2.4.1.7]{HA}, which also follows immediately from \cref{Cart cocart morphisms}.
\begin{lemma}\label{SM cart lemma}
  Suppose $\cC$ has $I$-indexed products, $\Tot_{\cT} \cO^{\otimes} \rightarrow \uFF_{I,*}$ is a cocartesian fibration, and \cref{O-monoid lemma} is satisfied.
  Then, the following conditions are equivalent.
  \begin{enumerate}[label={(\alph*)}]
    \item $\overline{\varphi}$ is a $\cT$-functor of cocartesian fibrations over $\uFF_{I,*}$.
    \item If $f\colon \bD \rightarrow \bD'$ is an active arrow in $\Tot \Tot_{\cT} \cO^{\otimes}$, then map $\varphi'(f)$ is an equivalence.
  \end{enumerate}
\end{lemma}

Now, we will also lift \cite[Prop~2.4.1.6]{HA}.
We have a fully faithful $\cT$-functor $\iota\colon \cO \hookrightarrow \Tot_{\cT} \cO^{\otimes}$.
\begin{lemma}\label{RKE removing domain lemma}
  Suppose $\varphi$ satisfies the conditions of \cref{SM cart lemma} and the action maps $\otimes^S \colon \cO_S \rightarrow \cO_V$ are right adjoint to the restriction maps $\Delta^S\colon  \cO_V \rightarrow \cO_S$.
  Then, the functor $\widetilde \varphi\colon \Tot_{\cT} \cO^{\otimes} \rightarrow \cC$ is right Kan-extended from the $I$-product-preserving functor $\cO \rightarrow \Tot_{\cT} \cO^{\otimes} \rightarrow \cC$ along $\iota$.
\end{lemma}
\begin{proof}
  The assumptions imply that there is a right $\cT$-adjoint $R\colon \Tot_{\cT} \cO^{\otimes} \rightarrow \cO$ to $\iota$, sending $(X_U) \mapsto \bigotimes_U^S X_U \simeq \prod_U^S X_U$.
  The $\cO$-monoid assumption shows that $\widetilde \varphi \sim \varphi \circ \iota \circ R$, which shows that $\widetilde \varphi$ is right Kan extended from $\varphi \circ \iota$ along $\iota$;
  moreover, the $\cO$-monoid assumption shows that $\varphi \circ \iota$ is $I$-product-preserving.
\end{proof}

\subsection{(Co)cartesian rigidity}
\begin{proposition}\label{Cocartesian rigidity pre-prop}
  Suppose $\cO^{\otimes}$ is an $I$-symmetric monoidal $\infty$-category such that $\otimes^S\colon \cO_S \rightarrow \cO_V$ is right adjoint to $\Delta^S\colon \cO_V \rightarrow \cO_S$ for all $S \in \FF_{I,V}$.
  Then, the forgetful functor
  \[
    U\colon \Fun_I^{\otimes}\prn{\cO^{\otimes}, \cC^{I-\times}} \rightarrow \Fun_{\cT}(\cO, \cC)
  \]
  is fully faithful with essential image spanned by the $I$-product preserving $\cT$-functors.
\end{proposition}
\begin{proof}
  Let $\widetilde \uFun_{\cT}(\Tot_{\cT} \cO^{\otimes}, \cC) \subset \uFun_{\cT}(\Tot_{\cT} \cO^{\otimes}, \cC)$ be the equivalent image of $\Fun_{I}^{\otimes}(\cO^{\otimes}, \cC^{\otimes})$.
  \cref{RKE removing domain lemma} constructs the solid portion of a diagram
\[\begin{tikzcd}[ampersand replacement=\&]
	{\Fun_{\cT}^{I-\times}(\cO,\cC)} \& {\widetilde\Fun_{\cT}\prn{\Tot_{\cT} \cO^{\otimes}, \cC}} \& {\Fun_{\cT}\prn{\Tot_{\cT} \cO^{\otimes}, \cC}} \\
	\& {\Fun_{\cT}^{I-\times}(\cO,\cC)} \& {\Fun_{\cT}\prn{\cO, \cC}}
	\arrow[dashed, from=1-1, to=1-2, hook']
	\arrow["{\iota_*}", curve={height=-18pt}, hook', from=1-1, to=1-3]
	\arrow[equals, from=1-1, to=2-2]
	\arrow[hook', from=1-2, to=1-3]
	\arrow["{\iota^*}", from=1-2, to=2-2]
	\arrow["U"{description}, from=1-2, to=2-3]
	\arrow["{\iota^*}", from=1-3, to=2-3]
	\arrow[hook, from=2-2, to=2-3]
\end{tikzcd}\]
  It suffices to verify that the dashed arrow exists, i.e. right Kan extensions of $I$-product-preserving functors along $\iota$ satisfy the conditions of \cref{SM cart lemma};
  but this follows by unwinding definitions.
\end{proof}

\begin{corollary}\label{many vops corollary}
  If $\cC$ has $I$-indexed products, then there exists a unique $I$-symmetric monoidal equivalence 
  $\cC^{I-\times} \simeq \prn{\prn{\cC^{\vop}}^{I-\sqcup}}^{\vop}$ lying over the equivalence 
  $\cC \simeq \prn{\cC^{\vop}}^{\vop}$;
  if $\cC$ has $I$-indexed coproducts, then there exists a unique $I$-symmetric monoidal equivalence 
  $\cC^{I-\sqcup} \simeq \prn{\prn{\cC^{\vop}}^{I-\times}}^{\vop}$ lying over the equivalence
  $\cC \simeq \prn{\cC^{\vop}}^{\vop}$.
\end{corollary}
\begin{proof}
  By conservativity of the underlying category (see \cref{Underlying conservative proposition}), it suffices to construct a unique $I$-symmetric monoidal \emph{functor} lying over the identity in each case. 
  For the first case, by \cref{Cocartesian rigidity pre-prop} it suffices to note that $\prn{\prn{\cC^{\vop}}^{I-\sqcup}}^{\vop}$ has $\Delta^S \dashv \otimes^S$.
  The second case follows from the first under the following equivalence of arrows, where 
  $\cD \deq \cC^{\vop}$.
  \[\begin{tikzcd}[ampersand replacement=\&, column sep=tiny]
	{\Fun_I^{\otimes}\prn{\prn{\cD^{\vop}}^{I-\sqcup}, \prn{\cD^{I-\times}}^{\vop}}} \& {\Fun_I^{\otimes}\prn{\prn{\prn{\cD^{\vop}}^{I-\sqcup}}^{\vop}, \prn{\prn{\cD^{I-\times}}^{\vop}}^{\vop}}} \& {\Fun_I^{\otimes}\prn{\prn{\prn{\cD^{\vop}}^{I-\sqcup}}^{\vop},\cD^{I-\times}}} \\
	{\Fun_{\cT}\prn{\cD^{\vop}, \cD^{\vop}}} \& {\Fun_{\cT}\prn{\prn{\cD^{\vop}}^{\vop}, \prn{\cD^{\vop}}^{\vop}}} \& {\Fun_{\cT}\prn{\cD,\prn{\cD^{\vop}}^{\vop}}}
	\arrow["\simeq"{description}, draw=none, from=1-1, to=1-2]
	\arrow[from=1-1, to=2-1]
	\arrow["\simeq"{description}, draw=none, from=1-2, to=1-3]
	\arrow[from=1-2, to=2-2]
	\arrow[from=1-3, to=2-3]
	\arrow["\simeq"{description}, draw=none, from=2-1, to=2-2]
	\arrow["\simeq"{description}, draw=none, from=2-2, to=2-3]
\end{tikzcd}\]
\end{proof}

\begin{corollary}\label{Cocartesian rigidity dual pre-prop}
  Suppose $\cO^{\otimes}$ is an $I$-symmetric monoidal $\infty$-category such that $\otimes^S\colon \cO_S \rightarrow \cO_V$ is left to $\Delta^S\colon \cO_V \rightarrow \cO_S$ for all $S \in \FF_{I,V}$.
  Then, the forgetful functor
  \[
    U\colon \Fun_I^{\otimes}\prn{\cO^{\otimes}, \cC^{I-\sqcup}} \rightarrow \Fun_{\cT}(\cO, \cC)
  \]
  is fully faithful with essential image spanned by the $I$-coproduct preserving $\cT$-functors.
\end{corollary}
\begin{proof}
  This follows by taking vertical opposites of \cref{Cocartesian rigidity pre-prop} in light of \cref{many vops corollary}. 
\end{proof}

We are now ready to prove our main generalization of \cref{Cocartesian rigidity main theorem generalized} (see p. \pageref{Cocartesian rigidity main theorem generalized}).
\begin{proof}[Proof of \cref{Cocartesian rigidity main theorem generalized}]
  We begin with the cartesian cases.
  To see that $(-)^{I-\times}$ is fully faithful, it suffices to combine \cref{CItimes is a preoperad,Cocartesian rigidity pre-prop}.
  The compatibility with $U$ is obvious, and the description of the image follows immediately from \cref{Cocartesian rigidity pre-prop}.
  The cocartesian case follows by the same argument using \cref{Cocartesian rigidity dual pre-prop}.
\end{proof}

\subsection{Wirtm\"uller maps}
Suppose $\cO^{\otimes}$ is an $I$-symmetric monoidal $\infty$-category and $\cC$ has $I$-indexed coproducts.
The equivalence $\Alg_{\cO}(\cC^{I-\sqcup}) \simeq \Fun_{\cT}(\cO, \cC)$ embeds $I$-symmetric monoidal functors $\cO^{\otimes} \rightarrow \cC^{I-\sqcup}$ as a full subcategory of $\cT$-functors $\cO \rightarrow \cC$.
We now record the property identifying this full subcategory.

\begin{observation}\label{Appendix wirthmuller}
  Let $F\colon \cD^{\otimes} \rightarrow \cE^{\otimes}$ be a lax $I$-symmetric monoidal functor.
  Then, the universal property for cocartesian arrows constructs, for each active arrow $\Ind_V^{\cT} S \rightarrow V$, an arrow
  \[
    \mu_S\colon \bigotimes_U^S F(-) \implies F\prn{\bigotimes_U^S -}
  \]
  such that $F$ is $I$-symmetric monoidal if and only if $\mu_S$ is an equivalence for all $S \in \uFF_I$.
  In particular, in the case of the lax $I$-symmetric monoidal functor $F\colon \cO^{\otimes} \rightarrow \cC^{I-\sqcup}$ classified by a functor $G\colon \cO \rightarrow \cC$, the arrow $\mu$ has the type
  \[
    \coprod_U^S G(-) \implies G\prn{\bigotimes_U^S -};
  \]
  moreover, unwinding definitions, in the case that $\cO^{\otimes} = \cC^{\otimes}$ is an $I$-symmetric monoidal structure on $\cC$ and, $G$ is the identity, and $1_\bullet \in \Gamma^{\upsilon(I)}\cO$ is initial, this map is precisely the $\otimes$-Wirthm\"uller map constructed in \cref{Wirthmuller construction}.
  In particular, we've observed that the identity classifies an $I$-symmetric monoidal equivalence $\cC^{\otimes} \xrightarrow{\;\;\;\; \sim \;\;\;\;} \cC^{I-\coprod}$ if and only if $\cC^{\otimes}$ has $I$-admissible $\otimes$-Wirthm\"uller isomorphisms.
\end{observation}

\subsection{A technical corollary on \texorpdfstring{$n$}{n}-Morita equivalences}\label{sec:morita}
A Segal morphism of algebraic patterns $\varphi\colon \fO \rightarrow \fP$ is called an \emph{$n$-Morita equivalence} if, for all complete $(n+1)$-categories $\cC$, the induced functor
\[
  f^*\colon \Seg_{\fP}(\cC) \rightarrow \Seg_{\fO}(\cC)
\]
is an equivalence;
in fact, it suffices to check this in the case $\cC  = \cS_{\leq n}$ \cite[Prop~2.1.9]{Barkan-arity}.
We have the following corollary.
\begin{corollary}\label{Morita equivalence vs L equivalence}
  Suppose $\varphi\colon \cP^{\otimes} \rightarrow \cO^{\otimes}$ is a morphism of $\Gamma$-$I$-preoperads such that the induced $\cT$-functor $U\cP \rightarrow U\cO$ is essentially surjective.
  Then, $\varphi$ is an $n$-Morita equivalence if and only if the associated map of $\cT$-operads $h_{n+1}L_{\Op_{\cT}} \cP^{\otimes} \rightarrow h_{n+1}L_{\Op_{\cT}} \cO^{\otimes}$ is an equivalence.
\end{corollary}
\begin{proof}
  There is an equivalence $\Seg_{\cO}(\cC) \simeq \Alg_{\cO}(\cC^{I-\times}) \simeq \Alg_{L_{\Op_{\cT}} \cO}(\cC^{I-\times})$
  natural in Segal morphisms over $\Span(\FF_{\cT})$, so the result follows from the recognition result for $n$-equivalences of $I$-operads \cite{EBV}. 
\end{proof}

Now, we define the \emph{$I$-preoperads} $\PreOp_I \deq \Cat_{\cT, /\Span_I(\FF_{\cT})}^{\Int-\cocart}$.
In \cite[\S~A.1]{EBV} we proved that, under the assumption that $\cO^{\otimes}$ is an $I$-operad, the canonical strong Segal morphism $s^* \cO^{\otimes} \rightarrow \cO^{\otimes}$ is a Morita equivalence and $s^* \cO^{\otimes}$ is soundly extendable;
in fact, we only used that $\cO^{\otimes}$ is an $I$-operad to conclude that $s^* \cO^{\otimes}$ is soundly extendable, and the same proof shows that $s^* \cO^{\otimes} \rightarrow \cO^{\otimes}$ underlies a natural Morita equivalence when $\cO^{\otimes}$ is an $I$-preoperad, and in particular, a natural $n$-Morita equivalence.

Moreover, this underlies the counit map of $I$-preoperads $\varepsilon\colon  s_! s^* \cO^{\otimes} \rightarrow \cO^{\otimes}$, which is an $L_{\Op_G}$-equivalence.
In particular, $\varepsilon$ becomes an equivalence after $\Alg_{(-)}(\ucS_{G, \leq n+1})$.
This yields a chain of natural equivalences.
\[
  \Alg_{\cO}(\ucS_{G, \leq n+1}) \simeq \Alg_{s_!s^*\cO}(\ucS_{G, \leq n+1}) \simeq \Alg_{s^* \cO}(\ucS_{G, \leq n+1}) \simeq \Seg_{s^* \cO}(\cS_{\leq n+1}) \simeq \Seg_{\cO}(\cS_{\leq n + 1})
\]
In particular, the same proof as \cref{Morita equivalence vs L equivalence} yields the following.
\begin{corollary}
  Suppose $\varphi\colon \cP^{\otimes} \rightarrow \cO^{\otimes}$ is a morphism of $I$-preoperads such that the induced $\cT$-functor $U\cP \rightarrow U\cO$ is essentially surjective.
  Then, $\varphi$ is an $n$-Morita equivalence if and only if the associated map of $\cT$-operads $h_{n+1}L_{\Op_{\cT}} \cP^{\otimes} \rightarrow h_{n+1}L_{\Op_{\cT}} \cO^{\otimes}$ is an equivalence.
\end{corollary}

\resumetocwriting
\section{\tI-operadic disintegration and assembly}\label{sec:disinte gration}
\stoptocwriting
In this appendix, we assume familiarity with the minutiae of \cite{EBV,Barkan} and 
prove \cref{thm:disintegration and assembly}.

\subsection{The algebraic pattern for families of \tI-operads}%
Fix $\cC$ a $\cT$-$\infty$-category.
\begin{construction}
  The \emph{totally inert pattern structure on $\cC$} has
  \[
    \prn{\Tot^{\Int} \cC}^{\el} \deq
    \prn{\Tot^{\Int} \cC}^{\Int} \deq 
    \Tot \cC; \hspace{50pt} \prn{\Tot^{\Int} \cC}^{\act} \deq \prn{\Tot \cC}^{\simeq}.
  \]
  We define the \emph{$\cC$-family pattern}
  \[
    \Tot^{\Int} \prn{\cC \times \uFF_{I,*}} \deq \Tot \cC^{\Int} \times_{\cT^{\op, \Int}} \Tot \uFF_{I,*}.\qedhere
  \]
\end{construction}

The following observation is as crucial as it is immediate.
\begin{observation}\label{T-space pattern}
  If $\cC$ is a $\cT$-space, then $\Tot^{\Int}(\cC \times \uFF_{I,*}) \rightarrow \Tot \uFF_{I,*}$ is an inert-cocartesian fibration, and the domain has the induced pattern structure.
  In particular, in this case, $\Tot^{\Int} (\cC \times \uFF_{I,*})$ is the pattern underlying the $\Gamma$-$I$-preoperad $\cC \times \uFF_{I,*}$.
\end{observation}

We begin by identifying $\cC$-indexed diagrams of $I$-operads.
\begin{proposition}\label{C-family pattern}
  There exists a natural equivalence
  \[
    \Fbrs\prn{\Tot^{\Int}\prn{\cC \times \uFF_{I,*}}} \simeq \Fun_{\cT}\prn{\cC,\uOp_I},
  \]
  so that when $\cC \simeq *$, this is the usual equivalence $\Fbrs(\uFF_{I,*}) \simeq \uOp_I$.
\end{proposition}
To prove this, we use the equifibered theory, focusing on the following lemmas.
\begin{lemma}\label{Segal family lemma}
  There exists a natural equivalence $\Seg_{\Tot(\cC \times \uFF_{I,*})}(\cD) \simeq \Fun_{\cT} \prn{\cC, \uCMon_I(\cD)}$.
\end{lemma}
\begin{proof}
  First off, we get an embedding
  \begin{align*}
    \Seg_{\Tot(\cC \times \uFF_{I,*})}(\cD) 
    &\subset \Fun(\Tot(\cC \times \uFF_{I,*}), \cD);\\
    &\simeq \Fun_{\cT}(\cC \times \uFF_{I,*}, \uCoFr^{\cT}(\cD));\\
    &\simeq \Fun_{\cT}\prn{\cC, \uFun_{\cT}\prn{\uFF_{I,*}, \uCoFr^{\cT} \cD}}
  \end{align*}
  characterized by the Segal condition that the restricted functor $\Tot \prn{\cC \times \uFF_{I,*}^{\Int}} \rightarrow \cD$ is right Kan extended from $\Tot \prn{\cC \times \uFF_{I,*}^{\el}}$.
  Now, unwinding conditions, this corresponds with the condition that the value $V$-functors $\uFF_{I_V,*}^{\Int} \rightarrow \uCoFr^{V} \cD$ are each right-Kan extended from $\uFF_{I_V,*}^{\Int}$, i.e. the corresponding functor factors through $\uCMon_I(\cD) \subset \uFun_{\cT}\prn{\uFF_{I,*}, \uCoFr^{\cT} \cD}$ (c.f. \cite[Ex~1.17]{Nardin_thesis}).
\end{proof}
\begin{lemma}
  $\Env\Tot(\cC \times \uFF_{I,*})$ corresponds naturally with the constant $\cT$-functor over $\uFF_{I,*}^{I-\sqcup}$;
  a natural transformation $F \rightarrow \Env(\cC \times \uFF_{I,*})$ is equifibered if and only if it is pointwise-equifibered.
\end{lemma}
\begin{proof}
  This follows by explicitly identifying the active arrows in $\Tot(\cC \times \uFF_{I,*})$ as products of equivalences and active arrows of $\Tot \uFF_{I,*}$.
\end{proof}

\begin{proof}[Proof of \cref{C-family pattern}]
  The above work constructs a string of natural equivalences
  \begin{align*}
    \Fbrs(\Tot(\cC \times \uFF_{I,*}))
    &\simeq \Fun_{\cT}\prn{\cC, \uCMon_I(\cD)}^{\mathrm{pointwise-equifibered}}_{/\Delta \uFF_{I,*}}\\
    &\simeq \Fun_{\cT}\prn{\cC, \uCMon_I(\cD)^{\mathrm{equifibered}}_{/\uFF_{I,*}}}\\
    &\simeq \Fun_{\cT} \prn{\cC, \uOp_I}.\qedhere
  \end{align*}
\end{proof}

Now, this is closely related to $(-)^{I-\sqcup}$, as described by the following construction.
\begin{construction}
  Fix $\cC$ a $\cT$-$\infty$-category.
  Then, pullback along the projection $(-) \times_{\Tot \uFF_{I,*}} \Tot \uGamma_I^* \rightarrow (-)$ determines a natural transformation
  \[
    \Fun_{/\Tot \uFF_{I,*}}\prn{-,\Tot \Tot_{\cT} \cC^{I-\sqcup}} \simeq \Fun_{/\cT^{\op}} \prn{- \times_{\Tot \uFF_{I,*}} \Tot \uGamma_I^*,\cC} \leftarrow \Fun_{/\cT^{\op}}\prn{-,\cC} \simeq \Fun_{/\Tot \uFF_{I,*}}\prn{-,\cC \times \Tot \uFF_{I,*}},
  \]
  which corresponds with a functor $\gamma \colon \Tot \cC \times \Tot \uFF_{I,*} \rightarrow \Tot \Tot_{\cT}\cC^{I-\sqcup}$ under Yoneda's lemma.
  Note that $\gamma(C,S) \simeq \Delta^SC$.
  Moreover, this is compatible with $\Gamma$-$I$-preoperadic structure in the case $\cC = \Tot_{\cT} \cP^{\otimes}$:
  \[
    \begin{tikzcd}[ampersand replacement=\&, column sep=large, row sep=small]
	{\uFF_{I,*} \times \Tot_{\cT} \cP^{\otimes}} \&\& {\uFF_{I,*}  \times \uFF_{I,*}} \\
	\& {\prn{\Tot_{\cT} \cP^{\otimes}}^{I-\sqcup}} \& {\prn{\uFF_{I,*}}^{I-\sqcup}} \& {\uFF_{I,*}} \\
	\& {\uFF_{I,*}}
	\arrow["{\id \times \pi_{\cP}}", from=1-1, to=1-3]
	\arrow["\gamma"{description}, from=1-1, to=2-2]
	\arrow["{\pr_1}"', curve={height=6pt}, from=1-1, to=3-2]
	\arrow["\wedge"{description}, from=1-3, to=2-4]
	\arrow["{\pi_{\cP}^{I-\sqcup}}"{description}, from=2-2, to=2-3]
	\arrow[from=2-2, to=3-2, "\pi"]
	\arrow["\vee"{description}, from=2-3, to=2-4]
\end{tikzcd}
  \]
  Pulling back to $\Span_I(\FF_{V})$, we acquire a simpler diagram
  \[
    \begin{tikzcd}[ampersand replacement=\&, column sep=large, row sep=small]
      {\Span_I(\FF_{V}) \times \Tot \cP^{\otimes}} \& {\Span_I(\FF_{\cT}) \times \Span_I(\FF_{\cT})} \\
      \& {\prn{\Tot_{V} \cP^{\otimes}}^{I-\sqcup}} \& {\Span_I(\FF_{\cT})} \\
      \& {\Span_I(\FF_{V})}
      \arrow["{\id \times \pi_{\cP}}", from=1-1, to=1-2]
      \arrow["\gamma"{description}, from=1-1, to=2-2]
      \arrow["{\pr_1}"', curve={height=12pt}, from=1-1, to=3-2]
      \arrow["\wedge"{description}, from=1-2, to=2-3]
      \arrow["\rho"{description}, from=2-2, to=2-3]
      \arrow["\pi", from=2-2, to=3-2]
    \end{tikzcd}
  \]
  so in particular, $\gamma$ becomes a \emph{bifunctor} under the alternative structure map $\rho$.
\end{construction}
In this paper, we mostly care about the case that $\cC$ is a generic $\cT$-space.
We will use the following specialization of Barkan's morita equivalence recognition result (c.f. \cite[Thm~2.3.3.23,Thm~2.3.3.26]{HA}).\footnote{To see this as a specialization of Barkan's result, note that by \cite{EBV}, $\Tot \uFF_{I,*}$ is soundly extendable, so $\Tot_{\cT}$ of an $I$-operad is soundly extendable by \cite[Lem~4.1.15]{Barkan}. The remaining modifications necessary are the observation that $\cO^{\el} \simeq \Tot U\cO$ (so condition (1) implies condition (1) of \cite[Thm~5.1.1]{Barkan}) and the identifications $\Fbrs(\cP) \simeq \Op_{I, /\cP^{\otimes}}$ of \cite{EBV} and \cite[Cor~4.1.17]{Barkan} as well as $\Mon_{\cP}(\cC) \simeq \Seg_{\cP}(\cC)$ of \cref{CMon is CAlg corollary}.}
\begin{proposition}[{\cite[Prop~3.1.16,Thm~5.1.1]{Barkan}}]\label{Equivalence theorem}
  Suppose $f\colon \cO^{\otimes} \rightarrow \cP^{\otimes}$ is a strong Segal morphism of $\Gamma$-$I$-preoperads such that $\cP^{\otimes}$ presents an $I$-operad and the following conditions hold:
  \begin{enumerate}[label={(\alph*)}]
    \item the $\cT$-functor $U\cO \rightarrow U\cP$ is an equivalence, and
    \item for every $O \in U\cO$, the map of spaces $\prn{\cO^{\act}_{/O}}^{\simeq} \rightarrow \prn{\cP^{\act}_{/f(O)}}^{\simeq}$ is an equivalence.
  \end{enumerate}
  Then, the functors $f^*\colon \Mon_{\cP}(\cC) \rightarrow \Mon_{\cO}(\cC)$ and $f^*\colon \Op_{I, /\cP^{\otimes}} \rightarrow \Fbrs(\cO)$ are equivalences.
\end{proposition}

\begin{proposition}\label{gamma is approx}
  When $X$ is a $\cT$-space, the functor $\gamma\colon \Tot^{\Int}\prn{X \times \Tot \uFF_{I,*}} \rightarrow  \Tot \Tot_{\cT} X^{I-\sqcup}$ satisfies the conditions of \cref{Equivalence theorem};
  in particular, $\gamma$ is an $L_{\Op_I}$ localization map.
\end{proposition}
\begin{proof}
  By unwinding definitions we find that $\Tot \gamma$ is an iso-segal morphism, and in particular it is a strong Segal morphism.
  Moreover, condition (a) follows simply by unwinding definitions.
 
  For condition (b), note that the elements of $\prn{\Tot_I X^{I-\sqcup}}^{\act}_{/\gamma(x,S)}$ correspond with maps $f\colon T \rightarrow S$ in $\FF_{\cT}$ together with elements $(y_U)_T \in X^T$ with distinguished paths $y_U \sim x$ within $X^U$;
  in particular, by contracting paths, we may construct a deformation retract of $\prn{\prn{\Tot_I X^{I-\sqcup}}^{\act}_{/\gamma(x,S)}}^{\simeq}$ onto the subspace $\FF_{I, /S}^{\simeq} \subset \prn{\prn{\Tot_I X^{I-\sqcup}}^{\act}_{/\gamma(x,S)}}^{\simeq}$ of identity paths.
  
  Similarly, we may perform a deformation retract of $\prn{\prn{X \times \uFF_{i,*}}^{\act}_{/(x,S)}}^{\simeq}$ onto the summand $\FF_{I, /S}^{\simeq} \subset \prn{\prn{X \times \uFF_{i,*}}^{\act}_{/(x,S)}}^{\simeq}$ of identity paths.
  It follows by unwinding definitions that these are taken isomorphically onto each other;
  alternatively, one may note that the induced endomorphism of $\FF_{I,/S}^{\simeq} \simeq \prod_{U \in \Orb(S)} \coprod_{T \in \FF_{I, U}} B \Aut_U(T)$ is a product of coproducts of maps classified by torsor maps $\Aut_U(T) \rightarrow \Aut_U(T)$, which are automatically isomorphisms, implying that our map of $1$-truncated spaces is an isomorphism on $\pi_0$ and on $\pi_1$ at all basepoints.
\end{proof} 

\begin{warning}
  A closely related analog of \cref{gamma is approx} is claimed in \cite[Rmk~2.4.3.6]{HA} in the case $\cT = *$ \emph{without the assumption that $\cC$ is a space};
  as pointed out in \cite[Rmk~2.3]{Krannich} Lurie's claim (and hence proof) is incorrect in general, but the claim was verified in \emph{loc. cit.} when $\cC$ is a space.
\end{warning}

We finish with the following proposition.
\begin{proposition}\label{Family left adjoint}
  If $\Tot_{\cT} \cP^{\otimes} \rightarrow \Tot^{\Int}\prn{ X \times \uFF_{I,*}}$ is a fibrous pattern, then $L_{\Op_I} \cP^{\otimes}$ is the $\cT$-colimit of the $\cT$-functor $X \rightarrow \uOp_I$ associated with $\cP^{\otimes}$.
\end{proposition}
\begin{proof}
  Note that $L_{\Op_I} \pi_{X^{I-\sqcup}!} \gamma_! \dashv \gamma^* \pi_{X^{I-\sqcup}}^*$, and the latter is equivalent to $\Delta\colon \Op_I \rightarrow \Fun_{\cT}(X, \uOp_I)$;
  the above presentation for the left adjoint is $L_{\Op_I} \cP^{\otimes}$, and indexed colimits are also left adjoint to $\Delta$, so the claim follows from uniqueness of left adjoints.
\end{proof}

We now apply this in the language of \emph{disintegration and assembly}.
\subsection{Disintegration and assembly}
Given  $X \in \cS_{\cT}$ and $\cO^{\otimes} \in \Op_{I, /X^{I-\sqcup}}$,  define the pullback $\Gamma$-$I$-preoperad
\[\begin{tikzcd}[ampersand replacement=\&]
	{\dis^I\prn{\cO^{\otimes}}} \& {\Tot_{\cT} \cO^{\otimes}} \\
	{X \times \uFF_{I,*}} \& {\Tot_{I} X^{I-\sqcup}}
	\arrow["\alpha", from=1-1, to=1-2]
	\arrow[from=1-1, to=2-1]
	\arrow["\lrcorner"{anchor=center, pos=0.125}, draw=none, from=1-1, to=2-2]
	\arrow[from=1-2, to=2-2]
	\arrow["\gamma", from=2-1, to=2-2]
\end{tikzcd}\]
We refer to $\dis^I(-)$ as the \emph{disintegration functor} and $\alpha$ as the \emph{assembly map}.
\begin{proposition}\label{assembly localization}
  $\alpha$ is an $L_{\Op_I}$-localization map.
\end{proposition}
\begin{proof}
  We verified in \cite[\S~A]{EBV} that the conditions of \cref{Equivalence theorem} are pullback-stable, so $\alpha$ is a Morita equivalence;
  by \cref{Morita equivalence vs L equivalence} it is then an $L_{\Op_I}$-equivalence.
  By assumption, $\Tot_{\cT} \cO^{\otimes}$ is $L_{\Op_I}$-local, proving the proposition.
\end{proof}

\begin{proposition}\label{assembly colimit}
  $\cO^{\otimes}$ is the $\cT$-colimit of the $\cT$-functor $X \rightarrow \uOp_I$ associated with $\dis^I(\cO^{\otimes})$.
\end{proposition}
\begin{proof}
  This is a straightforward application of \cref{Family left adjoint,assembly localization}.
\end{proof}

We spell out the following corollary, which summarizes the full power of what we've proved.
\begin{corollary}
  Let $X$ be a $\cT$-space.
  The assignment $x \mapsto \cO_x \deq  \Res_V^{\cT} \cO^{\otimes} \times_{\Res_V^{\cT} X^{I-\sqcup}} \cN_{I\infty}^{I-\sqcup}$,
  yields an equivalence
  \[
    \uOp_{I, /X^{I-\sqcup}} \simeq \uFun_{\cT}(X, \uOp_I).
  \]
  The unit of this equivalence specifies a natural equivalence.
  \[
    \cO^{\otimes} \simeq \ucolim_{x \in X} \cO_x^{\otimes}.
  \]
\end{corollary}
\begin{proof}
  The first claim follows by combining \cref{T-space pattern,C-family pattern,Equivalence theorem,gamma is approx}.
  The remaining claim is proved identically to \cref{assembly colimit}.
\end{proof}
\resumetocwriting

\section{Algebraic patterns and the \tI-symmetric monoidal structure on overcategories}
\stoptocwriting
In this appendix, we repeat the arguments of \cite[\S~2.2.2]{HA} in the setting of algebraic patterns, assuming familiarity with the minutiae of \cite[\S~2.2.2]{HA} and of \cite{Barkan}.
\subsection{The fibrous pattern case}
We fix $\fO$ an algebraic pattern and make the following temporary definitions.
\begin{definition}
  An \emph{$\fO$-monoidal $\infty$-category} is a fibrous pattern $\fC \rightarrow \fO$ which is also a cocartesian fibration.
  If $\fC \rightarrow \fO$ is a fibrous pattern, then the $\infty$-category of \emph{$\fO$-algebras in $\fC$} is
  \[
    \Alg_{/\fO}(\fC) \deq \Fun_{/\fO}^{\Int-\cocart}(\fO, \fC).\qedhere
  \]
\end{definition}
\begin{remark}
  When $\fO$ is sound, 
  \cite[Prop~4.1.7]{Barkan} shows that $\fO$-monoidal $\infty$-categories are synonymous with Segal fibrations to $\fO$ in the sense of \cite{Chu}.
\end{remark}

\begin{warning}
  If $\cO^{\otimes}$ underlies a $\cT$-operad, $\Alg_{\cO}(\cC)$ and $\Alg_{/\cO}(\cC)$ should not be confused;
  the latter consists of algebras \emph{over $\cO^{\otimes}$}.
  Nevertheless, in the case that $\cO^{\otimes} \rightarrow \Comm_{\cT}^{\otimes}$ is a monomorphism (i.e. $\cO^{\otimes}$ is a weak $\cN_\infty$-operad), these agree.
\end{warning}

For the duration of this appendix, we fix $q\colon \fC \rightarrow \fO$ a fibrous pattern and $p\colon K \rightarrow \Alg_{/\fO}(\fC)$ a $K$-indexed diagram of $\fO$-algebras in $\fC$.
Let $q''\colon \fC_{/p_{\fO}} \longrightarrow \fO$ be the construction made in \cite[Def~2.2.2.1]{HA}, interpreted as a functor of $\infty$-categories via \cite[Lem~2.2.2.6]{HA}.
The proof of the following theorem will involve essentially no new ideas over that of \cite[Thm~2.2.2.4]{HA}.
\begin{theorem}\label{Slice pattern theorem}
  $q''$ exhibits $\fC_{/p_{\fO}}$ as a fibrous $\fO$ pattern, which is an $\fO$-monoidal $\infty$-category if $\fC$ is.
\end{theorem}
\begin{proof}
  We apply \cite[Lem~2.2.2.7-9]{HA} on opposite categories.
  In particular, given an arrow $g\colon X \rightarrow q''(Y)$ in $\fO$ over which $\fC$ has a cocartesian lift, \cite[Lem~2.2.2.8]{HA} supplies a lift
  \[\begin{tikzcd}[ampersand replacement=\&, column sep=large]
	{*} \& {\fC_{/p_{\fO}}} \\
	{*^{\triangleright}} \& \fO
	\arrow["{\cbr{Y}}", from=1-1, to=1-2]
	\arrow[hook, from=1-1, to=2-1]
	\arrow["{q''}", from=1-2, to=2-2]
	\arrow["{\overline{g}}"{description}, dashed, from=2-1, to=1-2]
	\arrow["g"', from=2-1, to=2-2]
  \end{tikzcd}\]
  such that $\overline{g}$ is a $q$-colimit diagram;
  moreover, \cite[Lem~2.2.2.7]{HA} guarantees that $\overline{g}$ is a $q''$-cocartesian lift of $g$.
  Since $\fC$ has inert-cocartesian lifts, so does $\fC_{/p_{\fO}}$, and when $q$ is a cocartesian fibration, so is $q''$.

  We're left with verifying the Segal condition(s) for fibrous patterns;
  we use that of \cite[Prop~4.1.6]{Barkan}.
  As in the proof of \cite[Thm~2.2.2.4]{HA}, it follows from a simple application of \cite[Lem~2.2.2.9]{HA} that each of the relevant diagrams are limit diagrams, as they project to limit diagrams in $\fC$.
\end{proof}

\begin{remark}\label{cocartesian slice remark}
  As in \cite[Thm~2.2.2.4.(2)]{HA}, it follows from the above diagram that, given an arrow $f$ in $\fC_{/p_{\fO}}$ such that $\fC$ admits a $q''(f)$-cocartesian arrow, $f$ is $q''$-cocartesian if and only if its image in $\fC$ is $q$-cocartesian.
  In particular the inert arrows in $\fC_{/p_{\fO}}$ are the preimages of the inert morphisms of $\fC$, and if $\fC$ is $\fO$-monoidal, then the cocartesian arrows in $\fC_{/p_{\fO}}$ are the preimages of cocartesian arrows in $\fC$.
\end{remark}

It is worthwhile to explicitly record following immediate corollary of \cref{cocartesian slice remark}, in part because it establishes the ``pointwise'' nature of the coherences for slice $\fO$-monoidal structure.
\begin{corollary}\label{Unslicing is symmetric monoidal}
  If $\fC$ is an $\fO$-monoidal $\infty$-category, then the unslicing functor $\fC_{/p_{\fO}} \rightarrow \fC$ is an $\fO$-monoidal functor, i.e. it is a functor of cocartesian fibrations over $\fO$.
\end{corollary}

As claimed in \cite{HA}, we may pass through \emph{opposite categories} to establish the following result about undercategories without additional argument, noting that the additional assumption comes from the asymmetric assumptions of \cite[Lem~2.2.2.7,8]{HA}.
Let $q'\colon \fC_{p\fO/} \rightarrow \fO$ be the construction of \cite[Def~2.2.2.1]{HA}.
\begin{theorem}
  $q'$ exhibits $\fC_{p\fO/}$ as a fibrous $\fO$-pattern;
  moreover, if $\fC$ is an $\fO$-monoidal $\infty$-category and the value functors $\fO \rightarrow \fC$ are all $\fO$-monoidal functors, then $q'$ exhibits $\fC_{p\fO/}$ as an $\fO$-monoidal $\infty$-category.
\end{theorem}

\subsection{The \tI-symmetric monoidal case}
We explicitly specialize \cref{Slice pattern theorem} to $\fO \deq \Span_I(\FF_{\cT})$. 
\begin{corollary}\label{I-symmetric monoidal slice corollary}
  Let $\cC^{\otimes}$ be an $I$-symmetric monoidal $\infty$-category and $A \in \CAlg_I(\cC)$ an $I$-commutative algebra in $\cC$.
  Then, \cref{Slice pattern theorem} supplies an $I$-symmetric monoidal $\infty$-category $\cC^{\otimes}_{/A}$ such that
  \begin{enumerate}
    \item The underlying $\cT$-$\infty$-category of $\cC^{\otimes}_{/A}$ agrees with Shah's slice $\cT$-$\infty$-category $\cC_{/(A,\cT^{\op})}$ \cite{Shah};
      moreover, if $\NN *_V \subset \FF_{I,V}$, then the induced symmetric monoidal structure on $\cC^{\otimes}_{/A, \Res_V^{\cT} A}$ agrees with Lurie's with respect to the restricted $\EE_\infty$-algebra $\Res_V^{\cT} \in \CAlg(\cC_V)$.
    \item The $S$-indexed tensor functor in $\cC^{\otimes}_{/A}$ takes a tuple of maps $\prn{f_U\colon X_U \rightarrow \Res_U^{\cT} A}_S$ to the map
      \[
        \bigotimes_U^S X_U \xrightarrow{\;\;\; \bigotimes_U^S f_U \;\;\;} \bigotimes_U^S \Res_U^{\cT} A \xrightarrow{\;\;\; \mu \;\;\;} \Res_V^{\cT} A.
      \]
  \end{enumerate}
\end{corollary}
\begin{proof}
  (1) is functoriality of the relative slice construction with respect to pullback of the base $\infty$-category;
  this follows straightforwardly from the defining universal property.

  For (2), we may apply the universal property of \cite[Def~2.2.2.1]{HA} along the functor $\Delta^1 \rightarrow \Span_I(\FF_{\cT})$ classifying an $I$-admissible active arrow $\psi\colon \Ind_V^{\cT} S = \Ind_V^{\cT} S \rightarrow V$:
 active arrows lying over $\psi$ are in correspondence with dashed arrows (and homotopies) making the following diagram commute.
  \[\begin{tikzcd}[ampersand replacement=\&]
	{\Delta^1} \& {\Delta^1 \times *^{\triangleleft}} \& {\Delta^1} \\
  {\Span_I(\FF_{\cT})} \& \cC^{\otimes} \& {\Span_I(\FF_{\cT})}
	\arrow[hook', from=1-1, to=1-2]
	\arrow["\psi", from=1-1, to=2-1]
	\arrow[two heads, from=1-2, to=1-3]
	\arrow[from=1-2, to=2-2, dashed]
	\arrow["\psi", from=1-3, to=2-3]
	\arrow["{\cbr{A}}"', from=2-1, to=2-2]
	\arrow["\pi"', from=2-2, to=2-3]
\end{tikzcd}\]
  That is, active arrows over $\psi$ correspond with commuting diagrams
  \[\begin{tikzcd}[ampersand replacement=\&]
	{(X_U)} \& Y \\
  {(\Res_U^{\cT} A)} \& \Res_V^{\cT} A
	\arrow[from=1-1, to=1-2]
	\arrow[from=1-1, to=2-1]
	\arrow[from=1-2, to=2-2]
	\arrow[from=2-1, to=2-2]
\end{tikzcd}\]  
(with homotopies witnessing commutativity, between the bottom arrow and the active arrow $(\Res_U^{\cT} A) \rightarrow \Res_V^{\cT} A$, and between the underlying arrows in $\Span_I(\FF_{\cT})$ and $\psi$).
To compute the $S$-indexed tensor functor, we are tasked with exhibiting \emph{cocartesian} active arrows, and by \cref{cocartesian slice remark} it suffices to construct an active arrow whose top arrow is cocartesian.
Indeed, the outer diagram of the following suffices.
\[\begin{tikzcd}[ampersand replacement=\&]
	{(X_U)} \& {\bigotimes\limits_U^S X_U} \& {\bigotimes\limits_U^S X_U} \\
	{(\Res_U^{\cT} A)} \& {\bigotimes\limits_U^S \Res_U^{\cT} A} \& {\Res_V^{\cT} A}
	\arrow[from=1-1, to=1-2]
	\arrow[from=1-1, to=2-1]
	\arrow[equals, from=1-2, to=1-3]
	\arrow[from=1-2, to=2-2]
	\arrow[from=1-3, to=2-3]
	\arrow[from=2-1, to=2-2]
	\arrow[from=2-2, to=2-3]
\end{tikzcd}\]
\end{proof}

\resumetocwriting
\section{The \tI-symmetric monoidal structure on \tI-operadic left Kan extension}
\stoptocwriting
In the first arXiv copy of this paper, we incorrectly claimed that operadic left Kan extension along an arbitrary $G$-operad map $\triv(\cA)^{\otimes} \rightarrow \triv(\cB)^{\otimes}$ yielded a $G$-symmetric monoidal structure on $G$-left Kan with respect to the pointwise structure.
This neglected a compatibility assumption between the target $G$-symmetric monodial category and the colimits appearing in the pointwise formula for left Kan extension, which is nevertheless satisfied in all examples of interest.
We rectify this now, asserting the following compatibility.
\begin{definition}
  Let $\cK$ be a collection of small $\cT$-categories.
  A $I$-symmetric monoidal $\infty$-category $\cC^{\otimes}$ is \emph{compatible with $\cK$-colimits} if, for all $I$-admissible $S$, the $S$-indexed tensor $V$-functor
  \[
    \bigotimes^S\colon \prod_{U \in \Orb(S)} \cC_{\uU} \rightarrow \cC_{\uV}
  \]
  preserves $\cK$-colimits.
\end{definition}
\begin{proposition}\label{prop:I-symmetric-monoidal-I-operadic-left-kan-extension}
  Let $J \subset I \subset \uFF_{I,*}$ be weakly extensive subcategories, suppose $\cC^{\otimes} \in \Cat_I^{\otimes}$ is compatible with $\cK$-colimits,
  and let $\varphi\colon \cO^{\otimes} \rightarrow \cP^{\otimes}$ be a $J$-operad map such that
  \begin{enumerate}[label={(\alph*)}]
    \item \label[condition]{cond:slice-K} the slice $V$-categories $\Env_J\prn{\cO}_{\uV, /P}$ lie in $\cK$ for all $P \in \cP_V$, and
    \item the pullback $\cT$-functor $\varphi^*\colon \uAlg_{\cP}(\cC) \rightarrow \uAlg_{\cO}(\cC)$ admits a $\cT$-left adjoint satisfying the Beck-Chevalley condition that the following diagram of $\cT$-functors commutes
      \[
        \begin{tikzcd}[ampersand replacement=\&,  column sep=huge, row sep=tiny]
          {\uAlg_{\cO}(\cC)} \& {\uAlg_{\cP}(\cC)} \\
	{\Fun^{\otimes}_J\prn{\Env_J(\cO)^{\otimes}, \Bor_J^I \cC^{\otimes}}} \& {\Fun^{\otimes}_J\prn{\Env_J(\cP)^{\otimes}, \Bor_J^I \cC^{\otimes}}} \\
	\\
	\\
	\\
	{\uFun_{\cT}(\Env_J(\cO),\cC)} \& {\uFun_{\cT}(\Env_J(\cP),\cC)}
	\arrow["{{\varphi_!}}", from=1-1, to=1-2]
	\arrow["\simeq"{marking, allow upside down}, draw=none, from=1-1, to=2-1]
	\arrow["\simeq"{marking, allow upside down}, draw=none, from=1-2, to=2-2]
	\arrow["{\varphi_!}", from=2-1, to=2-2]
	\arrow["{U_{\Env}}", from=2-1, to=6-1]
	\arrow["{U_{\Env}}", from=2-2, to=6-2]
	\arrow["{{\Env_I(\varphi)_!}}", from=6-1, to=6-2]
        \end{tikzcd}
      \]
  \end{enumerate}
  Then $\varphi^*\colon \uAlg_{\cP}^{\otimes}(\cC) \rightarrow \uAlg_{\cO}^{\otimes}(\cC)$ admits an $I$-symmetric left adjoint lifting $\varphi_!$, naturally in $\cK$-colimit preserving $I$-symmetric monoidal functors in $\cC$.
\end{proposition}
Before proving this, we note how it implies a desirable corollary.
\begin{corollary}\label{cor:T-sym-mon-kan}
  Suppose $\cT$ is atomic orbital, $J \subset I$ $\cT$-weak indexing categories, and $\cC^{\otimes} \in \Cat_{I}^{\otimes}$ distributive.
  \begin{enumerate}
    \item If $\varphi\colon \cO^{\otimes} \rightarrow \cP^{\otimes}$ is a $J$-operad map satisfying \cref{cond:slice-K} for $\cT$-sifted diagrams for some $J \subset I$, 
      then $\cT$-operadic left Kan extension lifts to an $I$-symmetric monoidal adjunction
      \[
        \begin{tikzcd}[ampersand replacement=\&]
          {\uAlg_{\cO}^{\otimes}(\cC)} \& {\uAlg_{\cP}^{\otimes}(\cC),}
          \arrow[""{name=0, anchor=center, inner sep=0}, "{\varphi_!}", curve={height=-6pt}, from=1-1, to=1-2]
          \arrow[""{name=1, anchor=center, inner sep=0}, "{\varphi^*}", curve={height=-6pt}, from=1-2, to=1-1]
          \arrow["\dashv"{anchor=center, rotate=-90}, draw=none, from=0, to=1]
        \end{tikzcd}
      \]
      natural in $\cT$-sifted $\cT$-colimit preserving $I$-symmetric monoidal functors in $\cC$.
    \item If $F\colon \cA \rightarrow \cB$ is a $\cT$-functor such that, for all $V \in \cT$ and $A \in \cA_{V}$, the slice $V$-category $\cB_{\uV,/A}$ is $\cT$-sifted, then the left Kan extension functor $F_!\colon \uFun_{\cT}(\cA,\cC) \rightarrow \uFun_{\cT}(\cB,\cC)$ underlies an $I$-symmetric monoidal left adjoint to $F^*$ under the pointwise $I$-symmetric monoidal structure, naturally on $\cT$-sifted colimit preserving $I$-symmetric monoidal functors in $\cC$.
\end{enumerate}
\end{corollary} 
\begin{proof}
  We begin by noting that statement (2) is a special case of statement (1) for the initial one-color weak indexing category $J = I^{\triv}$,
  as the $I^{\triv}$-symmetric monoidal envelope is simply the identity.

  Now, for statement (1), note 
  that the proof of \ns{Thm.}{5.1.4} shows that $\cC$ is compatible with $\cT$-sifted $\cT$-colimits.
  Additionally, \ns{Thm.}{4.3.4} constructs a $\cT$-left adjoint to $\varphi^*$.
  We're left with proving that the Beck-Chevalley transformation
  \[
    \Env_I(\varphi)_! U \implies U \varphi_!
  \]
  is an equivalence;
  this follows by unwinding \ns{Rmk.}{4.3.6}, or by mimicking the proof of \llp{Lem.}{3.40}.
\end{proof}

Now, we'll prove \cref{prop:I-symmetric-monoidal-I-operadic-left-kan-extension} in two steps:
first we construct an \emph{oplax} $I$-symmetric monoidal structure, then we prove that it is $I$-symmetric monoidal.
The oplax symmetric monoidal structure comes from general considerations around the doctrinal adjunction, so we briefly increase generality.
\subsection{A doctrinal adjunction for \texorpdfstring{$\fO$-monoidal $\infty$-categories}{O-monoidal infinity-categories}}
Fix $\fO$ an algebraic pattern.
If $\fC,\fD \rightarrow \fO$ are $\fO$-monoidal $\infty$-categories, we refer to morphisms of fibrous $\fO$-patterns $\fC \rightarrow \fD$ as \emph{lax $\fO$-monoidal functors}.
Given $F\colon \fC \rightarrow \fD$ a lax $\fO$-monoidal functor, $f\colon P \rightarrow O$ an active arrow in $\fO$, and $X \in \fC_P$ a $P$-object in $\fC$, we refer to the dashed arrow 
\[
  \begin{tikzcd}[ampersand replacement=\&, column sep=huge]
	FX \\
	{f_{\otimes} F_X} \& {F f_{\otimes} X}
	\arrow["{c_{f,FX}}"', from=1-1, to=2-1]
  \arrow["{Fc_{f,X}}", from=1-1, to=2-2, curve={height=-6pt}]
  \arrow[dashed, "\ell_{F,f,X}", from=2-1, to=2-2]
\end{tikzcd}
\]
supplied by the universal property for cocartesian arrows as the \emph{$f$-indexed lax structure map of $F$ at $X$};
here, $f_{\otimes}$ is cocartesian transport along $f$ and $c_{f,Y}$ is the cocartesian arrow lifting $f$ with domain $Y$.
A lax $\fO$-monoidal functor $F$ is $\fO$-monoidal if and only if $\ell_{F,f,X}$ is an equivalence for all $f,X$.

Now, postcomposition with the involution $(-)^{\op}\colon \Cat \rightarrow \Cat$ yields an involution $(-)^{\vop}\colon \Cat_{\fO}^{\otimes} \rightarrow \Cat_{\fO}^{\otimes}$ called the \emph{vertical opposite}.
We will refer to lax $\fO$-monoidal functors $\fC^{\vop} \rightarrow \fD^{\vop}$ as \emph{oplax $\fO$-monoidal functors from $\fC$ to $\fD$};
an oplax $\fO$-monoidal functor $F\colon \fC \rightarrow \fD$ is $\fO$-monoidal if and only if the $f$-indexed oplax structure map of $F$ at $X$, $o_{F,f,X}\colon F f_{\otimes} X \longrightarrow f_{\otimes} F X$,
is an equivalence for all $f,X$.

We refer to adjunctions of $\fO$-monodial $\infty$-categories $F\colon \fC \rightleftarrows \fD\colon R$ relative to $\fO$ as \emph{lax $\fO$-monoiadal} (resp. oplax $\fO$-monoidal, $\fO$-monoidal) if $F$ and $R$ are lax $\fO$-monoidal (oplax $\fO$-monoidal, $\fO$-monodial).
\begin{proposition}[Doctrinal adjunction for algebraic patterns]\label{prop:doctrinal}
  Suppose $F^{\otimes}\colon \cD^{\otimes} \rightarrow \cC^{\otimes}$ is an $\fO$-monoidal functor whose underlying $\fO^{\el,\op}$-functor $F$ map admits an $\fO^{\el,\op}$-right adjoint $R$.
  Then, there is a lax $\fO$-monoidal right adjoint $R^{\otimes} \vdash F^{\otimes}$ lifting $R$, such that
  \begin{enumerate}[label={(\alph*)}]
    \item For all $O \in \fO$, the following commutes
      \[
        \begin{tikzcd}[ampersand replacement=\&,  column sep=huge, row sep=tiny]
          {\cC_O} \& {\cD_O} \\
          {\lim_{E \in \fO^{\el}_{/O}} \cC_E} \& {\lim_{E \in \fO^{\el}_{/O}} \cD_E}
          \arrow["{R^{\otimes}_O}", from=1-1, to=1-2]
          \arrow["\simeq"{marking, allow upside down}, draw=none, from=1-1, to=2-1]
          \arrow["\simeq"{marking, allow upside down}, draw=none, from=1-2, to=2-2]
          \arrow["{\lim_E R_E}"', from=2-1, to=2-2]
        \end{tikzcd}
    \]
  \item Given $f\colon P \rightarrow O$ in $\fO^{\act}$, the $f$-indexed lax structure map of $R$ at $X \in \cC_P$
    is an adjunct arrow
    \[
      \begin{tikzcd}[ampersand replacement=\&, column sep=small]
	\& RX \\
	{f_{\otimes} RX} \& {} \&\&\&\& {R f_{\otimes} X} \\
	{F f_{\otimes} RX} \& {f_{\otimes} FRX} \&\&\&\& {f_{\otimes} X}
	\arrow["c", from=1-2, to=2-1]
	\arrow["{{{g c}}}", from=1-2, to=2-6]
  \arrow[dashed, from=2-1, to=2-6, "\ell_{R,f,X}"]
	\arrow[""{name=0, anchor=center, inner sep=0}, draw=none, from=2-2, to=2-6]
	\arrow[""{name=0p, anchor=center, inner sep=0}, phantom, from=2-2, to=2-6, start anchor=center, end anchor=center]
	\arrow["\simeq"{marking, allow upside down}, draw=none, from=3-1, to=3-2]
	\arrow[""{name=1, anchor=center, inner sep=0}, "{f_{\otimes} \epsilon_{F,X}}", from=3-2, to=3-6]
	\arrow[""{name=1p, anchor=center, inner sep=0}, phantom, from=3-2, to=3-6, start anchor=center, end anchor=center]
	\arrow[between={0.3}{0.6}, squiggly, tail reversed, from=0p, to=1p]
\end{tikzcd}
\]
where $\epsilon_{F,X}$ is the counit of the adjunction $R_{P} \vdash F_P$ applied at $X$.
  \end{enumerate}
  If alternatively, $F$ admits an $\fO^{\el,\op}$-left adjoint $L$, then there is a oplax $\fO$-monoidal left adjoint $L^{\otimes} \dashv f^{\otimes}$ lifting $L$ which is vertical-opposite to the lax $\fO$-monoidal left adjoint of $f^{\otimes,\vop}$;
  in particular,
  \begin{enumerate}[label={(co-\alph*)}]
    \item for all $O \in \fO$, the following commutes
      \[
        \begin{tikzcd}[ampersand replacement=\&,  column sep=huge, row sep=tiny]
          {\cC_O} \& {\cD_O} \\
          {\lim_{E \in \fO^{\el}_{/O}} \cC_E} \& {\lim_{E \in \fO^{\el}_{/O}} \cD_E}
          \arrow["{L^{\otimes}_O}", from=1-1, to=1-2]
          \arrow["\simeq"{marking, allow upside down}, draw=none, from=1-1, to=2-1]
          \arrow["\simeq"{marking, allow upside down}, draw=none, from=1-2, to=2-2]
          \arrow["{\lim_E L_E}"', from=2-1, to=2-2]
        \end{tikzcd}
    \]
  \item The oplax structure map
    $o_{L,f,X}\colon Lf_{\otimes} X \rightarrow f_{\otimes} LX$
    is adjunct to the unit
    \[
      f_{\otimes} X \xrightarrow{\;\;\; f_{\otimes} \eta_{F,X} \;\;\;} f_{\otimes} FL X \simeq F f_{\otimes} LX.
    \]
\end{enumerate}
\end{proposition}
\begin{proof}
  First note that the second half of this proposition follows from the first half by taking vertical opposites.
  For the first half, we perform a similar argument to \ha{Cor.}{7.3.2.7};
  the Segal condition implies that the arrow $R_O$ constructed above is a right adjoint to $F_O$ (e.g. with unit given by $\lim_E \eta_E$) and these are compatible with inert-cocartesian transport, so \ha{Prop.}{7.3.2.6} constructs a relative right adjoint $R^{\otimes} \dashv L^{\otimes}$ over $\fO$ which preserves inert-cocartesian arrows and satisfies condition (a);
  in particular, we've constructed a lax $\fO$-monoidal adjunction.

  Now, to prove part (b), it suffices to note that the specified adjunct makes the triangle diagram commute;
  the universal property for cocartesian arrows then constructs a homotopy between $\ell_{R,f,X}$ and the adjunct to $f_{\otimes},\epsilon_{F,X}$.
\end{proof}

Before moving on, we point out a quick corollary for adjunctions between $\cT$-categories of algebras.
\begin{corollary}
  Suppose $F^{\otimes}\colon \cC^{\otimes} \rightarrow \cD^{\otimes}$ is an $I$-symmetric monoidal functor whose underlying $\cT$-functor admits a right adjoint $R\colon \cD \rightarrow \cC$.
  Then, there is a canonical lax $I$-symmetric monoidal lift $R^{\otimes}$ for $R$, inducing a natural lax $I$-symmetric monoidal adjunction
  \[
    F_*^{\otimes}\colon \uAlg^{\otimes}_{\cO}(\cC) \longrightleftarrows \uAlg^{\otimes}_{\cO}(\cD)\colon R_*^{\otimes}
  \]  
  by postcomposition for any $\cO^{\otimes}$, so that $F_*^{\otimes}$ is $I$-symmetric monoidal.
  If $R$ is $I$-symmetric monoidal, then so is $R_*^{\otimes}$;
  if $R$ is also fully faithful, then so is $R_*^{\otimes}$.
\end{corollary}
\begin{proof}
  The lax $I$-symmetric monoidal right adjoint is supplied by \cref{prop:doctrinal}, and the (lax) $I$-symmetric functors on $\cO$-algebras are supplied by functoriality of $\uAlg_{\cO}(-)$ developed in \EBVsubsec{3.2}.

  Now, the counit for $F^{\otimes} \dashv R^{\otimes}$ is a lax $I$-symmetric monoidal natural transformation $\epsilon\colon F^{\otimes} R^{\otimes}\implies \id_{\cD^{\otimes}}$, i.e. a functor $\Delta^1 \rightarrow \Alg_{\cD}(\cD)$.
  A similar construction holds for $\eta$.

  The construction of lifts to algebras yields a lift $\epsilon^{\otimes}\colon \Delta^1 \rightarrow \Alg_{\uAlg_{\cO}^{\otimes}(\cC)}(\uAlg_{\cO}^{\otimes}(\cC))$ for $\epsilon$ and similarly a lift $\eta^{\otimes}$ for $\eta$;
  functoriality of this construction shows the triangle identities for $\eta^{\otimes}$ and $\epsilon^{\otimes}$, yielding the desired lax $I$-symmetric monoidal adjunction.
  Symmetric monoidality of $R_*^{\otimes}$ is proved in \EBVsubsec{3.2};
  fully faithfulness follows by the counit expression for fully faithful right $\cT$-adjoints.
\end{proof}
For instance, this applies to easily construct a lax $I$-symmetric monoidal right adjoint to 
\[
  \SS[-] = \Sigma^\infty_*\colon \uAlg_{\cO}^{\otimes} \prn{\ucS_G} \rightarrow \uAlg_{\cO}^{\otimes}\prn{\uSp_G}.
\]

\subsection{Assembly for \tI-operadic left Kan extension}
We're ready to conclude.
\begin{proof}[Proof of \cref{prop:I-symmetric-monoidal-I-operadic-left-kan-extension}]
  Let $\varphi^{\otimes}_!$ be the oplax $I$-symmetric monoidal left adjoint supplied by 
  \cref{prop:doctrinal}.
  Now, by the Beck-Chevalley condition, the underlying natural transformation of the oplax structure map, $U_{\Env}o_{\varphi_!,\Ind_{V}^{\cT} S \rightarrow V,(X_U)_{U \in \Orb(S)}}(P)$, is given by the map
    \[
      \ucolim_{O \in \Env_J(\cO)_{\uV, /P}}\bigotimes_U^S X_U(O) \longrightarrow \bigotimes_U^S \ucolim_{O \in \Env_J(\cO)_{\uV, /P}}X_U(O)
    \]
    adjunct the unit of the adjunction $\Env_J(\varphi)_! \dashv \Env_J(\varphi)^*$;
    unwinding definitions, it is the assembly map for the indicated $\cT$-colimits in $\cC$ along the functor $\otimes^S$.
    This is an equivalence by the assumption that $\cC^{\otimes}$ is compatible with $\Env_J(\cO)_{\uV, /P}$-indexed colimits, so conservativity of $U$ implies that $o_{\varphi_!,\Ind_V^{\cT} S \rightarrow V, (X_U)}$ is a natural equivalence for all $S,(X_U)$.
    This shows that $\varphi_!^{\otimes}$ is $I$-symmetric monoidal, and the result then follows from 
  \cref{prop:doctrinal}.
\end{proof}

\resumetocwriting
\end{appendix}

\printbibliography

\end{document}